\documentclass [11pt]{amsart}
\usepackage {amsmath, amssymb, amscd, graphicx, color, url,epsf}
\usepackage[all,knot,arc]{xy}
\usepackage[text={6.5in,9in},centering,letterpaper,dvips]{geometry}
\usepackage{hyphenat}

\newtheorem {theorem}{Theorem}[section]
\newtheorem {lemma}[theorem]{Lemma}
\newtheorem {prop}[theorem]{Proposition}

\theoremstyle{definition}
\newtheorem {definition}[theorem]{Definition}

\theoremstyle{remark}
\newtheorem {remark}[theorem]{Remark}

\def\Cayley{\Gamma}
\def\straight{e}

\renewcommand\th{^{\text{th}}}
\newcommand\st{^{\text{st}}}
\newcommand\Z{\mathbb{Z}}
\newcommand\R{\mathbb{R}}
\newcommand\Q{\mathbb{Q}}

\newcommand\Xs{\mathbb{X}}
\newcommand\Os{\mathbb{O}}
\newcommand\Pent{\mathrm{Pent}}
\newcommand\Hex{\mathrm{Hex}}
\newcommand\Rect{\mathrm{Rect}}
\newcommand\Ring{R}

\DeclareMathOperator\rank{rank}
\def\halfz{\frac{1}{2}\Z}
\def\halfzl{(\halfz)^{\ell}}

\def\Interior{\mathrm{Int}}
\def\Ca{\widehat{C}}
\def\CLm{\textit{CL}^-}
\def\HLm{\textit{HL}^-}
\def\CLa{\widehat{\textit{CL}}}
\def\HLa{\widehat{\textit{HL}}}
\def\CLaa{\widetilde{\textit{CL}}}
\def\HLaa{\widetilde{\textit{HL}}}
\def\Caa{\widetilde{C}}
\def\Cm{C^-}
\def\CFm{\textit{CF}^-}

\def\HFLa{{\widehat{\textit{HFL}}}}
\def\HFLm{\textit{HFL}^-}

\def\Os{\mathbb O}

\def\fin\qedhere
\def\pr {{\text{pr}}}

\def\cm{\cdot}
\def\s{\mathbf s}

\def\x{\mathbf x}
\def\y{\mathbf y}
\def\z{\mathbf z}
\def\Torus{\mathcal T}
\def\Down {\mathcal D}

\newcommand\orL{\vec{L}}

\def\S{\mathbf S}

\def\NE{\mathcal I}
\def\NESW{\mathcal J}
\def\Edges{\mathrm{Edges}}
\def\B {\mathcal B}

\def\gr{\mathrm{gr}}

\def\Field{\mathbb F_2}
\newcommand\Filt{\mathcal F}

\newcommand\dm{\partial^-}
\newcommand\daa{\tilde{\partial}}

\newcommand\EmptyRect{\Rect^o}
\newcommand\tRect{\mathrm{tRect}}
\newcommand\stRect{\mathrm{tRect}^*}

\newcommand\Square{\Sigma}

\newcommand\FL{F^L}
\newcommand\FR{F^R}

\newcommand\FaQ{\widetilde F_Q}

\newcommand\FLaQ{\widetilde F^L_Q}
\newcommand\FRaQ{\widetilde F^R_Q}

\newcommand\Left{\mathcal L}
\newcommand\Right{\mathcal R}

\def\CI{C^I}
\def\CNI{C^\sNI}
\def\CNN{C^\sNN}
\def\CQ{\Caa_Q}
\def\CIQ{\Caa^I_Q}
\def\CNIQ{\Caa^\sNI_Q}
\def\CNNQ{\Caa^\sNN_Q}
\def\pL{\pi^L}
\def\pR{\pi^R}
\def\pF{\pi^F}
\def\oL{\vec{L}}

\def\Id{\mathbb I}

\def\Interval{\mathbf I}
\def\NN{\mathbf{NN}}
\def\NI{\mathbf{NI}}
\def\sNN{{\text{\sl NN}}}
\def\sNI{{\text{\sl NI}}}
\def\EmptyRect{\Rect^\circ}
\def\EmptyPent{\Pent^\circ}
\def\EmptyHex{\Hex^\circ}

\newcommand{\sign}{\mathcal S}

\def\mathcenter#1{\vcenter{\hbox{$#1$}}}
\def\mfig#1{\mathcenter{\includegraphics{#1}}}
\def\mfigb#1{\mathcenter{\includegraphics[trim=-1 -1 -1 -1]{#1}}}

{\makeatletter\gdef\reallynopagebreak{\nopagebreak\@nobreaktrue}}
\newcommand{\step}[2]{\subsection*{Step #1: #2}}

\begin{document}

\title{On combinatorial link Floer homology}

\author[Ciprian Manolescu]{Ciprian Manolescu}
\thanks {CM was supported by a Clay Research Fellowship.}
\address {Department of Mathematics, Columbia University\\ New York, NY 10027}
\email {cm@math.columbia.edu}

\author[Peter S. Ozsv\'ath]{Peter Ozsv\'ath}
\thanks {PSO was supported by NSF grant number DMS-0505811 and FRG-0244663}
\address {Department of Mathematics, Columbia University\\ New York, NY 10027}
\email {petero@math.columbia.edu}

\author[Zolt{\'a}n Szab{\'o}]{Zolt{\'a}n Szab{\'o}}
\thanks{ZSz was supported by NSF grant number DMS-0406155 and FRG-0244663}
\address{Department of Mathematics, Princeton University\\ Princeton, New Jersey 08544}
\email {szabo@math.princeton.edu}

\author[Dylan P. Thurston]{Dylan Thurston}
\thanks {DPT was supported by a Sloan Research Fellowship.}
\address {Department of Mathematics, Barnard College, Columbia University\\ New York, NY 10027}
\email {dthurston@barnard.edu}

\begin {abstract} 
  Link Floer homology is an invariant for links defined using a
  suitable version of Lagrangian Floer homology. In an earlier paper,
  this invariant was given a combinatorial description with mod $2$
  coefficients. In the present paper, we give a self-contained
  presentation of the basic properties of link Floer homology,
  including an elementary proof of its invariance. We also fix signs
  for the differentials, so that the theory is defined with integer
  coefficients.
\end {abstract}


\maketitle
\section{Introduction}

Heegaard Floer homology~\cite{HolDisk} is an invariant for
three-manifolds, defined using holomorphic disks and Heegaard
diagrams.  In~\cite{Knots} and~\cite{RasmussenThesis}, this
construction is extended to give an invariant, {\em knot Floer
homology}, for null-homologous knots in a closed, oriented
three-manifold. This construction is further generalized
in~\cite{Links} to the case of oriented links.  The definition of all
these invariants involves counts of holomorphic disks in the symmetric
product of a Riemann surface, which makes them rather challenging to
calculate.

More recently, Sucharit Sarkar discovered a principle which ensures
that for Heegaard diagrams with a certain property, the counts of
holomorphic disks are combinatorial. In~\cite{MOS}, the Heegaard
diagrams of the needed form are constructed from grid presentations of
knots or links in $S^3$, leading to an explicit, combinatorial
description of the knot or link Floer complex, taken with coefficients
in $\Z/2\Z$, henceforth called~$\Field$. (See also~\cite{SarkarWang}
for a different application of this principle.)
 
The purpose of the present paper is to develop knot (or link) Floer
homology in purely elementary terms, starting from a grid
presentation, and establish its topological invariance without
appealing to the earlier theory. We also give a
sign-refinement of this description, leading to a homology theory with
coefficients in $\Z$.

We recall the chain complex from~\cite{MOS}; but first, we need
to review some topological notions.

A {\em planar grid diagram}~$G$ lies on an $n\times n$ grid of squares
in the plane. Each square is decorated
either with an $X$, an $O$, or nothing. Moreover, the decorations
are arranged so that:
\begin{itemize}
\item {every row contains exactly one $X$ and one $O$;}
\item {every column contains exactly one $X$ and one $O$.}
\end {itemize}
The number $n$ is called the {\em grid number} of
$G$. Sometimes we find it convenient to number the
$O$'s and $X$'s by $\{O_i\}_{i=1}^n$ and $\{X_i\}_{i=1}^n$. We denote
the set of all $O$'s and $X$'s by~$\Os$ and~$\Xs$, respectively.
As a point of
comparison: the $O_i$ correspond to the ``white dots'' of~\cite{MOS}
and the $w_i$ of~\cite{Links}, while the $X_i$ to the ``black dots''
of~\cite{MOS} and the $z_i$ of~\cite{Links}. We find the current
notation clearer for pictures.

\begin{figure}
\[
\mfig{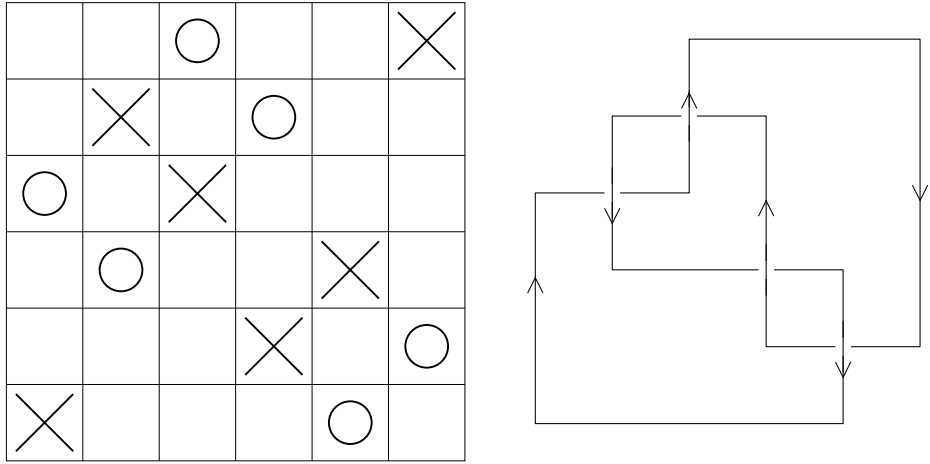}\qquad
\mfig{draws/gridlink.0}
\]
\caption {{\bf A grid presentation.}
Grid presentation for the figure eight knot.}
\label{fig:FigureEight}
\end{figure}

Given a planar grid diagram~$G$, we can place it in a
standard position on the plane as follows: the bottom left corner is
at the origin, and each cell is a square of edge length one. We then
construct an oriented, planar link projection by drawing horizontal
segments from the $O$'s to the $X$'s in each row, and vertical
segments from the $X$'s to the $O$'s in each column. At every
intersection point, we let the horizontal segment be the underpass and
the vertical one the overpass. This produces a planar diagram for an
oriented link $\orL$ in $S^3$. We say that $\orL$ has a grid
presentation given by~$G$. See Figure~\ref{fig:FigureEight}
for an example.

We transfer our grid diagrams to the torus $\Torus$ obtained by gluing
the topmost segment to the bottom-most one, and the leftmost segment
to the rightmost one. In the torus, our horizontal and vertical arcs
become horizontal and vertical circles.  The torus inherits its
orientation from the plane.  We call the resulting object
a {\em toroidal grid diagram}, or simply a grid diagram, for~$\orL$.  We
will again denote it by~$G$.

Given a toroidal grid diagram, we associate to it a chain complex
$\bigl(\Cm(G), \dm \bigr)$ as follows. The set of generators
of~$\Cm(G)$, denoted $\S$ or~$\S(G)$, consists of one-to-one
correspondences between
the horizontal and vertical circles. More geometrically, we can think
of the generators as $n$-tuples of intersection points between the
horizontal and vertical circles, with the property that no
intersection point appears on more than one horizontal (or vertical)
circle.

Before defining the differentials, we turn to a grading and a filtration 
on the complex, determined by two functions
$M \colon \S \longrightarrow \Z$ and
$A \colon \S \longrightarrow \halfzl$.

The function $M$ is defined as follows. Given two collections $A$,
$B$ of finitely many points in the plane, let $\NE(A,B)$ be the
number of pairs $(a_1,a_2)\in A$ and $(b_1,b_2)\in B$ with $a_1<b_1$
and $a_2<b_2$. Let $\NESW(A,B) = (\NE(A,B) + \NE(B,A))/2$.
Take a fundamental domain $[0,n)\times [0,n)$ for the torus, cut along a
horizontal and vertical circle, with the left and bottom edges
included.
Given a generator $\x\in\S$, we view $\x$ as a collection of points
with integer coordinates
in this fundamental domain.
Similarly, we view $\Os=\{O_i\}_{i=1}^n$ as a collection of
points in the plane with half-integer coordinates.
Define 
$$M(\x)=\NESW(\x,\x)-2\NESW(\x,\Os)
+\NESW(\Os,\Os)+1.$$
We find it convenient to write this formula more succinctly as
\begin{equation}
  \label{eq:MaslovFormula}
  M(\x)=\NESW(\x-\Os,\x-\Os)+1,
\end{equation}
where we extend $\NESW$  bilinearly over formal sums (or differences)
of subsets.
Note that the definition
of~$M$ appears to depend on which circles we cut along to create a
fundamental domain.  In fact, it does not (see
Lemma~\ref{lemma:MaslovWellDefined} below). Note also that this definition
of the Maslov grading is not identical with that given in~\cite{MOS}, but 
it is not difficult to see they agree.  See Lemma~\ref{lemma:MaslovChange}
below, and the remarks following it.

For an $\ell$-component link, we define an $\ell$-tuple of Alexander gradings
$A(\x)=(A_1(\x),\ldots,A_\ell(\x))$ by the formula
\begin{equation}
  \label{eq:AlexanderFormulaTwo}
A_i(\x) = \NESW(\x-\frac{1}{2}(\Xs+\Os),\Xs_i-\Os_i) - \Bigl(\frac{n_i-1}{2}\Bigr),
\end{equation}
where here $\Os_i\subset \Os$ is the subset corresponding to the
$i\th$ component of the link, $\Xs_i\subset \Xs$ is the set of $X$'s
belonging to the $i\th$ component of the link, and where we once again
use the bilinear extension of $\NESW$. For links, the $A_i$ may take
half-integral values. Again, this quantity is independent of how the
torus is cut up to form a planar rectangle (see
Lemma~\ref{lemma:AlexanderWellDefined} below).

\begin{figure}
\mbox{\vbox{\epsfbox{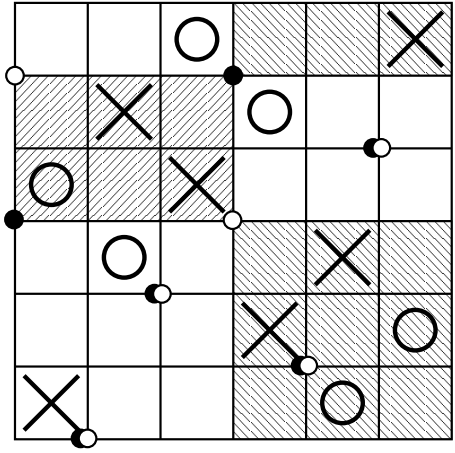}}}
\caption {{\bf Rectangles.}
The small dark circles describe the generator $\x$ and the hollow ones 
describe $\y$. There are two rectangles in $\Rect(\x, \y)$, shown here shaded 
by two types of diagonal hatchings. The rectangle on the left is in
$\EmptyRect(\x,\y)$ while the other one is not, because it 
contains a dark circle in its interior.}
\label{fig:Differential}
\end{figure}

Given a pair of generators $\x$ and $\y$, and an embedded rectangle $r$ in 
$\Torus$ whose edges are arcs in the horizontal and vertical circles, we 
say that $r$ connects $\x$ to $\y$ if $\x$ and $\y$ agree along all but 
two horizontal circles, if all four corners of $r$ are intersection points 
in $\x\cup\y$, and if, as we traverse each horizontal boundary 
component of $r$ in the direction dictated by the orientation that $r$ 
inherits from $\Torus$, then the arc is oriented from a point 
in $\x$ to the point in $\y$. (See Figure~\ref{fig:Differential} 
for an example.) Let $\Rect(\x,\y)$ denote the collection of 
rectangles connecting $\x$ to $\y$.
If $\x,\y\in \S$ agree
along all but two horizontal circles, then there are exactly two
rectangles in $\Rect(\x,\y)$; otherwise $\Rect(\x,\y)=\emptyset$.
Let $\Interior(r)$
denote the interior of the subset of $\Torus$ determined by $r$.
A rectangle~$r\in\Rect(\x,\y)$ is said to be \emph{empty} if
$\Interior(r)\cap\x = \emptyset$, or equivalently if
$\Interior(r)\cap\y = \emptyset$.  The space of empty rectangles
connecting~$\x$ and~$\y$ is denoted $\EmptyRect(\x,\y)$.

Let $\Ring$ denote the polynomial algebra over $\Field$ generated by
variables which are in one-to-one correspondence between the elements of 
$\Os$, and which we denote $\{U_i\}_{i=1}^n$. We think of
this ring as endowed with a {\em Maslov grading}, defined so that the
constant terms are in Maslov grading zero, and $U_i$ are in grading
$-2$.  The ring is also endowed with an {\em Alexander multi-filtration},
defined so that constant terms are in filtration level zero, while
the variables
$U_j$ corresponding to the $i\th$ component of the link
drop the $i\th$ multi-filtration level by one and preserve all others.

Let $\Cm(G)$ be the free $\Ring$-module with generating set $\S$.

We endow this module with an endomorphism $\dm\colon
\Cm(G)\longrightarrow \Cm(G)$ defined by
\begin{equation}
\label{eq:DefDm}
\dm(\x)=\sum_{\y\in\S}\,
\sum_{r\in\EmptyRect(\x,\y)}\!\!
U_1^{O_1(r)}\cdots U_n^{O_n(r)}\cm \y,
\end{equation} where $O_i(r)$ denotes the number of times $O_i$
appears in the interior of~$r$ (so $O_i(r)$ is either 0 or~1).

The results of~\cite{MOS} can be summarized by the following:

\begin{theorem} (Manolescu-Ozsv{\'a}th-Sarkar)
  \label{thm:MOS} The data $(\Cm(G),\dm)$ is a chain complex for the
  Heegaard-Floer homology
  $\CFm(S^3)$, with grading induced by $M$, and the filtration induced
  by $A$ coincides with the link filtration of
  $\CFm(S^3)$.\footnote{The reader should be warned: our conventions
  here on the Maslov grading are such that the total homology
  $H_*(\CFm(S^3))$ is isomorphic to a copy of the polynomial algebra
  in $U$, where the constants have grading equal to
  zero. In~\cite{HolDisk}, the convention is that the constants have
  grading equal to $-2$.}
\end{theorem}

In particular, appealing to the earlier theorem defined using
holomorphic disks~\cite{Links, Knots, RasmussenThesis}, the filtered
quasi-isomorphism type of this
chain complex $\Cm$ is a link invariant.  Other knot and link
invariants can be found by routine algebraic manipulations of $\Cm$ as
well (for example, by taking the homology of the associated graded
object).

Our main goal here is to prove the topological invariance of the filtered
quasi-isomorphism type of the resulting chain complex $\Cm(G)$,
without resorting to any of the holomorphic disk theory, and in
particular without resorting to Theorem~\ref{thm:MOS}. We prove the
following:

\begin{theorem}
\label{thm:WithSigns}
Let $\orL$ be an oriented, $\ell$-component link.  Number the elements
of $\Os=\{O_i\}_{i=1}^n$ so that $O_1,\ldots,O_\ell$ correspond to
different components of the link. Then the filtered quasi-isomorphism
type of the complex $(\Cm(G),\dm)$ over $\Z[U_1,\ldots,U_\ell]$ is an
invariant of the link.
\end{theorem}

We also give independent verification of the basic algebraic
properties of $\Cm(G)$ which, with $\Field$ (i.e., $\Z/2\Z$)
coefficients, follow
from Theorem~\ref{thm:MOS}, together with properties of the ``Heegaard
Floer homology package''.  Note that for technical reasons, for links
with more than one component the chain complex in~\cite{Links} was
originally defined only with coefficients in $\Field$.  

There are some related constructions one could consider. In one of
these, we set $U_1=\cdots=U_\ell=0$, let $\Ca(G)$ denote the resulting
chain complex, equipped with its Alexander filtration.  Taking the
homology of the associated graded object, we get a group whose
multi-graded Euler characteristic is the multi-variable Alexander
polynomial of $\orL$, times a suitable normalization factor (this is
proved in Equation~\eqref{Links:eq:EulerHFLa} of~\cite{Links}, see
also Theorem~\ref{thm:alexander} below).

We have endeavoured to separate the discussion of signs from the rest
of the body of the paper, to underscore the simplicity of the $\Field$
version which is sufficient for the knot-theoretic applications, and
also simpler to calculate. In particular, in
Section~\ref{sec:Invariance}, we establish
Theorem~\ref{thm:WithSigns}, working over coefficients in $\Field$,
where it could alternately be seen as an immediate consequence of
Theorem~\ref{thm:MOS}.  We hope, however, that the present
combinatorial proof of invariance has value in its own simplicity;
see also~\cite{Transverse} for another application.
The sign-refinements are dealt with in Section~\ref{sec:Signs}.

This paper is organized as follows. 
The algebraic
properties are established in Section~\ref{sec:BasicProperties}, and
topological invariance with coefficients in $\Field$
is established in Section~\ref{sec:Invariance}. In Section~\ref{sec:MoreProperties},
we describe some further properties of $\Cm$. In Section~\ref{sec:Signs},
we describe the sign conventions, and the modifications needed for the earlier
discussion to establish Theorem~\ref{thm:WithSigns} over $\Z$.
Finally, in Section~\ref{sec:Alexander}, we show that the Euler
characteristic of the homology is the Alexander polynomial.

\subsection*{Acknowledgements}

We would like to thank Dror Bar-Natan, Sergei Duzhin, Sergey Fomin, John Morgan, and
Sucharit Sarkar for helpful conversations.


\section{Properties of the chain complex $\Cm(G)$}
\label{sec:BasicProperties}

\subsection{Algebraic terminology} 

We recall some standard terminology from homological algebra.

For simplicity, we use coefficients in $\Field=\Z/2\Z$ for this
section, and also the next two. The definitions from algebra can be
made with $\Z$ coefficients with little change. Other
aspects of $\Z$ coefficients will be handled in
Section~\ref{sec:Signs}. (And in fact, the choices of signs in the formulas
below which, of course, are immaterial over $\Field$, have been 
chosen so as to work over $\Z$.)

\begin{definition}
  \label{def:Filtrations}
  We give $\Q^{\ell}$ its usual partial ordering,
  $(a_1,\ldots,a_\ell)\leq (b_1,\ldots,b_\ell)$ if for all
  $i=1,\ldots,\ell$, $a_i\leq b_i$.  Let $\Ring$ be the ring
  $\Field[U_1,\ldots,U_n]$. A function $g\colon
  \{1,\ldots,n\}\longrightarrow (\Q^{\geq 0})^\ell$ specifies a
  $\Q^{\ell}$ grading on $\Ring$. Fix a grading on $\Ring$.  Let $M$
  be a module over $\Ring$.  A {\em $\Q^\ell$-filtration} on a module
  $M$ is a collection of $\Ring$-submodules
  $\{\Filt_s(M)\}_{s\in\Q^\ell}$ of $M$ satisfying the following
  properties:
  \begin{itemize}
    \item $\Filt_s(M)\subset \Filt_t(M)$ if $s\leq t$
    \item multiplication by $U_i$ sends $\Filt_s(M)$ into 
      $\Filt_{s-g(i)}(M)$.
    \item  for
      all sufficiently large $s$ (with respect to $\leq$),
      $\Filt_s(M)=M$.
    \end{itemize} 
    A \emph{filtered $\Ring$-module map} $\phi\colon M \longrightarrow
    N$ is an
    $\Ring$-module map which carries $\Filt_s(M)$ into $\Filt_s(N)$.
    A {\em filtered chain complex} $(C,\partial)$ is a graded and
    filtered $\Ring$-module, equipped with a filtered endomorphism
    $\partial$ which drops grading by one.  Given filtered chain
    complexes $A$ and $B$, a {\em filtered chain map} is a chain map
    $\phi\colon A \longrightarrow B$ which is a grading-preserving,
    filtered $\Ring$-module map.  Given two filtered chain maps
    $\phi_i\colon A \longrightarrow B$ for $i=1,2$, a {\em filtered
      chain homotopy} is a filtered $\Ring$-module map $H\colon A
    \longrightarrow B$ which raises grading by one and satisfies the
    formula $$\partial_B \circ H + H\circ \partial_A =
    \phi_1-\phi_2.$$
    If a filtered chain homotopy exists between
    $\phi_1$ and $\phi_2$, then we say that $\phi_1$ and $\phi_2$ are
    {\em filtered chain homotopic}.  Let $\phi\colon A \longrightarrow
    B$ be a filtered chain map. We say that $\phi$ is a {\em filtered
      chain homotopy equivalence} if there is a map $\psi\colon B
    \longrightarrow A$ with the property that $\phi\circ \psi$ and
    $\psi\circ \phi$ are filtered chain homotopic to the identity
    maps.  A {\em filtered quasi-isomorphism} is a filtered map
    $\phi\colon A \longrightarrow B$ which induces an isomorphism from
    the homology groups $H_*(\Filt_s(A))$ to $H_*(\Filt_s(B))$.  The
    {\em associated graded object} of a filtered chain complex $C$ is
    the $\Q^\ell$-graded chain complex 
  $$\gr(C)= \bigoplus_{s\in\Q^\ell} \gr_s(C),$$
  where $\gr_s(C)$ is the quotient of $\Filt_s(C)$ by 
  the submodule generated by $\Filt_t(C)$ for all $t < s$,
  endowed with the differential induced from $\partial$.
\end{definition}

 A filtered chain homotopy equivalence is a
filtered quasi-isomorphism.  Moreover a map is a filtered
quasi-isomorphism if and only if it induces an isomorphism on the
homology of the associated graded object.

\begin{definition}
  \label{def:MappingCone}
  Given a filtered chain map $\phi\colon A \longrightarrow B$, we can
  form a new filtered chain complex, the {\em mapping cone} $M(\phi)$
  whose underlying module is $A\oplus B$, and which is endowed with
  the differential $D(a,b)=(\partial a, \phi(a)-\partial b)$, where
  here $\partial a$ and $\partial b$
  denotes the differentials of $a$ and $b$ within $A$ and $B$, respectively.
\end{definition}

The mapping cone fits into a short exact sequence of chain complexes
(where the maps are all filtered chain maps) $$\begin{CD} 0@>>> B @>>>
M(\phi)@>>> A @>>> 0\end{CD},$$ and whose connecting homomorphism
agrees with the map induced by $\phi$.

\begin{definition}
  Two filtered chain complexes $A$ and $B$ are {\em quasi-isomorphic}
  if there is a third filtered chain complex $C$ and filtered
  quasi-isomorphisms from $C$ to $A$ and to $B$.
\end{definition}

If $\phi_1\colon A \longrightarrow B$ and $\phi_2 \colon A
\longrightarrow B$ are chain homotopic, then their induced mapping
cones are quasi-isomorphic.

Our chain complexes will always be finitely generated over
$\Field[U_1,\ldots,U_n]$.

\subsection{The chain complex $\Cm$}
\label{sec:chain-complex-cm}

We verify that $\Cm(G)$ as defined in the introduction (using
coefficients in $\Field$) is a filtered chain complex in the above
sense, with (Alexander) filtration induced from the function $A$ and
(Maslov) grading induced from the function $M$.

\begin{lemma}
  \label{lemma:MaslovWellDefined}
  The function $M$ is well-defined, i.e., it is independent of the manner in 
  which a given generator $\x\in\S$ is drawn on the square.
\end{lemma}

\begin{proof}
  Fix $\x\in\S$, thought of as drawn in the usual fundamental domain
  with the bottom and left edges included, so there is one
  component~$a$ with coordinates~$(m,0)$.  Let $\x'$ denote the same
  generator in the fundamental domain with the top and left edges
  included, so there is now a component~$b$ with coordinates~$(m,n)$.
  For each~$i$ with $0\le i<n, i \ne m$, there is one component $c_i$
  in~$\x$ and~$\x'$ with first coordinate~$i$.  For $m < i < n$, the
  pair $(a,c_i)$ contributes $1$ to the count of $\NESW(\x,\x)$,
  whereas the corresponding pair $(c_i,b)$ does not contribute to
  $\NESW(\x',\x')$.  Symmetrically, for each $i$ with $0\le i<m$, the
  pair $(c_i,a)$ does not contribute to $\NESW(\x,\x)$,
  whereas~$(c_i,b)$ does contribute to $\NESW(\x',\x')$. It follows
  that $\NESW(\x,\x) + m =\NESW(\x',\x') + n-m-1$. We can similarly
  analyze $\NESW(\x',\Os)$ to find
  \begin{align*}
    \NESW(\x',\x') &= \NESW(\x,\x) + 2m - n + 1\\
    2\NESW(\x',\Os) &= 2\NESW(\x,\Os) + 2m-n\\
  \end{align*}
  In particular $M_\Os(\x') = M_\Os(\x) + 1$.

  To complete the rotation, we have to change~$\Os$ to~$\Os'$ by
  moving the $O$ in the bottom row, with
  coordinates~$(l-\frac12,\frac12)$, to $(l-\frac12,n+\frac12)$.  A
  similar analysis yields
  \begin{align*}
    2\NESW(\x',\Os') &= 2\NESW(\x',\Os) + 2l-n\\
    \NESW(\Os',\Os') &= \NESW(\Os,\Os) + 2l - n - 1.
  \end{align*}
  Thus $M_{\Os'}(\x') = M_\Os(\x') - 1 = M_{\Os}(\x)$, which is the
  desired cyclic invariance.

  The same reasoning also establishes invariance under horizontal rotation.
\end{proof}

The Maslov grading on $\Ring$ and the generating set $\S$ induces a
Maslov grading on the chain complex $\Cm$. Explicitly, the summand
$\Cm_d(G)$ is generated by expressions $U_1^{m_1}\cdots
U_n^{m_n}\cm \x$, with $\x\in\S$, where
$$d=M(\x)-2\sum_{i=1}^n m_i.$$

\begin{lemma}
  \label{lemma:MaslovChange}
  Suppose that $\x,\y\in\S$, and $r\in\Rect(\x,\y)$ is a rectangle with
  $\x \cap \Interior(r) = \emptyset$. Then
  \begin{equation}
    \label{eq:MaslovChange}
    M(\x)=M(\y)+1-2\sum_{i=1}^n O_i(r).
  \end{equation}
\end{lemma}

\begin{proof}
  Draw the torus $\Torus$ on a square in such a manner that the lower
  left corner of $r$ coincides with the lower left corner of the
  square.  Then it is clear that
  $\NESW(\x,\x)=\NESW(\y,\y)+1$ (since the two new coordinates $y_1$ and
  $y_2$ in $\y$ are the only pair counted in $\NESW(\x)$ which are
  not also counted in $\NESW(\x)$), while $\NESW({\mathbb
    O},\x) =\NESW(\Os,\y)+\#\{\Os\cap r\}$, since each
  $O_i\in r$ gives rise to exactly one pair $(x_1,O_i)$ counted
  in $\NESW(\Os,\x)$ which is not also counted in
  $\NESW(\Os,\x)$. Similarly, $\NESW(\x,{\mathbb
    O})=\NESW(\y,\Os)+\#\{\Os\cap r\}$.
  Equation~\eqref{eq:MaslovChange} now follows when $M$ is calculated
  with respect to a particular manner of lifting the data on $\Torus$ 
  to data on a square. But according to Lemma~\ref{lemma:MaslovWellDefined},
  the Maslov grading is independent of this data.
\end{proof}

The alert reader might notice that the definition of Maslov grading we
give here does not identically agree with that given in~\cite{MOS},
which we denote by $M'$.  However, by connecting any two generators
$\x\in\S$ by a sequence of rectangles satisfying
Lemma~\ref{lemma:MaslovChange} (the existence of which can be deduced
from the fact that the symmetric group is generated by
transpositions), we see at once that $M$ is uniquely characterized, up
to an additive constant, by Equation~\eqref{eq:MaslovChange}, which is
also satisfied by $M'$. It now remains to show that $M(\x_0)=M'(\x_0)$
for some $\x_0\in\S$.  To this end, we take $\x_0$ to be the generator
for which $x_i$ is on the lower left corner of the square marked with
$O_i$. According to the conventions from~\cite{MOS},
$M'(\x_0)=1-n$; it is easy to verify that $M(\x_0)=1-n$, as well.

For the Alexander gradings, we have the following analogue of Lemma~\ref{lemma:MaslovWellDefined}:

\begin{lemma}
  \label{lemma:AlexanderWellDefined}
  For a given link component $i$, the function $A_i$ is well-defined,
  i.e., it is independent of the
  manner in which a given generator $\x\in\S$ is drawn on the square.
\end{lemma}

\begin{proof}
  For a point $p \in \Z^2$, the quantities $\NE(p,\Xs_i-\Os_i)$ and
  $\NE(\Xs_i-\Os_i, p)$ both compute the winding number of the $i\th$
  component of the knot around the point~$p$.  This quantity is
  unchanged if $p$ is moved from the very bottom to the very top of
  the diagram (since in that case the winding number is~$0$), and if
  $\Xs_i$ and $\Os_i$ are rotated vertically once, it changes by $\pm
  1$ if $p$ is in between the $X$ and the $O$ that are moved, and is
  unchanged otherwise.  For a point $p$ with half-integer coordinates,
  the inequalities used in the definition of $\NE(p,\Xs_i-\Os_i)$
  effectively shift $p$ up and to the right by
  $(\frac{1}{2},\frac{1}{2})$ before computing the winding number.
  Similarly, $\NE(\Xs_i - \Os_i, p)$ computes the winding number
  around $p - (\frac{1}{2},\frac{1}{2})$.  Therefore $A_i(\x)$,
  defined as $\NESW(\x - \frac{1}{2}(\Xs+\Os), \Xs_i - \Os_i)$,
  computes the winding number of the $i\th$ component around a
  weighted sum of points which has total weight~$0$ in each row and
  column. This combination is therefore invariant under cyclic
  rotation of the whole diagram.
\end{proof}

The function $A\colon \S(G) \longrightarrow \halfzl \subset \Q^{\ell}$ endows $\Cm(G)$
with a $\Q^{\ell}$-filtration in the sense of
Definition~\ref{def:Filtrations}, for the function $g\colon
\{1,\ldots,n\}\longrightarrow \Z^{\ell}$ which associates to $i$ the
$j\th$ standard basis vector in $\Z^{\ell}$ if $O_i$ belongs to the
$j\th$ component of the link.  The element $(U_1^{m_1}\cdots
U_n^{m_n})\x$ has filtration level $a = (a_1,\ldots,a_\ell)$, where
\[
a = A(\x) - \sum_{i=1}^n m_i \cm g(i).
\]

It is sometimes useful to consider objects more general than 
rectangles, called domains. To define them, let us view the torus $\Torus$ as a two-dimensional cell complex, with the toroidal grid diagram inducing the cell decomposition with $n^2$ zero-cells, $2n^2$ one-cells and $n^2$ two-cells (the little squares). Let $U_{\alpha}$ be the one-dimensional subcomplex of $\Torus$ consisting of the union of the $n$ horizontal circles.

\begin {definition}
Given $\x,\y\in \S$, a \emph{path} from~$\x$ to~$\y$ is a
1-cycle $\gamma$ on the cell complex $\Torus$, such that
the boundary of the intersection of $\gamma$ with $U_{\alpha}$ is
$\y-\x$.  
\end {definition}

\begin {definition}
A {\em domain}
$p$ from $\x$ to $\y$ is a two-chain in $\Torus$ whose boundary
$\partial p$ is a path from $\x$ to $\y$. The \emph{support} of $p$ is the
union of the closures of the two-cells appearing (with nonzero multiplicity) in the two-chain $p$.
\end {definition}

Given $\x,\y\in\S$, let
$\pi(\x,\y)$ denote the space of domains from $\x$ to $\y$.
There is a natural composition law
$$* \colon \pi({\mathbf a},{\mathbf
b})\times \pi({\mathbf b},{\mathbf c})\longrightarrow \pi({\mathbf
a},{\mathbf c}).$$

For a domain $p\in\pi(\x,\y)$, we let $X_i(p)$ and $O_i(p)$ denote the 
multiplicity with which $X_i$ and $O_i$, respectively, appear in $p$.

\begin{prop}
  The differential $\dm$ drops Maslov grading by one, and respects the
  Alexander filtration.  Specifically, if $\x\in\S$ has $M(\x)=d$,
  then $\dm(\x)$ is written as a sum of elements in Maslov grading
  $d-1$.  Also, if $A(\x)=a$, then $\dm(\x)$ is a sum of elements with
  Alexander filtrations $\leq a$.
\end{prop}

\begin{proof}
  The fact that $\dm$ drops Maslov grading by one follows at once from
  Equation~\eqref{eq:MaslovChange}, together with the definition of
  $\dm$. 

  The fact that $\dm$ respects the Alexander filtration follows from
  basic properties of winding numbers. Specifically, given
  $\x,\y\in\S$ and $r\in \Rect(\x,\y)$, it is easy to see that
  $$A(\x)-A(\y)=\sum_{i} (X_i(r)-O_i(r)) \cdot g(i).$$

  Thus if
  $U_1^{m_1}\cdots U_n^{m_n}\cm \y$ appears with non-zero
  coefficient in $\dm (\x)$, then the Alexander filtration level of
  the corresponding term is smaller than the Alexander filtration
  level of $\x$ by $\sum_{i=1}^n X_i(r)\cdot g(i)$.
\end{proof}

With the terminology in place, we now verify that $\dm$ is the
differential of a chain complex.

\begin{prop}
  \label{prop:DSquaredZero}
  The endomorphism $\dm$ of $\Cm(G)$ is a differential,
  i.e., $\dm\circ\dm=0$.
\end{prop}

\begin{proof}
        Consider an element $\x\in\S$, viewed as a generator of
        $\Cm(G)$.  We can view $\dm\circ\dm(\x)$ as a count $$
        \dm\circ\dm(\x)
        =\sum_{\z\in\S}\,\sum_{\substack{p\in\pi(\x,\y)\\
        \x\not\in\Interior p}} N(p)\cm U_1^{O_1(p)}\cdots
        U_n^{O_n(p)}\cm \z.  $$ where here $N(p)$ denotes the number
        of ways of decomposing a domain as a composite of two empty 
        rectangles $p=r_1*r_2$, where $r_1\in\EmptyRect(\x,\y)$ and
        $r_2\in\EmptyRect(\y,\z)$ for some $\y\in\S$.

\begin{figure}
\begin{center}
\mbox{\vbox{\epsfbox{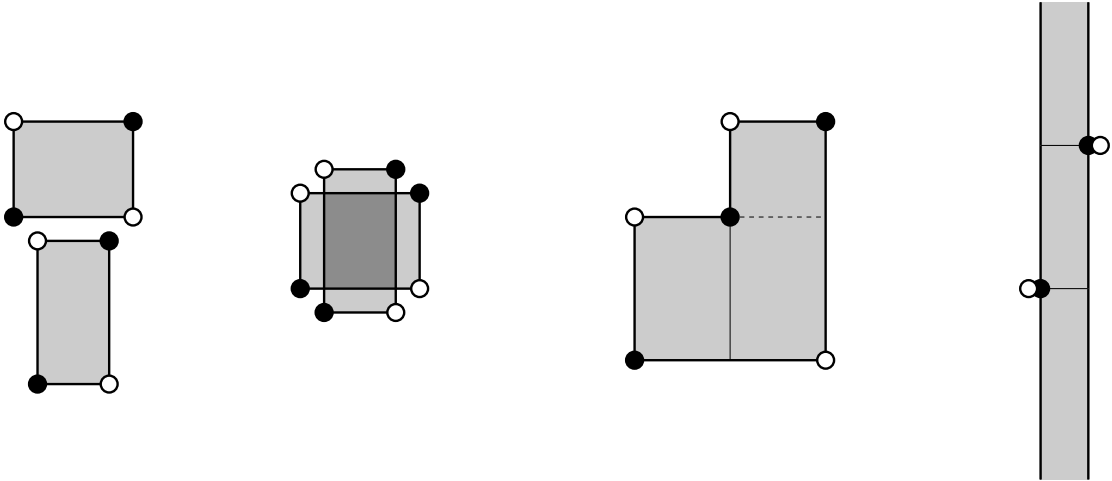}}}
\end{center}
\caption {{\bf $\dm\circ\dm=0$.}
The four combinatorially different ways the composite of two empty
rectangles $r_1*r_2$ can appear.  The initial point is indicated
by the dark circles, the final by the hollow ones.}
\label{fig:DSquaredZero}
\end{figure}

        If $\z\neq\x$, and if $p$ has a decomposition $p=r_1*r_2$,
        then we claim that there is a unique alternate decomposition
        $p=r_1'*r_2'$, where here $r_1'\in\Rect(\x,\y')$ and
        $r_2'\in\Rect(\y',\z)$. In fact, if $p=r_1*r_2$ is a domain
        obtained from two empty rectangles $r_1$ and $r_2$, then we
        claim that there are three possibilities for $p$:
        \begin{itemize} \item two disjoint rectangles; \item two
        rectangles with overlapping interiors (the darker region in
        Figure~\ref{fig:DSquaredZero}); and \item two rectangles which share a corner.
        \end{itemize} These three cases are illustrated in the first
        three diagrams in Figure~\ref{fig:DSquaredZero}.  In each
        case, there are exactly two decompositions of the obtained
        domain as a juxtaposition of empty rectangles: in the first
        two cases by taking the rectangles in the two possible orders,
        and in third case by decomposing either along the thin or
        dotted lines, cf.\ Figure~\ref{fig:Ell}. It follows at once that
        the $\z$ component of $\dm\circ\dm(\x)$ vanishes for $\z\neq \x$.

\begin{figure}
\begin{center}
\mbox{\vbox{\epsfbox{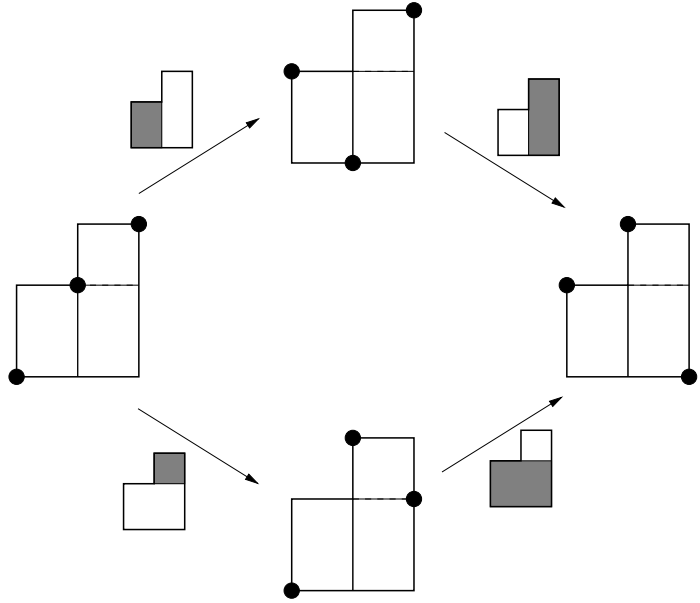}}}
\end{center}
\caption {{\bf The third case of Figure~\ref{fig:DSquaredZero}.}
The three black dots are permuted to give four different generators.
Each arrow represents a rectangle, which is shown shaded. There are two 
ways of connecting the initial generator $\x$ to the final generator $\z$: 
by following the top arrows, or the bottom ones. Each way gives a 
contribution to  $\dm\circ\dm$, and in the final count these contributions 
cancel out.} 
\label{fig:Ell} \end{figure}

        When $\z=\x$, however, the only domains $p\in\pi(\x,\x)$ which
        can be decomposed as a union of two empty rectangles are width
        one annuli, as in the fourth
        diagram in Figure~\ref{fig:DSquaredZero},
         or height one annuli in the torus.  There are $2n$ of
        these annuli. Each such annulus $p$ has a unique decomposition
        $p=r_1*r_2$ with $r_1\in\Rect(\x,\y)$ and $r_2\in\Rect(\y,\x)$
        (for some uniquely specified $\y$). The row or column
        containing $O_i$ contributes $U_i$ in the formula for
        $\dm\circ\dm(\x)$. Since $O_i$ appears in exactly one row and
        exactly one column, it follows now that the $\x$ component of
        $\dm\circ\dm(\x)$ vanishes, as well.
\end{proof}

The proof of the above proposition is elementary, depending on evident
properties of rectangles in the torus. However, it does deserve a few
extra words, since it is the starting point of this paper, and indeed
a recurring theme throughout.  Specifically, the alert reader will
observe that the remarks concerning juxtapositions of pairs of
rectangles is one of the last vestiges of Gromov's compactness
theorem, the foundation upon which Floer's theory of Lagrangian
intersections is built~\cite{FloerLag} (and knot Floer homology can be
viewed as a variant of that latter theory). The assertions
about annuli can also be seen as remnants of Gromov's theory, as they
are counting boundary degenerations.

In terms of combinatorics, we see a pattern that will be repeated
throughout the paper: in order to prove an identity with differentials
(e.g., that $(\dm)^2 = 0$, or that a map is a chain map) we consider
the composites of two domains; generally the composite domain will have
exactly two decompositions.  In some cases we need to add or delete
annuli of width or height one while taking care not to change the
factors of~$U_i$ that appear.

\subsection{Algebraic properties of $\Cm$}

We now turn to the basic algebraic properties of the chain complex.

In the following lemma, just as $U_i$ is a chain map which drops
filtration level by one, the filtered chain homotopy drops the
filtration level by one.

\begin{lemma}
  \label{lemma:ChainHomotopies}
  Suppose that $O_i$ and $O_k$ correspond to the same component of
  $\oL$.  Then multiplication by $U_i$ is filtered chain homotopic
  to multiplication by $U_k$.
\end{lemma}

\begin{proof}
  Since filtered chain homotopies can be composed, it suffices to show
  that if $O_i$ lies in the same row as some $X_j$ which in turn is in
  the same column as $O_k$, then multiplication by $U_i$ is filtered
  chain homotopic to multiplication by $U_k$. The filtered chain
  homotopy is furnished by counting rectangles which contain $X_j$. 

  Specifically, define
  $$H\colon \Cm(G)\longrightarrow \Cm(G).$$
  by the formula
  $$H(\x)=\sum_{\y\in\S}\,
   \sum_{\substack{r\in\EmptyRect(\x,\y)\\X_j\in r}}\!\!
    U_1^{O_1(r)}\cdots U_n^{O_n(r)}\cm \y.$$
  We claim that 
  $$\dm\circ H + H\circ\dm = U_i-U_k.$$
  This follows from the same argument as
  Proposition~\ref{prop:DSquaredZero}: Most composite domains on the
  left hand side can be decomposed in exactly two ways.  The exception
  are the horizontal and vertical annuli, necessarily containing $X_j$
  which contribute $U_i$ and $U_k$, respectively.
\end{proof}

\begin{prop}
  \label{prop:DoesntDependOnOrdering}
  Suppose that the oriented link $\oL$ has $\ell$ components.  Choose
  an ordering of $\Os=\{O_i\}_{i=1}^n$ so that for
  $i=1,\ldots,\ell$, $O_i$ corresponds to the $i\th$ component of
  $\oL$.  Then the filtered chain homotopy type of $\Cm(G)$,
  viewed as a chain complex over $\Field[U_1,\ldots,U_\ell]$, is
  independent of the ordering of $\Os$.
\end{prop}

\begin{proof}
  Different numberings can be connected via the filtered chain homotopies
  of Lemma~\ref{lemma:ChainHomotopies}.
\end{proof}

The basic link invariant is the filtered quasi-isomorphism class of $\Cm(\oL)$,
thought of as a complex of $\Field[U_1,\ldots,U_\ell]$ modules. But
there are some other natural constructions one can
consider.

For example, we can consider the chain complex $\Ca(G)$, which is
a chain complex over $\Field$, once again which is freely generated by
elements of $\S$, by setting the $U_i=0$ for $i=1,\ldots,\ell$.
We let~$\CLa(G)$ denote the graded object~$\gr(\Ca(G))$ associated
to the Alexander filtration, and let $\HLa(G)$ denote its homology.

\begin{lemma}
  \label{lemma:HLaFinGen}
  The group $\HLa(G)$ is a finitely-generated
  $\Field$-module.
\end{lemma}

\begin{proof}
  Clearly, $\CLa(G)$ is a finitely generated $\Ring$-module. It follows
  from Lemma~\ref{lemma:ChainHomotopies} that once we set $U_i=0$ for
  $i=1,\ldots,\ell$, then multiplication by $U_j$ is null-homotopic
  for all $j=1,\ldots,n$, and in particular it acts trivially on
  homology. It follows at once that $\HLa(G)=H_*(\CLa(G))$, which
  is clearly a finitely generated $\Ring$-module, is in fact a
  finitely generated $\Field$-module.
\end{proof}

There is another construction which is quite convenient to consider
for calculations~\cite{BaldwinGillam}.  This is the chain complex
$\Caa(G)$, which is obtained from $\Cm(G)$ by setting all the $U_i=0$,
and then taking the associated graded object.  (This complex is
denoted simply $C(G)$ in~\cite{MOS}, but we prefer to reserve this
notation for later use.)  Explicitly, this is the free $\Field$-module
generated by $\S$, endowed with the differential
$$\daa(\x)=\sum_{\y\in\S}\# \left\{r\in\Rect(\x,\y)\Big|
\begin{array}{l}
\forall x\in \x, x\not\in\Interior(r), \\
\forall i, O_i\not\in r~\text{and}~X_i\not\in r
\end{array}\right\} \cm
\y.$$
It is easy to relate the homology of $\CLaa(G) = \gr(\Caa(G))$
with the homology of
$\CLa(G)$, by some principles in homological algebra.

\begin{lemma}
        \label{lemma:UniversalCoefficients}
        Let $C$ be a filtered, graded chain complex of free modules over
        $\Field[U_1,\ldots,U_n]$, such that $U_i$ decreases the
        homological grading by two and the filtration by one, and such
        that multiplication by $U_i$ is
        chain homotopic to multiplication by $U_j$ for any $i,j$.
        Then $H_*(C/\{U_i=0\}_{i=1}^n)\cong H_*(C/U_1)\otimes
        V^{n-1}$, where $V$ is the two-dimensional bi-graded vector
        space spanned by one generator in bi-grading $(-1,-1)$ and
        another in bi-grading $(0,0)$.
\end{lemma}

\begin{proof}
  Suppose for notational simplicity that $n=2$.  Consider the chain
  map from the mapping cone of the chain map $U_1\colon C
  \longrightarrow C$ to $C/U_1$ gotten by taking the quotient on the
  second summand.  It follows easily from the five-lemma that this map
  is a quasi-isomorphism. Moreover, by iterating this observation, we
  see that $C/(U_1,U_2)$ is quasi-isomorphic to the mapping cone
  $$\begin{CD} C @>{U_1}>>C \\ @V{U_2}VV @VV{U_2}V \\ C@>{U_1}>> C,
  \end{CD} $$
  which in turn is quasi-isomorphic to the mapping cone of
  $$U_2 \colon C/U_1\longrightarrow C/U_1.$$
  But since $U_1$ and $U_2$
  are chain homotopic in $C$, we obtain an induced null-homotopy of
  the map induced by $U_2$ on $C/U_1$.  Thus, this latter mapping cone
  is isomorphic to the mapping cone of zero, i.e., to the direct sum
  $C/U_1\oplus C/U_1$, which in turn is quasi-isomorphic to
  $(C/U_1)\otimes V$. 
        
  We investigate now the filtrations and gradings.  In order for the
  quasi-isomorphism from $U_1\colon C \longrightarrow C$ to $C/U_1$ to
  be a filtered and graded map, we must shift gradings an filtrations
  on the mapping cone $M(U_1)$ appropriately. Specifically, let
  $C[a,b]$ denote the graded and filtered chain complex with the
  property that $\Filt_s(C_d[a,b])=\Filt_{s+b}(C_{d+a})$. Then the
  mapping cone $M(U_1)$ is $C[1,1]\oplus C$. Following through the
  above discussion, we see that the mapping cone $C/(U_1,U_2)$ is
  filtered and graded quasi-isomorphic to $C[1,1]/U_1\oplus C/U_1\cong
  (C/U_1)\otimes V$.

  This discussion generalizes readily to the case
  where $n>2$.
\end{proof}
        
\begin{prop}
  \label{prop:HLaaNoSigns}
  The homology groups $\HLa(G)$ determine $\HLaa(G)$;
  specifically,
  \[H_*(\CLaa(G)) \cong \HLa(G)\otimes
  \bigotimes_{i=1}^\ell V_i^{\otimes (n_i-1) }, \] where $V_i$ is the
  two-dimensional vector space spanned by two generators, one in zero
  Maslov and Alexander multi-gradings, and the other in Maslov grading
  minus one and Alexander multi-grading corresponding to minus the
  $i\th$~basis vector.
\end{prop}

\begin{proof}
        This follows easily from
        Lemma~\ref{lemma:UniversalCoefficients}, applied component by
        component.
\end{proof}

\noindent{\bf Notation.}
Perhaps the reader will find it convenient if we collect our
notational conventions here.  The chain complex $\Cm(G)$ refers to the
full chain complex (and indeed, we soon drop the minus from the
notation here), $\CLm(G)$ denotes its associated graded object, and
$\HLm(G)$ is the homology of the associated graded object. $\Ca(G)$
denotes the chain complex where we set one $U_i=0$ for each component of the
link, $\CLa(G)$ is its associated graded object, and $\HLa(G)$ is the
homology of the associated graded object.  $\Caa(G)$ is the chain
complex $\Cm(G)$ modulo the relations that {\em every} $U_i=0$, $\CLaa(G)$ is
the associated graded complex, and $\HLaa(G)$ is its homology.  Most
of these constructions have their analogues in Heegaard Floer
homology; for example, according to~\cite{MOS}, $\HLm(G)$ is
identified with $\HFLm(L)$, and $\HLa(G)$ with $\HFLa(L)$.  We find it
useful to distinguish these objects, especially when establishing
properties of the combinatorial complex which could alternatively be
handled by appealing to~\cite{MOS}, together with known properties of
Heegaard Floer homology.


\section{Invariance of combinatorial knot Floer homology}
\label{sec:Invariance}

Our goal in this section is to use elementary methods to show that
combinatorial knot Floer homology is independent of the grid diagram,
proving Theorem~\ref{thm:WithSigns} with coefficients in $\Field$.

\begin{figure}
\begin{center}
\mbox{\vbox{\epsfbox{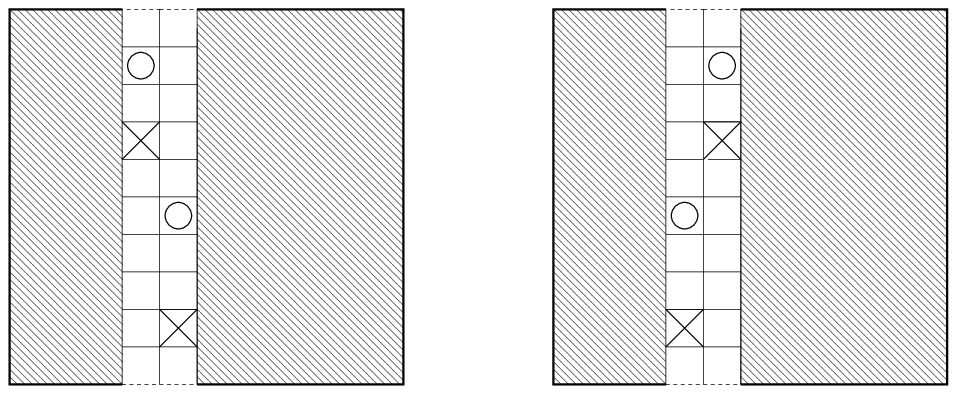}}}
\end{center}
\caption {{\bf Commutation.}
The two grid diagrams differ from each other by interchanging the two 
columns, but correspond to the same link.
}
\label{fig:CommutationPicture}
\end{figure}

\begin{figure}
\begin{center}
\input{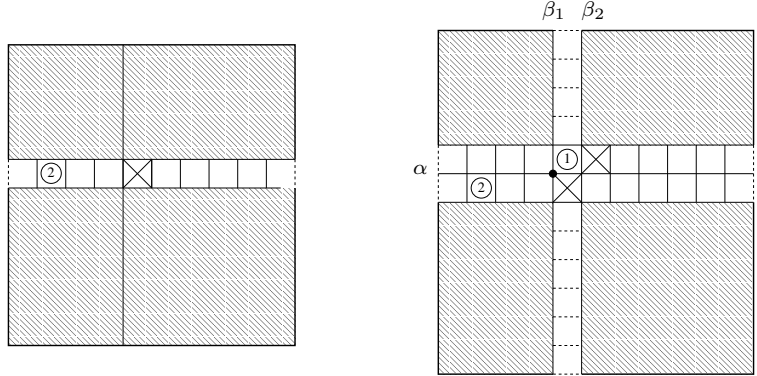}
\end{center}
\caption {{\bf Stabilization.}
On the left, we have an initial grid diagram; on the right, a new
diagram obtained by inserting the pictured row and column. Another
stabilization is given by switching the roles of the new middle two rows.}
\label{fig:StabilizationPicture}
\end{figure}

Following Cromwell~\cite{Cromwell} 
(compare also 
Dynnikov~\cite{Dynnikov}), any two 
grid diagrams for the same link can be connected by a sequence of
the following elementary moves:
\begin{enumerate}
\item {\bf Cyclic permutation.}
  This corresponds to cyclically permuting the rows and then the columns
  of the grid diagram.
\item {\bf Commutation.}  Consider a pair of consecutive columns 
  in the grid diagram $G$ with the following property: if we
  think of the $X$ and the $O$ from one column as separating the
  vertical circle into two arcs, then the $X$ and the $O$ from the
  adjacent column occur both on one of those two arcs. Under these
  hypotheses, switching the decorations of these two columns is a
  commutation move, cf.\ Figure~\ref{fig:CommutationPicture}.  There is 
  also a similar move where the roles  of columns and rows are 
  interchanged.
\item {\bf Stabilization/destabilization.}  Stabilization is gotten by adding two
  consecutive breaks in the link.  More precisely, if $G$ has arc
  index $n$, a stabilization $H$ is an arc index $n+1$ grid diagram
  obtained by splitting a row in $G$ in two and introducing a new
  column.  For convenience, label the original diagram so it has
  decorations $\{X_i\}_{i=2}^{n+1}$, $\{O_i\}_{i=2}^{n+1}$.  Let $O_i$
  and $X_i$ denote the two decorations in the original row. We copy
  $O_i$ onto one of the two new copies of the row it used to occupy,
  and copy $X_i$ onto the other copy.  We place decorations $O_1$ and
  $X_1$ in the new column so $O_1$ resp. $X_1$ occupy the same row as
  $X_i$ resp. $O_i$ in the new diagram, cf.\ Figure~\ref{fig:StabilizationPicture}. 
  Destabilization is the inverse move to stabilization.  Note that
  stabilization can be alternatively done by reversing the roles of
  rows and columns in the above description; however, such a
  stabilization can be reduced to the previous case, combined with a
  sequence of commutation moves.  In fact, we can consider only
  certain restricted stabilization moves, where three of the four
  squares $O_1$, $X_1$, $O_i$, and $X_i$ share a common vertex; i.e.,
  the new column is introduced next to $O_i$ or $X_i$.  However, there
  are now different types of stabilizations corresponding to the different
  ways of dividing the $O$'s and $X$'s among the two new rows.

\end{enumerate}

Of course, since our complex is associated not to the planar grid diagram, but
rather to the induced picture on the torus, the fact that it is invariant under
cyclic permutation is a tautology.

We turn to commutation invariance
next, and then stabilization invariance.

Note that all the chain complexes $\Ca(G)$, $\Caa(G)$ depend on the
quasi-isomorphism type of $\Cm(G)$; thus, the latter is the most basic
object. Thus, to streamline notation, we choose here to drop the
superscript `$-$' from the notation of this chain complex and its
differential.

\subsection{Commutation invariance}
\label{subsec:Commutation}

\begin{figure}
\begin{center}
\input {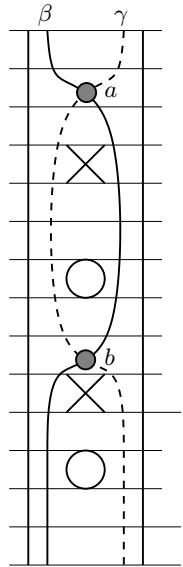}
\end{center}
\caption {{\bf Commutation.}
A commutation move, viewed as replacing one vertical circle ($\beta$, 
undashed) with another ($\gamma$, dashed).}
\label{fig:CommutePicture}
\end{figure}

Let $G$ be a grid diagram for $\oL$, and let $H$ be a different grid
diagram obtained by commuting two vertical edges. It is convenient to
draw both diagrams on the same torus, replacing a distinguished
vertical circle $\beta$ for $G$ with a different one $\gamma$ for $H$,
as pictured in Figure~\ref{fig:CommutePicture}.  The circles $\beta$
and $\gamma$ meet each other transversally in two points~$a$ and~$b$,
which are not on a horizontal circle.

\begin{figure}
\begin{center}
  \mbox{\vbox{\epsfbox{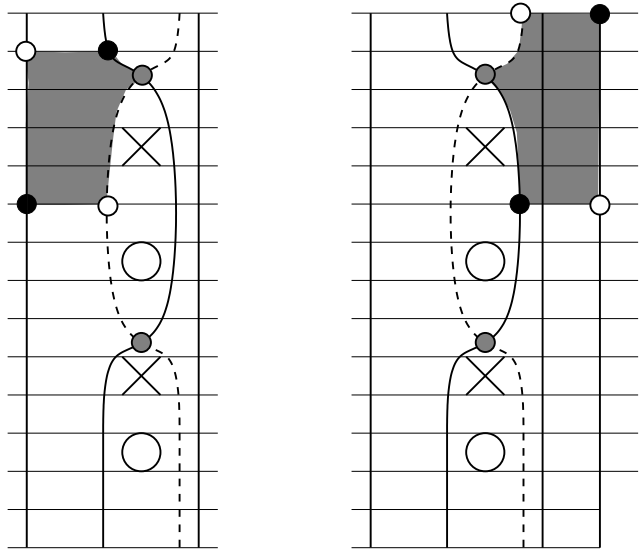}}}
\end{center}
\caption {{\bf Pentagons.} We have indicated here two
allowed pentagons in $\Pent_{\beta\gamma}(\x,\y)$, where components
of $\x$ are indicated by solid points, and those of $\y$ are indicated
by hollow ones.}
\label{fig:Pentagon}
\end{figure}

We define a chain map $\Phi_{\beta\gamma} \colon C(G) \longrightarrow
C(H)$ by counting pentagons.  Given $\x\in \S(G)$ and $\y\in\S(H)$, we
let $\Pent_{\beta\gamma}(\x,\y)$ denote the space of embedded
pentagons with the following properties. This space is empty unless
$\x$ and $\y$ coincide at $n-2$ points. An element of
$\Pent_{\beta\gamma}(\x,\y)$ is an embedded disk in $\Torus$, whose
boundary consists of five arcs, each contained in horizontal or
vertical circles. Moreover, under the orientation induced on the
boundary of $p$, we start at the $\beta$-component of $\x$, traverse
the arc of a horizontal circle, meet its corresponding component of
$\y$, proceed to an arc of a vertical circle, meet the corresponding
component of $\x$, continue through another horizontal circle,
meet the component of $\y$ contained in the distinguished circle
$\gamma$, proceed to an arc in $\gamma$, meet an intersection point of
$\beta$ with $\gamma$, and finally, traverse an arc in $\beta$ until
we arrive back at the initial component of $\x$.  Finally, all the
angles here are required to be less than straight angles. These
conditions imply that
there is a particular intersection point, denoted $a$, between $\beta$
and $\gamma$ which appears as one of the corners of any pentagon in
$\Pent_{\beta\gamma}(\x,\y)$.  The other intersection point $b$
appears in all of the pentagons in $\Pent_{\gamma\beta}(\y,\x)$.
Examples are pictured in Figure~\ref{fig:Pentagon}.  The space of
empty pentagons $p\in \Pent_{\beta\gamma}(\x,\y)$ with $\x \cap \Interior(p) =
\emptyset$, is denoted $\EmptyPent_{\beta\gamma}$.

Given $\x\in\S(G)$, define
\begin{align*}
\Phi_{\beta\gamma}(\x) =& \sum_{\y\in \S(H)}\,\,
    \sum_{p\in\EmptyPent_{\beta\gamma}(\x,\y)}\!\!
U_1^{O_1(p)}\cdots U_n^{O_n(p)} \cm \y
\in C(H).
\end{align*}

\begin{lemma}
  \label{lemma:PhiChainMap}
  The map $\Phi_{\beta\gamma}$ is a filtered chain map.
\end{lemma}

\begin{figure}
\begin{center}
  \mbox{\vbox{\epsfbox{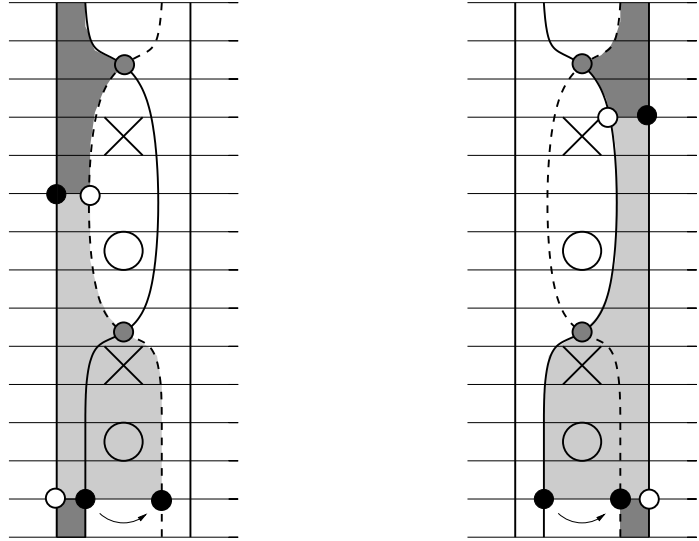}}}
\end{center}
\caption {{\bf Special case.}
  The generators $\x$ and $c(\x)$ are marked by dark circles; they differ 
from each other only on one row. The arrow indicates how the dark circle 
in $\x$ is replaced by a corresponding dark circle in $c(\x)$. On the left 
we have a (darkly-shaded) pentagon followed by a (lightly-shaded) 
rectangle, and on the right we have a rectangle followed by a pentagon. 
The intermediate generators 
are marked by hollow circles.} \label{fig:SpecialCase}
\end{figure}

\begin{proof}
  The fact that $\Phi_{\beta\gamma}$ preserves Alexander filtration and
  Maslov gradings is straightforward.
  Like the proof of
  Proposition~\ref{prop:DSquaredZero}, the proof that
  $\Phi_{\beta\gamma}$ is a chain map proceeds by considering domains
  which are obtained as a juxtaposition of a pentagon and a rectangle,
  representing terms in $\partial \circ \Phi_{\beta\gamma}$, and
  observing that such domains typically have an alternate decomposition
  to represent a term in $\Phi_{\beta\gamma}\circ \partial$.  One
  example is illustrated in Figure~\ref{fig:ChainMap}.  Other terms
  are more straightforward, consisting either of a disjoint rectangle
  and pentagon, a rectangle and pentagon with overlapping interior, or
  a rectangle and a pentagon which meet along a different edge; the
  pictures are similar to those in Figure~\ref{fig:DSquaredZero}.
\begin{figure}
\begin{center}
  \mbox{\vbox{\epsfbox{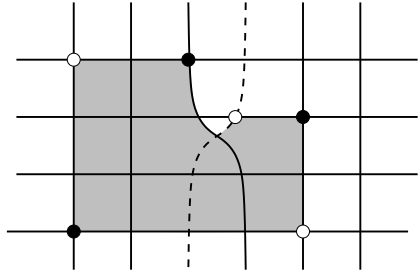}}}
\end{center}
\caption {{\bf Chain map.} 
  The given domain can be decomposed either as a pentagon followed by
  a rectangle, or a rectangle followed by a pentagon.  The first
  decomposition represents a term in $\partial \circ \Phi_{\beta\gamma}$,
  the second a term in $\Phi_{\beta\gamma}\circ \partial$.}
\label{fig:ChainMap}
\end{figure}
There is one special case, of a type of domain which has only one
decomposition: these are the domains obtained as the union of a width
one pentagon $p$ and a width one rectangle $r$. In this case, if we
let $\x\in\S(G)$, there is a canonical closest generator $c(\x)\in\S(H)$
(with the property that $\x$ and $c(\x)$ agree at all intersection
points away from $\beta\cup\gamma$). It is easy to see, then, that our
domain has the form $r*p$ or $p*r$ (depending on the local picture of 
$\x$), and it connects $\x$ to $c(\x)$. But then, such
domains are in one-to-one correspondence with domains of the form
$r'*p'$ or $p' * r'$, where if $p$ is a left pentagon, then $p'$ is a 
right pentagon, and vice versa. See Figure~\ref{fig:SpecialCase}.

\end{proof}

We can define chain homotopy operators analogously, only now counting hexagons.

\begin{figure}
\begin{center}
  \mbox{\vbox{\epsfbox{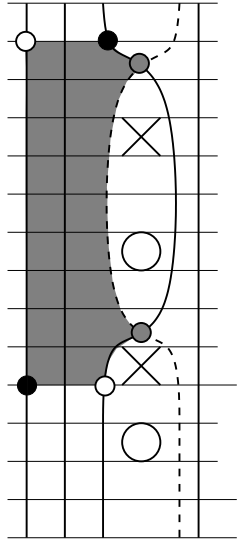}}}
\end{center}
\caption {{\bf Hexagon.} We have illustrated here 
a hexagon in  $\Hex_{\beta\gamma\beta}$.}
\label{fig:Hexagon}
\end{figure}

More specifically, given $\x,\y\in \S(G)$, we let
$\Hex_{\beta\gamma\beta}(\x,\y)$ denote the space of embedded hexagons
with the following property. This space, too, is empty unless $\x$ and
$\y$ coincide at $n-2$ points. Moreover, an element of
$\Hex_{\beta\gamma\beta}(\x,\y)$ is an embedded disk in $\Torus$,
whose boundary consists of six arcs, each contained in horizontal
or vertical circles. More specifically, under the orientation induced
on the boundary of $p$, we start at the $\beta$-component of $\x$,
traverse the arc of a horizontal circle, meet its corresponding
component of $\y$, proceed to an arc of a vertical circle, meet its
corresponding component of $\x$, continue through another
horizontal circle, meet its component of $\y$, which contained in the
distinguished circle $\beta$, continue along $\beta$ until the
intersection point $b$ of $\beta$, with $\gamma$, continue on $\gamma$
to the intersection point $a$ of $\beta$ and $\gamma$, proceed again
on $\beta$ to the the $\beta$-component of $\x$, which was also our
initial point.  Moreover, all corner points of our hexagon are
again required to be less than straight angles.  An example is given in
Figure~\ref{fig:Hexagon}. We define the space of empty
hexagons~$\EmptyHex_{\beta\gamma\beta}$, with
interior disjoint from $\x$, as before.  There is also a corresponding
notion $\Hex_{\gamma\beta\gamma}$. We now define the function
$H_{\beta\gamma\beta}\colon C(G) \longrightarrow C(G)$ by
\[
  H_{\beta\gamma\beta}(\x) = \sum_{\y\in \S(G)}\,
    \sum_{h\in\EmptyHex_{\beta\gamma\beta}(\x,\y)}\!
    U_1^{O_1(h)}\cdots U_n^{O_n(h)} \cm \y.
\]

\begin{prop}
  \label{prop:Commute}
  The map $\Phi_{\beta\gamma}\colon C(G)\longrightarrow C(H)$ is a chain
  homotopy equivalence; more precisely
  \begin{align*}
    \Id+\Phi_{\gamma\beta}\circ \Phi_{\beta\gamma} +  \partial \circ H_{\beta\gamma\beta} + H_{\beta\gamma\beta}\circ \partial  &= 0 \\
    \Id+\Phi_{\beta\gamma}\circ \Phi_{\gamma\beta} + \partial \circ H_{\gamma\beta\gamma} + H_{\gamma\beta\gamma}\circ \partial &=0.
\end{align*}
\end{prop}

\begin{figure}
\begin{center}
  \mbox{\vbox{\epsfbox{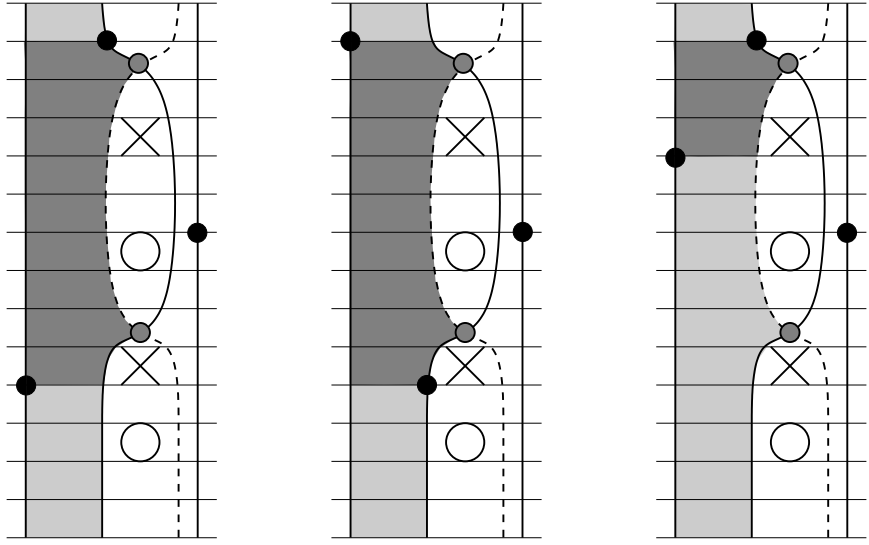}}}
\end{center}
\caption {{\bf Decomposing the identity map.} 
  Consider the three configurations in $C(G)$, indicated by dark
  circles.  The shaded region can be thought of as decomposed into a
  hexagon followed by a rectangle (as on the left), a rectangle
  followed by a hexagon (as in the middle), or a pair of pentagons as
  on the right. The first can be thought of as counting terms in
  $\partial \circ H_{\beta\gamma}$, the middle terms in
  $H_{\beta\gamma}\circ\partial$, and the right in
  $\Phi_{\gamma\beta}\circ \Phi_{\beta\gamma}$. There are three more
  cases, if the $\beta$-component of the configuration lies on the
  other arc in $\beta$; in this case, we must decompose the annulus on
  the right.}
\label{fig:DecomposeIdentity}
\end{figure}

\begin{proof}
  Juxtaposing two pentagons appearing in $\Phi_{\gamma\beta}\circ
  \Phi_{\beta\gamma}$, we generically obtain a composite domain which
  admits a unique alternative decomposition as a hexagon and a square,
  counted in $\partial\circ H_{\beta\gamma\beta}$ or
  $H_{\beta\gamma\beta}\circ\partial$.  Typically, the remaining terms
  in $\partial\circ H_{\beta\gamma\beta}$ cancel with terms
  $H_{\beta\gamma\beta}\circ\partial$.
  
  There is, however, one composite region which has a unique
  decomposition.  Specifically, the vertical circles
  $\beta_1$, $\beta_2$,  and $\gamma$ divide up $\Torus$ into a collection of
  components, two of which are annuli and do not contain any $X$.
  Depending on the initial point $\x$, exactly one of these annuli can
  be thought of as a juxtaposition of two pentagons, or a hexagon and
  a rectangle which is counted once in $\Phi_{\gamma\beta}\circ
  \Phi_{\beta\gamma}+ \partial \circ
  H_{\beta\gamma\beta}+H_{\beta\gamma\beta}\circ \partial$; but it is
  also counted in the identity map.  See Figure~\ref{fig:DecomposeIdentity}.
\end{proof}

\subsection{Stabilization invariance}
\label{subsec:Stabilization}

Let $G$ be a grid diagram and $H$ denote a stabilization.  We discuss
in detail the case where we introduce a new column with $O_1$
immediately above $X_1$ (and $X_2$ is immediately to the left or to
the right of $O_1$); the case where $X_1$ is immediately above $O_1$
can be treated symmetrically by a rotation of all diagrams by $180^\circ$.

More specifically, given a horizontal arc from $O_2$ to $X_2$, we
introduce a vertical segment (somewhere along the arc) consisting of a
new pair $O_1$ and $X_1$, where $O_1$ is on the square right above
$X_1$, which in turn is in the same row as the new copy of $O_2$, as
in Figure~\ref{fig:StabilizationPicture}. Indeed, do this in such a
manner that three of the four squares marked $O_1$, $O_2$, $X_1$, and
$X_2$ share a common vertex. Furthermore, by applying commutation, we
can assume without loss of generality that these three squares are
$O_1$, $X_1$, and $X_2$.  Thus, the grid of $H$ is gotten by inserting
a new column of squares, where two consecutive squares are marked by
$O_1$ and $X_1$. We let $\beta_1$ be the vertical circle on the left,
and $\beta_2$ the one on the right. Let $\alpha$ denote the new
horizontal circle in $H$ which separates $O_1$ from $X_1$.

Let $B=\Cm(G)$ and $C=\Cm(H)$. 
Let $C'$ be the mapping cone of 
$$U_2-U_1\colon B[U_1]\longrightarrow B[U_1],$$
i.e.,
$C'[U_1]=B[U_1]\oplus B[U_1]$, endowed with the differential
$\partial\colon C' \longrightarrow C'$ given by
$$
\partial'(a,b)=(\partial a, (U_2-U_1) \cm a-\partial b) $$
where
here $\partial$ denotes the differential within $C$ (actually, in the
sequel we drop the prime from the differential within $C'$, as well,
and hope that the differential is clear from the context).  Note that
$B$ is a chain complex over $\Field[U_2,\ldots,U_n]$, so that $B[U_1]$
denotes the induced complex over $\Field[U_1,\ldots,U_n]$ gotten by
introducing a new formal variable $U_1$.  Let $\Left$ and 
$\Right\cong B[U_1]$ be the subgroups of $C'$ of elements of the form
$(c,0)$ and $(0,c)$ for $c\in B[U_1]$, respectively.  The module
$\Right$ inherits
Alexander and Maslov gradings from its identification with $B[U_1]$,
while $\Left$ is given the Alexander and Maslov gradings which are one
less than those it inherits from its identification with $B[U_1]$.
With respect to these conventions, the mapping cone is a filtered
complex of $\Ring$-modules.

\begin{lemma}
  \label{lem:MappingCone}
  The map from $C'$ to $B$ that takes $(a,b)$ to $a/\{U_1 = U_2\}$ is
  a quasi\hyp isomorphism.
\end{lemma}

\begin{proof}
  In general, the mapping cone~$C'$ of a map $f: C_1 \rightarrow C_2$ fits
  into a short exact sequence on homology from $C_2$ to $C'$ to
  $C_1$.  The connecting homomorphism in the corresponding long exact
  sequence on homology is the map induced by~$f$.  In this case, $f$
  is $U_1 - U_2$, which is injective on the homology of $B[U_1]$, so
  the map from $C'$ to $B$ is a quasi-isomorphism.
\end{proof}

It therefore suffices to define a filtered quasi-isomorphism
\begin{equation}
  \label{eq:StabilizationMap}
  F\colon C \longrightarrow C'.
\end{equation}
To do this, we introduce a little more notation.

Let $\S(G)$ be the generating set of $B$, and $\S(H)$ be the generating
set of $C$. Let $x_0$ be the intersection point of $\alpha$ and
$\beta_1$ (the dark dot in Figure~\ref{fig:StabilizationPicture}). Let
$\Interval\subset \S(H)$ be the set of $\x\in\S(H)$ which
contain $x_0$. There is, of
course, a natural (point-wise) identification between $\S(G)$ and
$\Interval$, which drops Alexander and Maslov
grading by one.  More precisely, given $\x\in\S(G)$, let $\phi(\x)\in
\S(H)$ denote the induced generator in $\Interval$ which is gotten by
inserting $x_0$. We then have
\begin{align}
M_{C(G)}(\x)&=M_{C(H)}(\phi(\x))+1=M_{C'}(0,\phi(\x))=M_{C'}(\phi(\x),0)+1\\
A_{C(G)}(\x)&=A_{C(H)}(\phi(\x))+g(1)=A_{C'}(0,\phi(\x))=A_{C'}(\phi(\x),0)+g(1)
\end{align}
where $g$ is the function from Section~\ref{sec:chain-complex-cm},
mapping from $i$ to the basis vector corresponding to the component of
the link containing~$O_i$.
With this said, we will
henceforth suppress $\phi$ from the notation, thinking of $\Left$  
and $\Right$ as generated by configurations in $\Interval\subset \S(H)$.

As such, the differentials within $\Left$ and $\Right$ count
rectangles which do not contain $x_0$ on their boundary, although they
may contain $x_0$ in their interior.  Note however that the boundary
operator (in $\Left$ and $\Right$) for rectangles containing $x_0$
does not involve the variable~$U_1$.

\begin{figure}
\begin{center}
\mbox{\vbox{\epsfbox{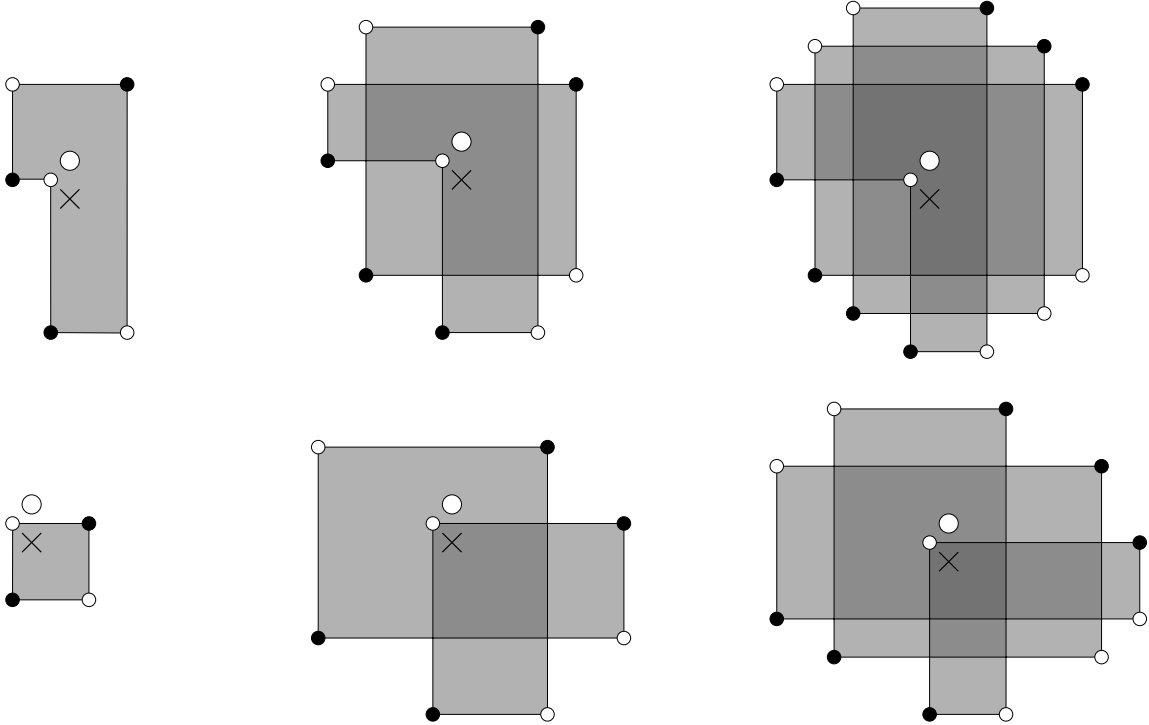}}}
\end{center}
\caption {{\bf Types of domains.}
We have listed here domains in the stabilized diagram, labeling the
initial points by dark circles, and terminal points by empty circles.
The top row lists domains of type $L$, while the second row lists
some of type $R$. The marked $O$ and $X$ are the new ones in the
stabilized picture. Complexities from the left on the first row are
$3$, $5$, and $7$ respectively; on the second, they are $2$, $4$, and
$6$. Darker shading corresponds to
higher local multiplicities.  Not shown is the trivial domain of
type~$L$, which has complexity~$1$.}
\label{fig:Domains}
\end{figure}

\begin{definition}
  For $\x \in \S(H)$ and $\y \in \Interval \subset \S(H)$, a domain
  $p\in \pi(\x,\y)$ is said to be of type $L$ or $R$ if either
  it is trivial, in which case $p$ has type~$L$, or it satisfies the
  following conditions:
  \begin{itemize}
  \item $p$ has only non-negative local multiplicities.
  \item For each $c\in\x\cup\y$, other than $x_0$, at least three of
    the four adjoining squares have vanishing local multiplicities.
  \item In a neighborhood of $x_0$ the local multiplicities in three
    of the adjoining rectangles are the same number~$k$.  When $p$ has
    type $L$, the lower left corner has local multiplicity~$k-1$,
    while for $p$ of type $R$ the lower right corner has
    multiplicity~$k+1$.
  \item $\partial p$ is connected.
  \end{itemize}
  The \emph{complexity} of the trivial domain is~1; the complexity of
  any other domain is the number of horizontal lines in its boundary.
  The set of type $L$ (or $R$) domains from $\x$ to $\y$ is
  denoted $\pL(\x,\y)$ (or $\pR(\x,\y)$).  We set $\pF(\x,\y) =
  \pL(\x,\y) \cup \pR(\x,\y)$, and call its elements domains of {\em
    type $F$}; see Figure~\ref{fig:Domains} for examples. We denote by $\pF$ the union of the sets $\pF(\x,\y),$ over all possible $\x$ and $\y.$

  The \emph{innermost height} (resp.\ \emph{width}) of a domain in
  $\pF$ is the vertical (resp.\ horizontal) distance from the corner
  adjacent horizontally (resp.\ vertically) to $x_0$ to the corner
  after that.
\end{definition}

We now define maps
\begin{align*}
  \FL&\colon C\longrightarrow \Left\\
  \FR&\colon C\longrightarrow \Right
\end{align*}
where $\FL$ (resp.\ $\FR$) counts domains of type $L$ (resp.\ $R$)
without factors of $U_1$.  Specifically, define
\begin{align*}
\FL(\x)&=\sum_{\y\in\S}\,\sum_{p\in\pL(\x,\y)}
U_2^{O_2(p)}\cdots U_n^{O_n(p)}\cm\y \\
\FR(\x)&=\sum_{\y\in\S}\,\sum_{p\in\pR(\x,\y)}
U_2^{O_2(p)}\cdots U_n^{O_n(p)}\cm\y.
\end{align*}

We put these together to define a map
\[
F=
\begin{pmatrix}
\FL \\ \FR 
\end{pmatrix} \colon C \longrightarrow C'.
\]

\begin{lemma}
  \label{lemma:FIsChain}
The map $F\colon C \longrightarrow C'$ preserves
Maslov grading, respects Alexander filtrations, and is a chain map.
\end{lemma}

\begin{proof}
  The fact that the gradings and filtrations are respected 
  is straightforward.
  For instance,
  the Alexander filtration shift of a region~$p$ is given by counting
  the number of $O$'s minus the number of $X$'s contained in~$p$.  A
  region of type $L$ contains $O_1$ and $X_1$ an equal number of
  times, and every other $O_i$ comes with a cancelling factor of
  $U_i$, so the Alexander filtration shift is negative.  The other
  shifts can be checked in a similar way.
  
  To prove that $F$ is a chain map, we consider all the terms in
  the expression $\partial \circ F$ or $F \circ \partial$.  Most of
  these are counts of composite domains $p * r$ or $r * p$, where $r$
  is a rectangle and $p$ is a type $L$ or $R$ domain.  A rectangle
  $r\in\pi(\x,\y)$ cannot contribute to this count if any component $x\in\x$
  is in the interior of $r$, except in the special case where $x=x_0$,
  and the rectangle is thought of as connecting two intersection points in
  $\Left$ or $\Right$, in which case we say it is of Type~2. All other 
  empty rectangles are said to be of Type~1.
  
  There are several cases of domains contributing to $\partial\circ F$
  or $F\circ \partial$ , which we group according to whether $r$ is a
  Type~1 or Type~2 rectangle, and to how many corners $p$ and $r$ have
  in common. We list the cases below; verifying that these are the only
  cases is a straightforward exercise in planar geometry.

If $r$ is of Type~1, we have the following possibilities:

\begin {enumerate}
\item[I(0)] A composition in either order of a domain~$p \in \pF$ and an 
        empty rectangle~$r$ of Type~1, with all corners distinct.  This domain
        appears in both $\partial\circ F$ and $F\circ\partial$ as
        compositions in two different orders, $p * r$ and $r' * p'$, where $r$ 
        has the same support as $r'$ and $p$ has the same support as $p'$. 
\item[I(1)] 
 A composition in either order of a non-trivial domain~$p \in
    \pF$ and an empty rectangle~$r$, with $p$ and $r$ sharing one
    corner and $r$ disjoint from $x_0$ (including the boundary).  The
    union of these two domains has a unique concave corner not at
    $x_0$, and we can slice this into a domain in $\pF$ and a
    rectangle of Type~1 in two ways by cutting in either way from this
    concave corner.  This gives the domain as a composition in exactly
    two ways. An example is shown in Figure~\ref{fig:CaseI1}.

\begin{figure}
\begin{center}
\mbox{\vbox{\epsfbox{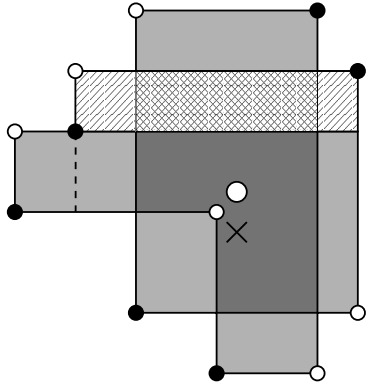}}}
\end{center}
\caption {{\bf Case I(1).}
  An example of a domain with two decompositions $r*p=r'*p'$, both accounted for in case I(1).}
\label{fig:CaseI1}
\end{figure}
\item[I($1'$)] A composition $r * p$ with $r$ and $p$ sharing one corner and  
$x_0$ appearing on
    the horizontal or vertical boundary of~$r$.  The composite looks
    again like a domain in $\pF$ or the rotation by $180^\circ$ of
    such a domain. See Figure~\ref{fig:CaseI1primeII1}. A special case worth 
\begin{figure}
\begin{center}
\mbox{\vbox{\epsfbox{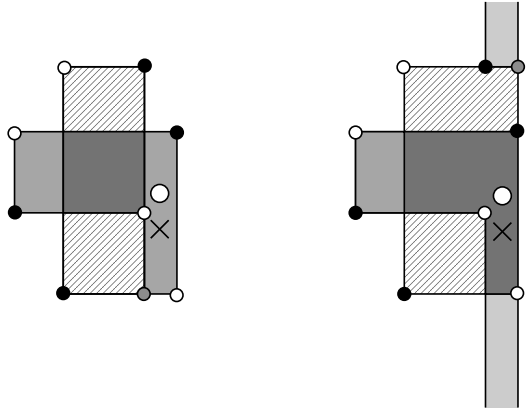}}}
\end{center}
\caption {{\bf Cases I(1$'$) and II(1).} 
There are two terms in $\partial\circ \FL + \FL\circ \partial$ starting at the black
dots and ending at the white dots, thought of as elements of $\Left$. The term on
the left is a juxtaposition $r*p$ (as in I($1'$)), while the second is $p'*r'$, where
$p'$ is of type $L$ and $r'$ is of Type 2 (as in II(1)).}
\label{fig:CaseI1primeII1}
\end{figure}
    mentioning is when $r\in\EmptyRect(\x,\y)$ with $\y\in\Interval$; in 
    this case $p$ is trivial, with complexity~1, as in
    Figure~\ref{fig:CaseI1primeI3simple}.
  \item[I(2)] A composition in either order of $p\in\pF$ and
    $r\in\EmptyRect$, where $p$ and $r$ share two corners other than
    possibly $x_0$, see Figures~\ref{fig:CaseI2II0simple}
    and~\ref{fig:CaseI2II0}.  In this case $p$ has complexity at
    least~3.
\item[I(3)] A domain that wraps around the torus
    with a decomposition as $p * r$ or $r * p$,
    where $r$ is an empty rectangle of Type~1 and $p\in \pF$ has
    innermost height or width equal to~1, and~$r$ and~$p$ share three
    corners other than possibly~$x_0$.  This decomposition is
    unique.  The total domain contains a
    unique vertical or horizontal annulus of height or width equal
    to~1.  When the complexity~$m$ of $p$ is equal to 2, the domain is
    just this annulus.  Examples are shown in
    Figures~\ref{fig:CaseI1primeI3simple} ($m = 3$, horizontal),
\begin{figure}
\begin{center}
\mbox{\vbox{\epsfbox{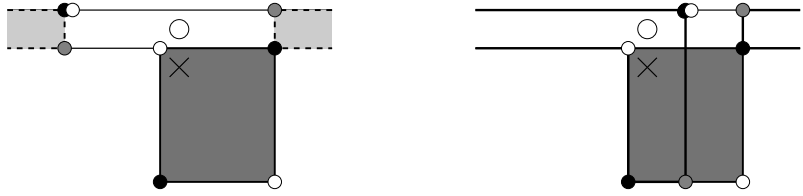}}}
\end{center}
\caption {\textbf{Cases I(1$'$) and I(3), where $m=1$.}
  In both pictures, the
  darkly-shaded rectangle represents a map from the black generator to
  the white one, followed by the natural map (induced by the trivial
  domain, which has complexity 1) to the white generator thought of as
  an element of in~$\Left$.  This is accounted for in
  I($1'$). Depending on the placement of the black dot in the top row,
  we can cancel this either with a term in $\FL\circ \partial$ or
  $\partial\circ \FL$. In the first case (on the left), we have the
  domain $r*p$, where $r$ is the height one
  (lightly-shaded) rectangle in the row through $O_1$, to the
  intermediate generator (labelled by the shaded circle), thought of
  as a differential within $C(H)$, followed by a complexity 3 domain $p$
  with innermost height equal to one, which 
  we trust the reader can spot. In the second case (on the right), we have
  the decomposition $p*r$, where $p$ is the
  complexity $3$ domain with innermost height equal to one from the
  black generator to the intermediate
  generator, which is bounded by the dark line, followed by a
  rectangle to the white generator, which again we leave to the reader to find.
  In both cases the alternate term is accounted for in case~I(3).
\label{fig:CaseI1primeI3simple}}
\end{figure}
    \ref{fig:CaseI1primeI3horiz}~($m = 5$, horizontal),
\begin{figure}
\begin{center}
\mbox{\vbox{\epsfbox{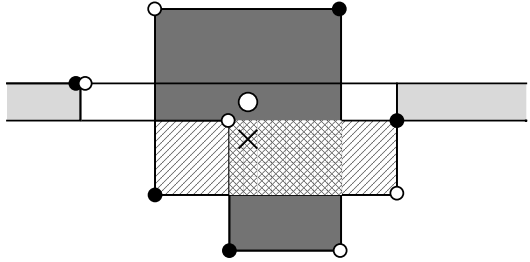}}}
\end{center}
\caption {{\bf Cases I(1$'$) and I(3), horizontal annulus.}
  There are two terms in $\partial \circ \FL+\FL\circ \partial$
  starting at the black dots and ending at the white dots.  One of
  them counts the composite domain $r'*p'$ where $r'$ is the hatched
  rectangle containing $X$, and $p'$ is the darkly-shaded complexity
  $3$ domain (accounted for in I($1'$)); and the other is a count of
  $r*p$, where $r$ is the height one, lightly shaded rectangle,
  followed by a complexity $5$ domain with innermost height equal to
  one (accounted for in I(3)).}
\label{fig:CaseI1primeI3horiz}
\end{figure}
    \ref{fig:CaseI1primeI3vert}~($m = 5$, vertical), and
\begin{figure}
\begin{center}
\mbox{\vbox{\epsfbox{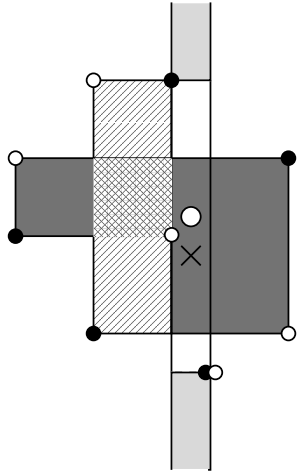}}}
\end{center}
\caption {{\bf Cases I(1$'$) and I(3), vertical annulus.}
  There are two terms in $\partial \circ \FL+\FL\circ \partial$
  starting at the black dots and ending at the white dots.  One of
  them counts the composite domain $r'*p'$ where $r'$ is the hatched
  rectangle containing the white dot $x_0$ in its boundary, and $p'$
  is the darkly-shaded complexity $3$ domain (accounted for in
  I($1'$)); and the other is a count of $r*p$, where $r$ is the height
  one, lightly shaded rectangle, followed by a complexity $5$ domain
  with innermost height equal to one (accounted for in I(3)).}
\label{fig:CaseI1primeI3vert}
\end{figure}
    \ref{fig:typeS}~($m = 4$, horizontal). 
\begin{figure}
\begin{center}
\mbox{\vbox{\epsfbox{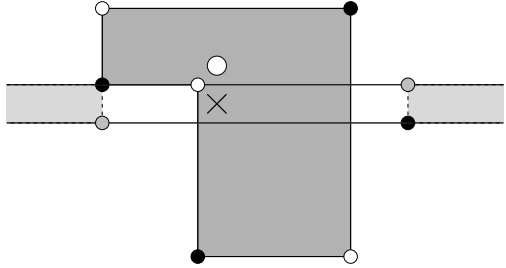}}}
\end{center}
\caption {\textbf{Cases I(1$'$), I(3), and (S).}
  This case is similar to those in
  Figures~\ref{fig:CaseI1primeI3horiz}
  and~\ref{fig:CaseI1primeI3vert}, except that it also involves a
  domain of type~(S).
  We count terms in $\partial \circ F+F\circ \partial$ starting at the
  black dots and ending at the white dots (thought of as representing
  an element of $\Right$).  The darkly-shaded polygon represents a
  domain of type $L$ from the black to the white generator.
  Post-composing with the differential from $\Left$ to $\Right$, we
  get $(U_2-U_1)$ times the white generator. Alternatively, the region
  can be decomposed as a rectangle containing $O_1$ (with a factor of
  $U_1$), composed with the rectangle
  containing~$X_1$, thought of as a polygon of type~$R$.
  Alternatively, there is a term induced by the height one
  (lightly-shaded) rectangle, followed by a complexity 4 domain of
  type~$R$, which the reader can easily spot.  One of these two
  domains contains~$O_2$, and hence the composite will count with a
  factor of $U_2$.}
\label{fig:typeS}
\end{figure}
\end {enumerate}

If $r$ is of Type 2, the composition must be of the form $p * r$, because 
Type 2 rectangles only appear in the differential of the target complex 
$C'$. We only have two possibilities:

\begin {enumerate}
\item[II(0)]  All the corners of $p$ and $r$ are disjoint. 
\item[II(1)]  A domain that wraps around the torus
    with a decomposition as $p * r$, where $r$ is a rectangle of
    Type~2 that shares one corner with~$p$.  This decomposition is
    unique, and the total domain again contains a unique thin (i.e., width one or height one) 
    annulus. See 
Figure~\ref{fig:CaseI1primeII1}.

\end {enumerate}

Apart from these, there is one other special contribution to $F \circ
\partial$, which does not come from a decomposition of a domain into
$p*r$ or $r'*p'$:
\begin {enumerate}
\item[(S)]
A domain $p \in \pL$ followed by the differential from
$\Left$ to $\Right$, which multiplies by $U_2-U_1$.
\end {enumerate}

Contributions from case I(0) cancel each other out, and the same goes for 
those from case I(1). In fact, these cases are the exact analogs of the 
first three cases in Figure~\ref{fig:DSquaredZero} for the proof of 
Proposition~\ref{prop:DSquaredZero}.  See Figure~\ref{fig:CaseI1} for
an example.

We claim that contributions from case I($1'$) cancel with 
contributions from case~II(1) or~I(3), together with possibly a
contribution from case~(S). Indeed, for each domain of 
type~I($1'$) made of a rectangle~$r_1 \in \EmptyRect(\x,\y)$ and a domain 
$p_1 \in \pF(\y,\z)$ of complexity~$m$, let $p_0 = r_1 * p_1$.  We can make a
new domain~$p_0'$ by adding a thin annulus abutting~$x_0$ on the
opposite side of $x_0$ from~$r_1$.  (For instance, if the right side
of $r$ touches $x_0$, add a vertical annulus of width one whose left
side touches $x_0$.)  In the case when $m=1$, when $r_1$ touches
$x_0$ at a corner, we attach a horizontal annulus if $r_1$ contains
$X_1$ and a vertical annulus otherwise, as in
Figure~\ref{fig:CaseI1primeI3simple}. If the innermost height or width
of $p_1$ is~1, then
$p_0$ decomposes as $p_2 * r_2$, where $p_2\in\pF$ has
complexity~$m$. This corresponds to a contribution from case~II(1), as
in Figure~\ref{fig:CaseI1primeII1}.
If, on the other hand, the innermost height or
width (as appropriate) of $p_1$ is not 1, the new domain~$p_0'$ is
of type I(3) and in turn decomposes as $p_2 *
r_2$ or $r_2 * p_2$, depending on the placement of the generator on
the new row or column, where $p_2\in\pF$ has complexity $m+2$.
See Figures~\ref{fig:CaseI1primeI3simple}--\ref{fig:typeS}.

  In these cases involving annuli, if $p_i \in \pR$ and the annulus is
  horizontal, the rectangle~$r_1$
  contains~$O_1$ and so has a contribution which is multiplied by
  $U_1$, while the domain $p_0'$ contains $O_2$ and so has a
  contribution which is multiplied by $U_2$.  Thus these two terms
  contribute $U_1 - U_2$ to the composite map from $\x$ to $\z$.  On
  the other hand, in this case the domain $p_0$ is itself in
  $\pL(\x,\z)$, and so we get a cancelling contribution of
  type~(S), as in Figure~\ref{fig:typeS}.
  In other cases the two domains~$p_0$ and~$p_0'$ give the
  same contribution to the boundary map.

Compositions $r*p$ or $p*r$ 
from case I(2), with $p$ of complexity $m \geq 3$, cancel out 
compositions $r' * p'$ from case II(0), with $p'$ of complexity $m-2$, 
as illustrated in Figures~\ref{fig:CaseI2II0simple} ($m = 3$) 
and~\ref{fig:CaseI2II0} ($m = 5$).
\begin{figure}
\begin{center}
\mbox{\vbox{\epsfbox{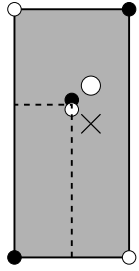}}}
\end{center}
\caption {{\bf
Cases I(2) and II(0), with complexity $m=3$.}
The simplest case of the pairing between cases I(2) and II(0).}
\label{fig:CaseI2II0simple}
\end{figure}
\begin{figure}
\begin{center}
\mbox{\vbox{\epsfbox{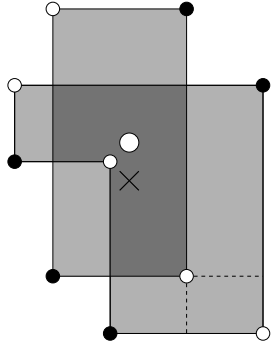}}}
\end{center}
\caption {{\bf Cases I(2) and II(0).}
The illustrated domain can be decomposed as a complexity $3$ domain of type $L$ followed
by a Type 2 rectangle (accounted for in II(0)), or alternatively a complexity $5$ domain of type
$L$ followed by a Type 1 rectangle (accounted for in I(2)).
\label{fig:CaseI2II0}}
\end{figure}

The only domains left to cancel are those of type~I(3) with $m=2$ and
type~(S) with $m=1$. There are two kinds of domains of 
type~I(3) with $m=2$: a vertical and a horizontal annulus,
containing~$O_1$ and~$O_2$, respectively, and in both cases
containing~$X_1$. These
domains map a generator $\x\in\Interval$
to itself, and so cancel
the remaining contribution from 
the maps of type~(S).
\end{proof}

In order to see that $F$ is a quasi-isomorphism, we will introduce an
appropriate filtration.  Consider $\Caa(H)$.  Let $Q$ be a collection
of $(n-1)^2$ dots, one placed in each square which do not appear in
the row or column through $O_1$.  Given $h\in\halfzl$, let
$\Caa(H,h)$ denote the summand generated by generators $\x$ with
Alexander gradings equal to $h$.

Note that for fixed $\x,\y\in\S_h$, for any two domains $p,
p'\in\pi(\x,\y)$ with $O_i(p)=X_i(p)=O_i(p')=X_i(p')$ for all~$i$,
we have that $\#(Q\cap p) = \#(Q\cap p')$. Thus, we can find 
a function~$\Filt$ so that for any
$\x,\y\in\S$, if $p\in\pi(\x,\y)$ is a domain with
$O_{i}(p)=X_{i}(p)=0$ for all $i$, then
$$\Filt(\x)-\Filt(\y) = \#(Q\cap p).$$
The function~$\Filt$
determines a filtration on $\Caa(H,h)$, whose associated
graded object counts only those rectangles which contain no $O_i$,
$X_i$, or points in $Q$. Thus, these rectangles must be supported in
the row or column through $O_1$. We let $\Caa_Q$ denote this associated
graded object, and typically drop $h$ from the notation.

We recall now a well-known principle from homological algebra (see for example
Theorem~3.2 of~\cite{UsersGuide}).

\begin{lemma}
        \label{lemma:FilteredIsomorphism}
  Suppose that $F \colon C\longrightarrow C'$ is a filtered chain map
  which induces an isomorphism on the homology of the associated
  graded object. Then $F$ is a filtered quasi-isomorphism.
\end{lemma}

We decompose $\S=\Interval \cup (\NI)\cup (\NN)$, where
$\NI$ consists those configurations whose $\beta_2$~component is
$\alpha\cap\beta_2$ and whose $\beta_1$~component is not in $\alpha$,
while $\NN$ consists of those whose $\beta_2$~component and
$\beta_1$~component are not on~$\alpha$.
We have corresponding decompositions of modules:
$C=\CI\oplus\CNI\oplus\CNN$.

\begin{lemma}
  \label{lemma:Splittings}
  $H_*(\CQ)$ is isomorphic to the free $\Field$-module generated
  by elements of $\Interval$ and $\NI$.
\end{lemma}

\begin{proof}
   There are two cases, according to whether the $X_2$ marks the square
   to the left or the right of $O_1$.

 Suppose $X_2$ is in the square just to the right of the square marked
 $O_1$.  Then we have a direct sum splitting $\CQ\cong \CNIQ\oplus B$,
 where the differentials in $\CNIQ$ are trivial, hence its homology
 is the free $\Field$-module generated by elements of $\NI$; and
 where $B$ is a chain complex fitting into an exact sequence
 \[
 \begin{CD}
   0@>>>\CIQ@>>>B@>>>\CNNQ@>>>0.
 \end{CD}
 \]
 Moreover, it is easy to see that $H_*(\CNNQ)=0$.
 Finally, the differentials in $\CIQ$ are trivial, so its homology
 is the free $\Field$-module generated by elements of~$\Interval$.

 Suppose on the other hand that $X_2$ is just to the left of $O_1$.
 Then there is a direct sum splitting $\CQ\cong\CIQ\oplus B'$,
 where once again the differentials on $\CIQ$ are trivial
 and $B'$ fits into an exact sequence
 \[
 \begin{CD}
   0@>>>\CNNQ@>>>B'@>>>\CNIQ@>>>0,
 \end{CD}
 \]
 where $H_*(\CNNQ)=0$ and the differentials on $\CNIQ$ are trivial.
\end{proof}

\begin{prop}
  \label{prop:Stabilization}
  The map $F$ is a filtered quasi-isomorphism.
\end{prop}

\begin{proof}
  We consider the map induced by $F$:
  $$\FaQ \colon \Caa_Q \longrightarrow \Caa_Q'.$$
  $\Caa_Q'$ splits as a direct sum of chain
  complexes $\Left_Q\oplus \Right_Q$, both of which are freely
  generated by elements in $\Interval$.

   There are two cases.
   First take the case where $X_2$ is in the square just to the right
   of the square marked $O_1$.
   Consider the subcomplex $\CIQ\oplus\CNIQ\subset \CQ$. By
   Lemma~\ref{lemma:Splittings}, this subcomplex
   carries the homology, and hence it suffices to show that the restriction
   of $\FaQ$ to this subcomplex induces an isomorphism in homology.
  
   To this end observe that $\FLaQ$ restricted to $\CIQ$ is an
   isomorphism.  Moreover, $\FRaQ$ restricted to $\CNIQ$ counts
   rectangles supported in the row and column through $O_1$ and which
   contain $X_1$ in their interior and end up in $\Interval$ (since no
   other domains of type~$R$ is disjoint from $Q$). But for each
   element of $\NI$, there is a unique such rectangle. Thus $\FaQ$ is
   a quasi-isomorphism when $X_2$ is just to the right of $O_1$.
  
   In the second case, where $X_2$ is just to the left of $O_1$, we
   proceed as follows. In this case $\CIQ$ is a direct summand of the
   complex $\CQ$ (cf.\ the proof of Lemma~\ref{lemma:Splittings}).
   Moreover, it is easy to see that $\FLaQ$ restricted to $\CIQ$ is an
   isomorphism of chain complexes. It remains to show that the
   restriction of $\FRaQ$ is a quasi-isomorphism.  This is true
   because the only domains of type $R$ which do not contain $X_2$
   are rectangles, and those which are supported in the allowed region
   connect configurations of type $\NI$ to $\Interval$. Once again,
   the result now follows from the fact that there is a unique
   rectangle of type $R$ connecting a given element of $\NI$ to an
   element of $\Interval$. This completes the verification that $\FaQ$
   is a quasi-isomorphism.

  We now appeal to Lemma~\ref{lemma:FilteredIsomorphism} to conclude
  that ${\widehat F}$ is quasi-isomorphism; and another application of
  the same principle gives that $F$ is a quasi-isomorphism, as well.
\end{proof}

\begin{remark}
  The chain complex~$C'$ used in this stabilization proof can be
  viewed as the chain complex associated to the Heegaard diagram where
  the vertical circle~$\beta_1$ is replaced by a small circle
  enclosing~$O_1$ and~$X_1$.  In this Heegaard diagram it is
  straightforward to check that the counts of holomorphic disks are
  still combinatorial and equivalent to the boundary operator in~$C'$.
\end{remark}

\subsection{Completion of topological invariance, without signs}

We have now all the pieces needed to establish Theorem~\ref{thm:WithSigns},
with coefficients in $\Field=\Z/2\Z$.

\begin{proof}[Proof of Theorem~\ref{thm:WithSigns}]
This result now is an immediate consequence of Cromwell's theorem, our
earlier remarks on cyclic permutation, and
Propositions~\ref{prop:Stabilization} and~\ref{prop:Commute}. 
\end{proof}


\section{Signs}
\label{sec:Signs}

\begin{definition}
  \label{def:SignAssignment}
  A {\em true sign assignment}, or simply a {\em sign assignment},
  is a function
  $$\sign\colon \EmptyRect\longrightarrow \{\pm 1\}$$
  with the following properties:
  \begin{enumerate}
  \item[(Sq)] For any four distinct $r_1,r_2,r_1',r_2'\in\EmptyRect$
    with $r_1*r_2=r_1'*r_2'$, we have that
    $$\sign(r_1)\cm\sign(r_2)=-\sign(r'_1)\cm\sign(r'_2).$$
  \item[(V)] If $r_1,r_2\in\EmptyRect$ have the property that
    $r_1*r_2$ is a vertical annulus, then 
    $$\sign(r_1)\cm\sign(r_2)=-1.$$
  \item[(H)]
    If $r_1,r_2\in\EmptyRect$ have the property that
    $r_1*r_2$ is a horizontal annulus, then 
    $$\sign(r_1)\cm\sign(r_2)=+1.$$
  \end{enumerate}
\end{definition}

\begin{theorem}
  \label{thm:SignAssignments}
  There is a sign assignment in the sense of
  Definition~\ref{def:SignAssignment}.  Moreover, this sign assignment
  is essentially unique: if~$\sign_1$
  and~$\sign_2$ are two sign assignments, then there is a function
  $f\colon \S\longrightarrow \{\pm 1\}$ so that for all
  $r\in\EmptyRect(\x,\y)$,
  $\sign_1(r)=f(\x)\cm f(\y)\cm \sign_2(r)$.
\end{theorem}

We turn to the proof of this theorem in
Subsection~\ref{subsec:SignExists}.  We can use the sign assignment
from Theorem~\ref{thm:SignAssignments} to construct the chain complex
over~$\Z$ as follows.  Fix a true sign assignment $\sign$. Define
$\Cm(G)$ to be the free $\Z[U_1,\ldots,U_n]$-module generated by
$\x\in\S(G)$, endowed with Maslov grading and Alexander filtration as
before.  We endow this with the endomorphism
\begin{align*}
\dm_\sign&: \Cm(G) \longrightarrow \Cm(G)\\
\dm_\sign(\x)&=\sum_{\y\in\S}\,
    \sum_{r\in\EmptyRect(\x,\y)}\!\!
\sign(r) \cm U_1^{O_1(r)}\cdots U_n^{O_n(r)}\cm \y.
\end{align*}

We will check that this endomorphism gives the sign refinement of $\Cm(G)$
needed in Theorem~\ref{thm:WithSigns}. In turn, the proof of that
theorem involves reexamining the invariance proof from
Section~\ref{sec:Invariance}, and constructing sign refinements for
the chain maps and homotopies used there. We turn to this task in
Subsection~\ref{subsec:DiscussionOverZ}.  However, first we construct
the sign assignments, proving
Theorem~\ref{thm:SignAssignments}.

\subsection{The existence and uniqueness of sign assignments}
\label{subsec:SignExists}

\begin{definition}
        A {\em thin rectangle} is a rectangle with width one.  We
        denote the set of thin rectangles $\tRect$; given
        $\x,\y\in\S$, we let
        $\tRect(\x,\y)=\tRect\cap\Rect(\x,\y)$.
        For fixed $\x$ and~$\y$ and $n > 2$, there can be at most one
        element in $\tRect(\x,\y)$.
\end{definition}

Sign assignments as in Theorem~\ref{thm:SignAssignments}
are constructed in the following six steps.

\begin{enumerate}
\item 
Define sign assignments in a more restricted sense, 
{\em sign assignments for the Cayley graph}. These are analogues of sign assignments
defined only for thin rectangles supported
in an $(n-1)\times (n-1)$ subsquare of the torus,
satisfying a suitable restriction of Property~(Sq) from Definition~\ref{def:SignAssignment}.
\item
Show that sign assignments for the Cayley graph satisfy a uniqueness property.
Establish existence by giving an explicit formula.  (It is also possible to give a more abstract existence
argument, but the formula is needed in the next step.)
\item
Extend the formula to include all thin rectangles on the torus,
and show that it satisfies, once again, axioms gotten by restricting
Properties~(Sq) and~(V) to thin rectangles.
\item 
Show that a sign assignment on thin rectangles can be extended
uniquely to a function
satisfying Properties~(Sq) and~(V),
but not necessarily Property~(H).
\item
Given the sign assignment on thin rectangles chosen in Step 3, establish 
a formula for the values of the function from Step 4 on empty rectangles 
supported in the $(n-1)\times (n-1)$ subsquare of the torus.
\item
With our choices of signs, use the explicit formulas from Step~5 to
show that the function from Step~4 satisfies
Property~(H), thus giving a sign assignment in the 
sense of Definition~\ref{def:SignAssignment}.
\end{enumerate}

\step{1}{Define sign assignments on the Cayley graph}

We denote by $\Square = [0,n-1] \times [0,n-1]$ the $(n-1)\times 
(n-1)$-subsquare of the torus with the lower left corner at the origin.

\vskip.2cm
\begin{definition}
        Given $\x,\y\in\S$, a {\em thin rectangle in $\Square$ from
        $\x$ to $\y$} is a rectangle $r\in\tRect(\x,\y)$ supported
        inside $\Square$.  A thin rectangle in $\Square$ {\em connects
        $\x$ and $\y$} if it is a thin rectangle in $\Square$ from
        $\x$ to $\y$ or from $\y$ to $\x$. The set of all thin
        rectangles in $\Square$ is denoted $\stRect$.
\end{definition}

The set $\stRect$ has an interpretation in terms
of a Cayley graph of the symmetric group in the following sense.

Consider the graph $\Cayley$ whose vertices are elements in the
symmetric group on $n$ letters ${\mathfrak S}_n$, and whose edges are
labeled by the $n-1$ adjacent transpositions $\{\tau_i\}_{i=1}^{n-1}$
in ${\mathfrak S}_n$, with an edge labelled $\tau_i$ connecting
$\sigma_1,\sigma_2\in {\mathfrak S}_n$ precisely when
$\sigma_2=\sigma_1\cm \tau_i$.  When $\sigma_2=\sigma_1\cm \tau_i$, we
join $\sigma_1$ and $\sigma_2$ by exactly one edge, i.e. we do not
draw an additional one for the relation $\sigma_1=\sigma_2\cm
\tau_i$. $\Cayley$ is the Cayley graph of ${\mathfrak S}_n$
with respect to the generators $\tau_i$.

There is a one-to-one correspondence between elements in ${\mathfrak
S}_n$ and generators~$\S$, which is obtained by viewing elements of~$\S$
as graphs of permutations~$\sigma_\x$. (To this end,
we think of ${\mathfrak S}_n$ as permutations of the letters $\{0,\ldots,n-1\}$.)
This can be extended to a
one-to-one correspondence between edges in the Cayley graph $\Cayley$ and 
elements of $\stRect$, sending a rectangle $r\in\stRect$ which connects $\x$
and $\y$ to the corresponding edge in the Cayley graph connecting
$\sigma_\x$ and $\sigma_\y$:
\begin{equation}
        \label{eq:stRectCayley}
 \stRect\cong \Edges(\Cayley).
\end{equation}

\begin{definition}
  \label{def:SignAssignmentSn}
  A \emph{sign assignment on the Cayley graph} is a function
        \[
                \sign_0\colon \Edges(\Cayley) \longrightarrow \{\pm 1\} 
        \]
  with the following properties:
  \begin{enumerate}
  \item[(Sq)] If $\{e_1,\ldots,e_4\}$ are four
    edges which form a square, then 
    $$\sign_0(e_1)\cdots\sign_0(e_4)=-1.$$
  \item[(Hex)] If $\{e_1,\ldots,e_6\}$  are six edges
    which form a hexagon, then
    $$\sign_0(e_1)\cdots\sign_0(e_6)=1.$$
  \end{enumerate}
\end{definition}

\begin{figure}
\begin{center}
\input{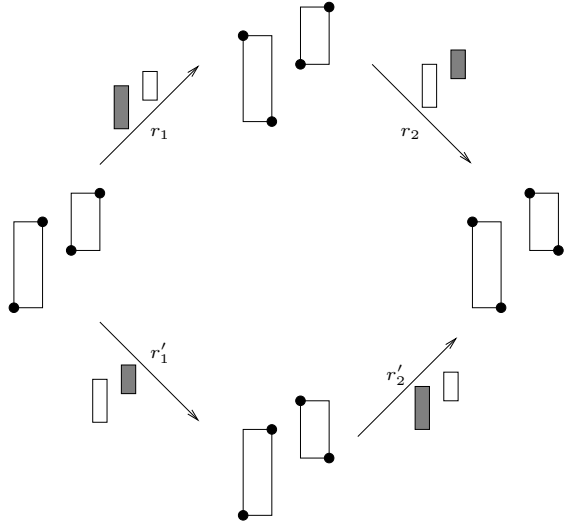}
\end{center}
\caption {{\bf The rectangles in the square rule.}}
\label{fig:CayleySq}
\end{figure}

Note that a square in the Cayley graph corresponds to two pairs of
disjoint rectangles $r_i\in\tRect(\x_i,\x_{i+1})$ and
$r_i'\in\tRect(\x_i',\x_{i+1}')$ for $i=1,2$, with $\x_1=\x_1'$ and
$\x_2=\x_2'$, as pictured in Figure~\ref{fig:CayleySq}. Similarly, a
hexagon in the Cayley graph corresponds to six thin rectangles
$r_i\in\tRect(\x_i,\x_{i+1})$ and $r'_i\in\tRect(\x'_i,\x'_{i+1})$ for
$i=1,2,3$ with $\x_1=\x'_1$ and $\x_3=\x'_3$, such that the union of
the support of $r_1$, $r_2$, and $r_3$ is a rectangle (with width
two), as pictured in Figure~\ref{fig:CayleyHex}.

\begin{figure}
\begin{center}
\input{CayleyHex.pstex_t}
\end{center}
\caption {{\bf The rectangles in the hexagon rule.} Two applications of the square rule (Sq) in Definition~\ref{def:SignAssignment} give $\sign(r_1)\cdot \sign(r_2) = \sign(r_1') \cdot \sign(r_4)$ and  $\sign(r_4)\cdot \sign(r_3) = \sign(r_2') \cdot \sign(r_3')$. These imply the hexagon rule $\sign(r_1) \cdot \sign(r_2) \cdot \sign(r_3) = \sign(r_1') \cdot \sign(r_2') \cdot \sign(r_3')$. }
\label{fig:CayleyHex}
\end{figure}

We can relate the above notion of a sign assignment to the earlier
notion of sign assignments (Definition~\ref{def:SignAssignment}):

\begin{lemma}
  \label{lemma:CayleyGraphSignChar} The restriction of a sign
  assignment in the sense of Definition~\ref{def:SignAssignment} to
  the Cayley graph of~${\mathfrak S}_n$ is a sign assignment on the
  Cayley graph as defined in Definition~\ref{def:SignAssignmentSn}.
\end{lemma}

\begin{proof}
  The first property follows from the corresponding property in
  Definition~\ref{def:SignAssignment}.  For the second property, note
  that there is a rectangle $r_4$ of width~2 that cuts across a diagonal of
  the hexagon, cf. Figure~\ref{fig:CayleyHex}.  Two applications of
  Property~(Sq) (both involving $r_4$)
  now shows that the total number of sign changes around
  the hexagon must be even.
\end{proof}

\step{2}{Signs assignments on the Cayley graph exist and are unique}
\begin{prop}
  \label{prop:sign-cayley}
  A sign assignment on the Cayley graph
  exists and is unique up to equivalence given by changing the sign of
  the basis elements (as in Theorem~\ref{thm:SignAssignments}).
\end{prop}

\begin{proof}
  Recall that for $A, B \subset \R^2$, we defined $\NE(A,B)$ to be
  the number of pairs $(a_1,a_2)\in A$ and $(b_1,b_2)\in B$ with
  $a_1<b_1$ and $a_2<b_2$.

  Given an edge of the Cayley graph, let $r\in\stRect(\x,\y)$ denote
  its corresponding rectangle. Let $h=h(r)$ denote the height of the top
  edge of $r$ (i.e., the four corners of $r$ are $(i,a)$, $(i,h)$,
  $(i+1,a)$ and $(i+1,h)$, where $(i,a)$ and $(i+1,h)$ belong to $\x$
  and $(i,h)$ and $(i+1,a)$ belong to $\y$).

  We then define
  \begin{equation}
        \label{eq:ThinSignAssignment}
        \sign(r)=(-1)^{\NE(\x,\{(x_1,x_2)\in\x \mid x_2\leq h(r)\})}.
  \end{equation}

  We check that each square anti-commutes. To this end, observe that a
  square in the Cayley graph corresponds to four rectangles
  $r_1\in\stRect(\x,\y)$, $r_2\in\stRect(\y,\z)$,
  $r_1'\in\stRect(\x,\y')$ and
  $r_2'\in\stRect(\y',\z)$, where $r_1$ and $r_2$ have distinct corners 
  and $r_1 * r_2 = r_1' * r_2'$ as in Figure~\ref{fig:CayleySq}.
  Number the rectangles so that $h(r_1) = h(r_2') <
  h(r_2) = h(r_1')$. It is easy to see that
  \begin{align*}
        \sign(r_1)=\sign(r_2'), \\
        \sign(r_2)=-\sign(r_1').
  \end{align*}

  We can similarly check that a hexagon commutes.  Consider the six
  thin rectangles corresponding to a a hexagon in the Cayley graph
  $r_i\in\tRect(\x_i,\x_{i+1})$ and $r'_i\in\tRect(\x'_i,\x'_{i+1})$
  for $i=1,2,3$ with $\x_1=\x'_1$ and $\x_3=\x'_3$, as pictured in
  Figure~\ref{fig:CayleyHex}. One can check that \begin{align*}
  \sign(r_1)&=\sign(r_3') \\ \sign(r_2)&=\sign(r_2') \\
  \sign(r_3)&=\sign(r_1'), \end{align*} so that in particular
  $$\sign(r_1)\sign(r_2)\sign(r_3)=\sign(r_1')\sign(r_2')\sign(r_3'),$$
  as needed.

  We will now prove uniqueness.  Let $\sign$ and $\sign'$ be two
  sign assignments on~${\mathfrak S}_n$.  Define a new function~$T$ on the Cayley
  graph by
  \[
  T(\x;\tau_i) = \sign(\x;\tau_i)\sign'(\x;\tau_i)
  \]
  Then the product of~$T$ around any square or hexagon is equal
  to~$1$.

  Let~$W$ be the Cayley complex of~${\mathfrak S}_n$: the 2-complex whose
  edges and vertices form the Cayley graph of~${\mathfrak S}_n$, and whose 2-cells
  are the squares connecting $\{\x, \x\tau_i,
  \x\tau_i\tau_j, \x\tau_j\}$ for $|i-j| > 1$ and the
  hexagons connecting $\{\x, \x\tau_i, \x\tau_i\tau_{i+1},
  \x\tau_i\tau_{i+1}\tau_i, \x\tau_{i+1}\tau_i,
  \x\tau_{i+1}\}$.  Since these squares and hexagons (together with the relations $\tau_i^2 = 1$, which are suppressed in the definition of $\Cayley$) form a
  complete set of relations for ${\mathfrak S}_n$, the complex~$W$ is 
  simply connected.

  Now consider~$T$ as an element of $C^1(W; \{\pm 1\})$.  The
  conditions on~$T$ are equivalent to saying that it is a cocycle:
  $\delta T = 0$.  Since~$W$ is simply connected, there is therefore a
  function $f \in C^0(W;\{\pm 1\})$ so that $\delta f = T$.  This
  function~$f$ gives the desired choice of signs on the basis.
\end{proof}

\begin{remark}
  We could prove Proposition~\ref{prop:sign-cayley} without explicitly
  exhibiting the sign assignment: In general, suppose we have a
  2-complex~$W$ and are looking for an assignment of $\pm1$ to the
  edges of the 2-complex so that the number of $-1$ signs is odd
  around a prescribed set of 2-cells.  Such an assignment is unique
  (if it exists) iff~$H^1(W;\{\pm 1\})$ is trivial, as in the proof of the
  Proposition.  Furthermore, such an assignment
  exists if there is a 3-complex $W'$ with $W$ as
  its 2-skeleton so that $H^2(W';\{\pm1\}) = 0$ and the set of faces
  with an odd number of $-1$~signs, considered as a 2-cocycle on $W'$,
  is coclosed.  In the
  case at hand, we can take $W'$ to be the 3-skeleton of the
  \emph{permutahedron}~\cite{Permutahedron}, which can be defined as
  the convex hull of the vectors obtained by permuting the coordinates
  of $(1,2,\ldots,n)$.  This 3-skeleton is $W$ with the following
  types of 3-cells attached:
  \begin{itemize}
  \item Cubes corresponding to 3 disjoint transpositions, an ${\mathfrak S}_2
    \times {\mathfrak S}_2 \times {\mathfrak S}_2 \subset {\mathfrak S}_n$;
  \item Hexagonal prisms corresponding to ${\mathfrak S}_3 \times {\mathfrak S}_2 \subset
    {\mathfrak S}_n$; and
  \item Truncated octahedra corresponding to ${\mathfrak S}_4 \subset {\mathfrak S}_n$.
  \end{itemize}
  (For the last case, note that the Cayley graph of~${\mathfrak S}_4$ is the
  boundary of a truncated octahedron.)  In each case the number of
  squares on the boundary of the 3-cell is even, so an assignment of
  signs exists.  The permutahedron is convex (hence contractible),
  so $W'$ is 2-connected.
\end{remark}

\step{3}{Extend sign assignments to all thin rectangles in the torus}

\begin{definition}
        \label{def:VerticalSignAssignmentThin}
        A {\em vertical sign assignment for thin rectangles} is a
        function $$\sign\colon \tRect\longrightarrow \{\pm 1\},$$
        which satisfies the following properties:
        \begin{enumerate} 
        \item[(Sq)]
        Given thin rectangles $r_1\in\tRect(\x,\y)$ and
        $r_2\in\tRect(\y,\z)$ with distinct corners, if we let
        $r_1'\in\tRect(\x,\y')$ and $r_2'\in\tRect(\y',\z)$ be two other 
        rectangles such that $r_1 * r_2 = r_1' * r_2'$, we have that
        $$\sign(r_1)\sign(r_2)=-\sign(r_1')\sign(r_2').$$
        (See Figure~\ref{fig:CayleySq}.)

        \item[(Hex)]
          Given six thin rectangles $r_i\in\tRect(\x_i,\x_{i+1})$ 
        and $r'_i\in\tRect(\x'_i,\x'_{i+1})$ for $i=1,2,3$ with
        $\x_1=\x'_1$ and $\x_3=\x'_3$, such that the union of the support
        of $r_1$, $r_2$, and $r_3$ is a rectangle (with width two), 
        we have that 
        
        $$\sign(r_1)\sign(r_2)\sign(r_3)=\sign(r_1')\sign(r_2')\sign(r_3').$$
        (See Figure~\ref{fig:CayleyHex}.)

        \item [(V)]
              If $r_1\in\tRect(\x,\y)$ and $r_2\in\tRect(\y,\x)$, then 
        $\sign(r_1)=-\sign(r_2)$.
        \end{enumerate}
\end{definition}

\begin{prop}
        \label{prop:VerticalSignAssignmentThin}
        There is a vertical sign assignment for thin rectangles.
\end{prop}

\begin{proof}
        We extend Equation~\eqref{eq:ThinSignAssignment}, as follows.

        Note that $\Torus$ is obtained from $\Square$ by adding one more
        row of squares, which are 
        of the form 
        $$\{[i,i+1]\times [n-1,n]\}_{i=0}^{n-1},$$ 
        and one more column of squares which are of the form 
        $$\{[n-1,n]\times [j,j+1]\}_{j=0}^{n-1}.$$

        Consider a thin rectangle $r$ in the torus.  If $r$ is
        contained in $\Square\subset \Torus$, then $\sign(r)$ is as in
        Equation~\eqref{eq:ThinSignAssignment}. If $r\in\tRect(\x,\y)$
        is a thin rectangle which is supported in the new column, but
        which is disjoint from the new row, so that it is of the form
        $[n-1,n]\times [a,b]$ with $0\leq a<b<n$ we define
        \begin{equation} \label{eq:LastColumn}
        \sign(r)=(-1)^{\NE\big(\x,\{(x_1,x_2)\in\x| x_2 \leq a\}\big)
        + \NE\big(\x,\{(x_1,x_2)\in\x| a<x_2<b \text{ and }
        x_2~\text{even}\}\big)+b}.  \end{equation} The thin rectangles
        not covered by the above two cases are those whose interiors
        meet the new row of squares. It is easy to see that for each
        such rectangle $r\in\tRect(\x,\y)$, there is a unique other
        thin rectangle $r'\in\tRect(\y,\x)$, whose interior does not
        meet the row of squares, and hence whose sign is defined
        either by Equation~\eqref{eq:ThinSignAssignment} or by
        Equation~\eqref{eq:LastColumn}. We then define
        $$\sign(r)=-\sign(r').$$

This ensures that Property~(V) in 
Definition~\ref{def:VerticalSignAssignmentThin} is satisfied. We must now 
verify Properties~(Sq) and~(Hex).

If all the thin rectangles involved are contained in $\Sigma$, then the 
two conditions were already checked in the proof of 
Proposition~\ref{prop:sign-cayley}. Let us consider the cases when all the 
rectangles involved are disjoint from the new row, but at least one of 
them is supported in the new column, so that its sign is given by 
Equation~\eqref{eq:LastColumn}.

Let us consider the square rule, with the support of $r_1$ (which is the 
same as the support of $r_2'$) being in the last column, and the support 
of $r_2$ (the same as that of $r_1'$) in $\Sigma$. Let $h(r_1)$ and 
$h(r_2)$ be the heights of the top edges of $r_1$ and $r_2$, respectively. 
If $h(r_1) \leq h(r_2)$, then, just as in the proof of 
Proposition~\ref{prop:sign-cayley}, we have: 
$$  \sign(r_1)=\sign(r_2'), \  \sign(r_2)=-\sign(r_1').$$

If $h(r_1) \geq h= h(r_2)$, cf.\ Figure~\ref{fig:WidthOne}, let $r_1 = 
[n-1,n]\times [a,b]$. If $h < a$, we obtain:
$$ \sign(r_1) = -\sign(r_2'), \ \sign(r_2)= \sign(r_1').$$
If $h \in (a,b)$ as in Figure~\ref{fig:WidthOne}
then, in comparing $\sign(r_2)$ and $\sign(r_1')$ using 
Equation~\eqref{eq:ThinSignAssignment}, there is a discrepancy coming from 
$a$ pairs of points where the second 
point has coordinates $(n-1, a)$, and $h-a$ pairs of points where the 
first point has coordinates $(0,a)$. Therefore,
$$ \sign(r_2) = (-1)^h\sign(r_1').$$

On the other hand, in comparing $\sign(r_1)$ and $\sign(r_2')$ using
Equation~\eqref{eq:LastColumn}, there can only be a discrepancy of one 
extra pair (two corners of $r_2$), which appears in case $h$ is even.  
Thus,
$$ \sign(r_1) = (-1)^{h+1}\sign(r_2'),$$
which implies that $\sign(r_1)\sign(r_2)=-\sign(r_1')\sign(r_2')$, as 
desired.

\begin {figure}
\begin {center}
\input {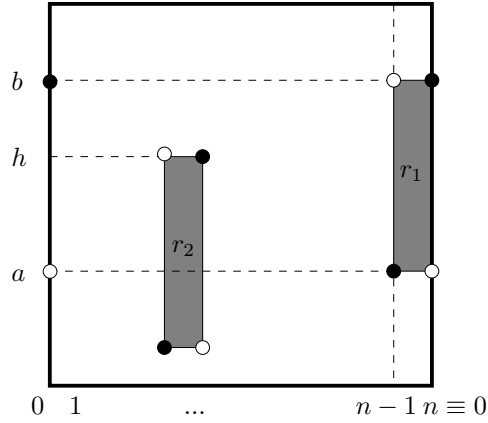}
\end {center}
\caption {{\bf Anti-Commutation of width one rectangles.}
In computing the sign of the rectangle $r_2$ we use the black dot on
the horizontal line of height $a$, while for the rectangle $r_1'$, which
has the same support, we use the white dot on the
leftmost vertical edge.}
\label {fig:WidthOne}
\end {figure}

\begin {figure}
\begin {center}
\input {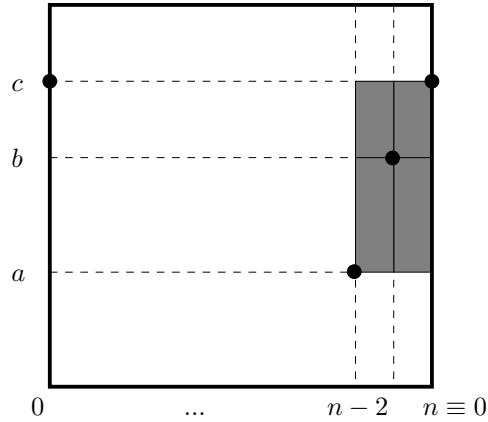}
\end {center}
\caption {{\bf Hexagon rule in the last two columns.}
We compute the signs in terms of the generator $\x;$ some of its
components are the black dots shown here.}
\label {fig:HexLastTwo}
\end {figure}

Let us now consider the hexagon rule, where the rectangles are supported 
in the last two columns, as in Figure~\ref{fig:HexLastTwo}. We denote by 
$a$, $b$, and $c$  (with $a<b<c$) the three possible heights at which the 
relevant rectangles have their horizontal edges. We use the notation 
from Definition \ref{def:VerticalSignAssignmentThin}, where $\x$ is the 
initial generator of $r_1$ and $r_1'$. When applying formulas 
\eqref{eq:ThinSignAssignment} and \eqref{eq:LastColumn}, we must keep 
in mind that the initial point of the relevant rectangle may differ  
from $\x$ (at heights $a$, $b$ and $c$), cf.\ Figure~\ref{fig:HexLastTwo}.  
After finding the relevant discrepancies, we can express everything in 
terms of $\x:$ 
\begin{align*}
        \sign(r_1)&= (-1)^{\NE\big(\x,\{(x_1,x_2)\in\x| x_2 \leq 
        b\}\big) + \NE\big(\x,\{(x_1,x_2)\in\x| b<x_2<c, \
        x_2~\text{even})\}\big)+c }\\
        \sign(r_1')&= (-1)^{\NE\big(\x,\{(x_1,x_2)\in\x| x_2 \leq
        b\}\big)} \\
        \sign(r_2) &= (-1)^{\NE\big(\x,\{(x_1,x_2)\in\x| x_2 \leq
        a\}\big) + \NE\big(\x,\{(x_1,x_2)\in\x| a<x_2<c \}\big)+1} \\
        \sign(r_2')&= (-1)^{\NE\big(\x,\{(x_1,x_2)\in\x| x_2 \leq
        a\}\big) + \NE\big(\x,\{(x_1,x_2)\in\x| a<x_2<c, \
        x_2~\text{even})\}\big) + b +c +1} \\
        \sign(r_3) &=  (-1)^{\NE\big(\x,\{(x_1,x_2)\in\x| x_2 \leq
        a\}\big) + \NE\big(\x,\{(x_1,x_2)\in\x| a<x_2\leq b, \ 
        x_2~\text{even})\}\big) + b} \\
        \sign(r_3')&= (-1)^{\NE\big(\x,\{(x_1,x_2)\in\x| x_2 \leq
        a\}\big) + \NE\big(\x,\{(x_1,x_2)\in\x| a<x_2<c \}\big)}
\end{align*}

Putting these relations together, we obtain the required identity:
$$\sign(r_1)\sign(r_2)\sign(r_3)=\sign(r_1')\sign(r_2')\sign(r_3').$$

There is a similar computation that needs to be done for the hexagon rule 
when the rectangles are supported in the first and the last column. We 
leave this case as an exercise for the reader.

Finally, we need to check the square and the hexagon rule when some of the 
rectangles involved wrap vertically around the torus, i.e., their support 
has nontrivial intersection with the horizontal line $l$ of height 
$n-\frac{1}{2}$. We call such rectangles {\em vertically wrapped}. For the 
square rule, either two or all four of the four rectangles involved are 
vertically wrapped. For the hexagon rule, exactly four out of the six 
rectangles involved are vertically wrapped. The square and the hexagon 
rules now follow from the corresponding ones when we replace 
the vertically wrapped rectangles $r$ with their counterparts $r'$ such 
that $r * r'$ are vertical annuli. Indeed, from equations of the form 
$\sign(r) = -\sign(r')$ we always pick up an even number of minus 
signs (either two or four), so the overall signs are unchanged.
\end{proof}

\step{4}{Extend vertical sign assignment to all empty rectangles}
We weaken the notion of sign assignments from
Definition~\ref{def:SignAssignment} a little.

\begin{definition}
  \label{def:VerticalSignAssignment}
  A {\em vertical sign assignment} is a function
  $$\sign\colon \EmptyRect\longrightarrow \{\pm 1\}$$
  which satisfies
  Properties~(Sq) and~(V)
  from Definition~\ref{def:SignAssignment}.  Sometimes, we call this a
  {\em vertical sign assignment on all empty rectangles}, to
  distinguish it from the seemingly weaker vertical sign assignments
  on thin rectangles.
\end{definition}

Our goal in this step is to show that a vertical sign assignment
for thin rectangles can be uniquely extended to a vertical sign
assignment on all empty rectangles.

\begin{figure}
\begin{center}
\mbox{\vbox{\epsfbox{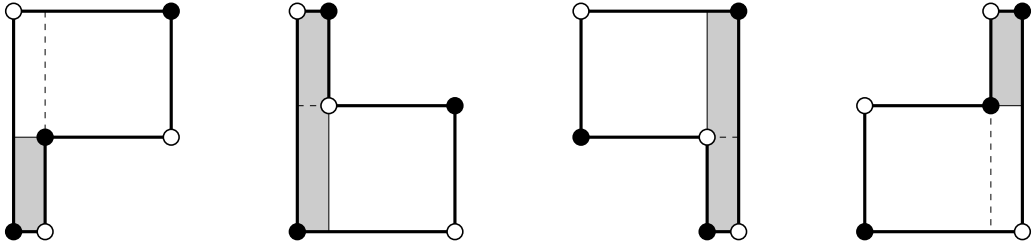}}}
\end{center}
\caption {{\bf The four cases in Lemma~\ref{lemma:narrow-rects}.}
  Each case is a decomposition $r_1 * r_2$, where $r_1$ is shown
  shaded.  The dotted line gives the alternate decomposition $r_1' *
  r_2'$.
\label{fig:FourCases}}
\end{figure}

Given a vertical sign assignment for thin rectangles $\sign_0$, we
define an extension
$$\sign\colon\EmptyRect\longrightarrow \{\pm 1\},$$
by extending the
definition inductively on the width $w$ of the rectangle.  Explicitly,
if $r\in\EmptyRect(\x,\y)$ is a rectangle with width $1$, then
$\sign(r)=\sign_0(r)$. Suppose next that $\sign$ is defined for all
rectangles of width less than $w$ for some $w > 1$.  Given
$r\in\EmptyRect(\x,\y)$ of width $w$,
there is exactly one rectangle $r_1$ ending at $\x$ with width one
whose upper left corner coincides with the lower left corner of $r$,
as in the first diagram in Figure~\ref{fig:FourCases}.
Then $r_1*r$ has an alternate decomposition as $r_2*r_3$, where
$r_2$ has width $w-1$, and $r_3$ has width one. We can then define
\begin{equation}
  \label{eq:Extension}
  \sign(r)=-\sign_0(r_1)\sign(r_2)\sign_0(r_3).
\end{equation}
(The right hand side is defined, since the width of $r_2$ is $w-1$.)

We verify that~$\sign$ is a vertical sign assignment in the sense of
Definition~\ref{def:SignAssignment}, which we do in stages.

\begin{lemma}
  \label{lemma:narrow-rects}
  Suppose that we have four rectangles
  $r_1,r_2,r_1',r_2'\in\EmptyRect$ with $r_1*r_2=r_1'*r_2'$, where
  $r_1$ and $r_2'$ have width one, and $r_1$ and $r_2$ share exactly
  one corner. Then
  \begin{equation}
    \label{eq:AntiCommutation}
    \sign(r_1)\sign(r_2)=-\sign(r_1')\sign(r_2').
  \end{equation}
\end{lemma}

\begin{proof}
  There are four cases in the proof, as illustrated in
  Figure~\ref{fig:FourCases}.  In the first case, when the upper-left
  corner of~$r_1$ is the lower-left corner
  of~$r_2$, the conclusion holds by the definition of the extension
  of~$\sign$, Equation~\eqref{eq:Extension}.
  In the second case, when the lower-right corner of~$r_1$ is the
  lower-left corner of~$r_2$, the conclusion follows
  from Property~(V) of the vertical sign assignment:
  if we relate both thin rectangles, $r_1$ and~$r_2'$, to the other
  rectangle in the same vertical column, we get the same rectangles
  as in the previous case.

  Otherwise, we will prove the result by induction on the maximum of
  the widths of~$r_2$ and~$r_1'$.  We treat the base case first, where
  this maximum is equal to two.  There are two cases:
  \begin{itemize}
  \item If the upper-left corner of~$r_1$ is the upper-right corner
    of~$r_2$, Equation~\eqref{eq:AntiCommutation} follows from
    Property~(Sq) of the vertical sign assignment,
    the hexagon rule, and the definition of~$\sign$.
  \item If the lower-right corner of~$r_1$ is the upper-right corner
    of~$r_2$, Equation~\eqref{eq:AntiCommutation} follows
    from Property~(V) applied to the case above.
  \end{itemize}

We now treat the induction on the maximum of the widths of~$r_2$
and~$r_1'$, which we may assume is greater than two.  If the
upper-left corner of~$r_1$ is the upper-right corner of~$r_2$ (the
third case in the figure), we may
again apply Property~(V) to change to the last
case.  So we may assume
that the lower right corner of~$r_1$ is shared with the upper right
corner of~$r_2$.  Since the width of~$r_2$ is greater than two, we can
find a thin
rectangle~$r_0$ to~$\x$, the initial point of~$r_1$, disjoint
from~$r_1$ with the property that $r_0$ and~$r_2$ share a corner.  We
consider now the composite $p=r_0*r_1*r_2$.  We organize the various
decompositions of~$p$ into the graph in Figure~\ref{fig:Cube}.
\begin{figure}
  \mbox{\vbox{\includegraphics[width=5in]{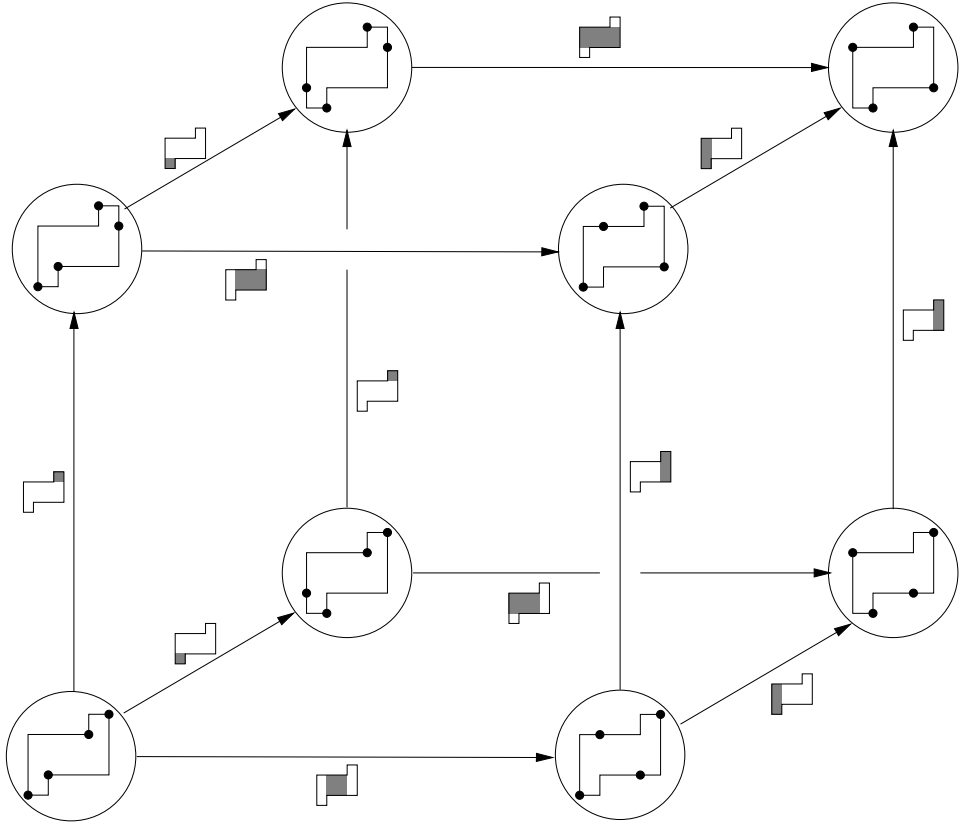}}}
  \caption{\textbf{The cube of decompositions in
      Lemma~\ref{lemma:narrow-rects}}.  Each vertex is a configuration
  (shown by black dots); each edge is a rectangle (shown shaded).  The
  different decompositions are the different ways to go from the
  lower-left corner to the upper-right corner following three edges.
  We are trying to show that the back face anti-commutes; the front
  face is the inductive case.}
  \label{fig:Cube}
\end{figure}
Each edge corresponds to some rectangle in some decomposition of $p$,
and each vertex corresponds to one of the various elements of $\S$ which can
be connected by rectangles in some decomposition of $p$. Each face
corresponds to two decompositions of $p$ which have some rectangle in
common. Thus we have organized the decompositions of
$p$ in a cube.  A face is said to {\em anti-commute} if the product of
the signs associated to its four edges is $-1$. Our goal is to show
that the face belonging to the two decompositions $r_0*r_1'*r_2'$ and
$r_0*r_1*r_2$ anti-commutes.  To see this, we observe that the other
five faces anti-commute: two of them do, as they correspond to rearranging
two disjoint thin rectangles (Property~(Sq) of the
vertical sign assignment for thin rectangles), two of the faces
anti-commute by the definition of the extension,
Equation~\eqref{eq:Extension}, and the fifth face anti-commutes by
induction on the width of~$r_2$.
\end{proof}

\begin{lemma}
        \label{lemma:ArbitraryWidthOne}
  Suppose that we have four distinct rectangles
  $r_1,r_2,r_1',r_2'\in\EmptyRect$ with $r_1*r_2=r_1'*r_2'$, where
  $r_1$ and $r_2'$ have width one, then
\[
    \sign(r_1)\sign(r_2)=-\sign(r_1')\sign(r_2').
\]
\end{lemma}

\begin{proof}
        We prove the result by induction on the width of $r_2$.

        The case where $r_2$ shares exactly one corner with~$r_1$ was
        handled in Lemma~\ref{lemma:narrow-rects}, so we may assume
        the corners of~$r_1$ and~$r_2$ are distinct.
        According to Property~(V), we can assume
        without loss of generality that the supports of~$r_1$ and~$r_2$
        are disjoint.  If~$r_2$ also has width one, we are done by
        Property~(Sq).

        Otherwise, we proceed by induction in a way similar to the
        last case
        of Lemma~\ref{lemma:narrow-rects}.  Specifically, we can find
        another width one rectangle~$r_0$ ending at the initial point
        of~$r_1$ which shares one corner with~$r_2$. Consider the polygon
        $p=r_0*r_1*r_2$. We can once again organize the various
        decompositions of this polygon into a cube.
        In the case where~$r_0$ and~$r_1$
        are disjoint, two of the faces anti-commute according
        to Property~(Sq) of $\sign_0$, two
        anti-commute by Lemma~\ref{lemma:narrow-rects}, and a fifth
        anti-commutes by induction on the width of $r_2$. Thus,
        the sixth must anti-commute, as well. 

        In the other case, where~$r_0$ and~$r_1$ share a corner
        (which we need to consider only when the width of $r_2$
        is three), four of the squares in the cube
         anti-commute according to Lemma~\ref{lemma:narrow-rects}.
         The fifth involves a domain $r_1'*r_2'$ where~$r_1'$
         and~$r_2'$ are both rectangles of width two, with $r_1'$
         containing~$r_1$ and $r_2'$ contained in~$r_2$.  To show that
         this last face anti-commutes, we can find
         another rectangle~$r_0'$ that can be pre-composed with
         these two to give a domain $r_0'*r_1'*r_2'$.  The various
         decompositions of this domain can again be arranged into a
         cube, in which four of the squares anti-commute by
         Lemma~\ref{lemma:narrow-rects}, the fifth anti-commutes by
         Property~(Sq), and the sixth is the
         face with domain $r_1'*r_2'$, as desired.
\end{proof}

\begin{prop}
        \label{prop:VerticalSignAssingnments}
  Given any four empty rectangles
  $r_1,r_2,r_1',r_2'\in\EmptyRect$ with $r_1*r_2=r_1'*r_2'$, we have that
\[    \sign(r_1)\sign(r_2)=-\sign(r_1')\sign(r_2').
\]
\end{prop}

\begin{proof}
        Using suitably-placed width one rectangles as before, we can
        narrow $r_1$ and/or $r_2$ by induction, until one or the other has
        width one and hence is covered by
        Lemma~\ref{lemma:ArbitraryWidthOne}.
\end{proof}

In effect, the above proposition shows that a vertical sign assignment for 
thin rectangles~$\sign_0$ can be canonically extended to a vertical sign 
assignment for arbitrary rectangles~$\sign$ in the sense of
Definition~\ref{def:VerticalSignAssignment}.

\step{5}{Signs for rectangles supported in $\Square$}

In Proposition~\ref{prop:VerticalSignAssignmentThin} we constructed a 
vertical sign assignment for thin rectangles. According to Step 4 
above, this gives a vertical sign assignment $\sign:\EmptyRect \to 
\{\pm1\} $ on all empty rectangles. We aim to give an explicit formula for 
$\sign(r)$ in the case when $r$ is supported in the subsquare $\Sigma = 
[0,n-1] \times [0,n-1]$.

\begin {figure}
\begin {center}
\input {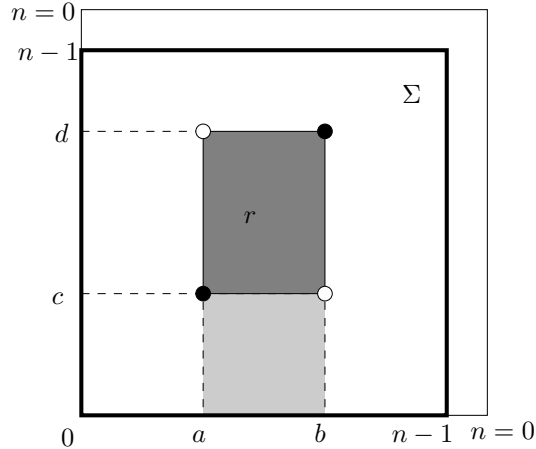}
\end {center}
\caption {{\bf Signs for rectangles on the square.}
The subsquare $\Square$ of the torus is shown bordered by the thick 
lines. Inside we have a rectangle $r$. The quantity $\Down(r)$ counts the 
number of points in $\x$ which lie in the lightly shaded region below~$r$.}
\label {fig:SquareSigns}
\end {figure}

\begin {prop}
\label {prop:RectanglesSquare}
Let $r = [a,b] \times [c,d] \in \EmptyRect(\x, \y)$, with $0 \leq a < b 
\leq n-1$ and $0 \leq c < d \leq n-1$. Denote by $\Down (r)$ the number of 
points $(x_1, x_2) \in \x$ which lie strictly below~$r$, i.e., $x_1 \in 
(a,b)$ and $x_2 \in [0,c)$, cf.\ Figure~\ref{fig:SquareSigns}. Then:
\begin {equation}
\label {eq:SquareSign}
 \sign(r) = (-1)^{\NE\big(\x,\{(x_1,x_2)\in\x| x_2 \leq
d\}\big) + \Down(r) \cdot \big (\NE\big(\x,\{(x_1,x_2)\in\x| c<x_2\leq 
d\} \big) + 1 \big)}.  
\end {equation}
\end {prop}

\begin{proof}
We use induction on the width~$w$ of~$r$. When $r$ has width one, 
\eqref{eq:SquareSign} is just the formula~\eqref{eq:ThinSignAssignment}.

Assume $w > 1$. We distinguish two cases, according to the position of
the point $(a+1,l) \in \x$ in the vertical line $x_1 = a+1$. 
When it is below the support of~$r$, the inductive step follows by using the 
anti-commutation rule for a decomposition of the form $r * r_1 = r_2 * r'$, 
where the rectangles $r_1$ and $r_2$ have width one and lie in the $a\th$ 
column, so that their signs can be computed using 
\eqref{eq:ThinSignAssignment}, and $r'$ has width $w-1$. The latter case is 
similar, but the decompositions are of the form $r_1 * r = r' * r_2$. In 
both cases the computations are straightforward. 
\end{proof}

\step{6}{From vertical sign assignments to true sign assignments}

\vskip.2cm

\begin{lemma}
        \label{lemma:RowFunction}
        If $\sign$ is a vertical sign assignment, then there is
        some function $\rho\colon \{1,\ldots,n\}\rightarrow \{ \pm 1\}$
        with the property that if $r_1$ and $r_2$ are two height one
        rectangles, with $r_1*r_2$ connecting some $\x\in\S$ to itself,
        then $\sign(r_1)\sign(r_2)=\rho(i)$, where the support of $r_1
        * r_2$ consists of the $i\th$ row. 
\end{lemma}

\begin{proof}
	We claim first that when the support of $r_1*r_2$ is a horizontal
    annulus $[0,n)\times [i,i+1]$, $\sign(r_1)\sign(r_2)$ depends only
    on~$i$ and the components of $\x$ in the $i\th$ and
    $(i+1)\st$ rows. This follows
	readily from the square rule; if $\x$ and $\x'$ are two generators
    whose components on the $i\th$ and $(i+1)\st$ agree and there is some
	$r\in\EmptyRect(\x,\x')$, then two applications of the
	square rule establish the claim. More generally, we can always
	get between two generators $\x$ and $\x'$ that agree on these two
    rows by a sequence of squares whose corners
	are disjoint from the $i\th$ row, verifying the claim.

	We next verify that if $r_1\in\EmptyRect(\x,\y)$ are
	$r_2\in\EmptyRect(\y,\x)$ are two rectangles which, together,
	form the $i\th$ row, then $\sign(r_1)\sign(r_2)$ depends
	only on the row in which $r_1$ and $r_2$ are supported. To
	this end, observe that there is another rectangle
	$r_3\in\Rect(\x,\z)$ with height one supported in the row
	$i+1$, and $r_4\in\Rect(\z,\x)$.  We now claim that
	$r_1*r_2*r_3*r_4$ differs from another decomposition
	$r_1'*r_2'*r_3'*r_4'$ in two steps, in such a way that the
	supports of $r_3$ and $r_4$ agree with those of $r_3'$ and
	$r_4'$ (hence together they occupy the $i+1^{st}$ row), and
	the supports of $r_1'$ and $r_2'$ also occupy the $i\th$
	row, but the support of $r_1'$ and $r_2'$ are different from
	the supports of $r_1$ and $r_2$. Thus, it follows from the
	square rule that
	$\sign(r_1)\sign(r_2)=\sign(r_1')\sign(r_2')$. It
	is easy to see that any two $r_1\in\EmptyRect(\x,\y)$ and
	$r_2\in\EmptyRect(\y,\x)$ which occupy the $i\th$ row can be
	connected by a finite sequence of such moves.
\end{proof}

We now specialize to the vertical sign assignment $\sign$ constructed (in 
Step 3) from the vertical sign assignment on thin rectangles exhibited in 
Proposition~\ref{prop:VerticalSignAssignmentThin}, which was based on the 
formulas \eqref{eq:ThinSignAssignment} and \eqref{eq:LastColumn}. We claim 
that the function $\rho$ from Lemma~\ref{lemma:RowFunction} is
identically~$1$.

\begin{lemma}
        \label{lemma:RhoSquare}
$\rho(i) = 1$ for $i = 1, \dots, n-1$.
\end{lemma}

\begin{proof} It suffices to check that $\sign(r_1) \sign(r_2) = 1$ when 
$r_1 \in \Rect(\x, \y)$ is of the form $[0, n-1] \times [i-1, i]$, and $r_2 
\in \Rect(\y, \x)$ is the square $[n-1, n] \times [i-1, i]$ in the last 
column, cf. Figure~\ref{fig:RhoSquare}. Proposition~\ref{prop:RectanglesSquare} 
gives
$$ \sign (r_1) = (-1)^{\NE\big(\x,\{(x_1,x_2)\in\x| x_2 \leq
i\}\big) + i + 1} = (-1)^{\NE\big(\x,\{(x_1,x_2)\in\x| x_2 \leq
i-1\}\big) + 1} .$$
On the other hand, from Equation \eqref{eq:LastColumn} we get
$$ \sign(r_2) = (-1)^{\NE\big(\y,\{(y_1,y_2)\in\y| y_2 \leq
i-1\}\big) + i} = (-1)^{\NE\big(\x,\{(x_1,x_2)\in\x| x_2 \leq
i-1\}\big) + 1}. $$

\begin {figure}
\begin {center}
\input {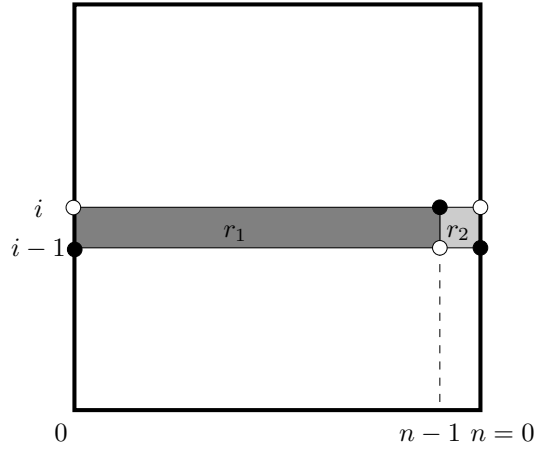}
\end {center}
\caption {{\bf Computing $\rho(i)$ for $i \leq n$.}
The generator $\x$ is shown in black dots, and $\y$ in hollow dots. When 
computing the sign of $r_2$ with formula \eqref{eq:LastColumn}, we need 
to use $\y$ in the place of $\x$.} \label {fig:RhoSquare}
\end {figure}

\end{proof}

\begin {figure}
\begin {center}
\input {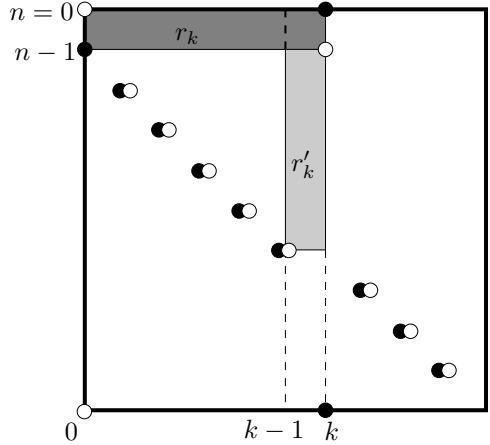}
\end {center}
\caption {{\bf Computing $\rho(n)$.}
The rectangles $r_k$ is darkly shaded, and $r_k'$ lightly shaded. The 
generators $\x_k$ and 
$\y_k$ are represented by the black and white dots, respectively. There is 
an alternate decomposition of $r_k * r_k'$ given by cutting along the 
dashed line in the top row.} \label {fig:LastRho}
\end {figure}

To check that $\rho (n) = 1$, we first prove the following:

\begin {lemma}
Let $k\in {1, 2, \dots, n-1}$. Denote by $\x_k, \y_k \in \S$ the 
configurations
\begin {eqnarray*}
 \x_k &=& \{(i, n-1-i) | 0 \leq i < k\} \cup  \{(k,0)\} \cup \{(i, 
n-i) | k < i < n\}, \\
\y_k &=& \{(0, 0)\} \cup  \{(i, n-1-i) | 1 \leq i < k\} \cup \{(k, 
n-1)\} 
\cup \{(i,n-i) | k < i < n\}.
\end {eqnarray*}
Let $r_k \in \Rect(\x_k, \y_k)$ be the rectangle of width $k$ and height 
$1$ supported in the last row, cf. Figure~\ref{fig:LastRho}. Then, its sign 
is given by:
$$ \sign (r_k) = (-1)^n. $$
\end {lemma}

\begin {proof}
Induction on $k$. When $k=1$, the sign of $r_1$ is minus the sign 
of a thin rectangle of width one and height $n-1$ supported in the first 
column. The latter can be computed using
Equation \eqref{eq:ThinSignAssignment}, which gives the answer 
$(-1)^{n-1}$; therefore, $\sign(r_1) = (-1)^n$. 

For $k > 1$, let $r_k' \in \Rect(\y_k, \z_k)$ be the rectangle of 
width one supported in the $k\th$ column, cf. Figure~\ref{fig:LastRho}. Its 
sign is $(-1)^{n+k}$ by formula \eqref{eq:ThinSignAssignment}. The domain 
$r_k * r_k'$ has an alternate decomposition as $p_k * r_{k-1}$, where 
the rectangle $p_k$ is again supported in the $k\th$ column and has a 
counterpart $p_k'$ such that $p_k * p_k'$ is a vertical annulus. Formula 
\eqref{eq:ThinSignAssignment} gives $\sign(p_k') = (-1)^{n+k}$, so that 
$\sign(p_k) = (-1)^{n+k+1}$. The inductive step now follows form the 
anti-commutation relation $\sign(r_k)\sign(r_k') = -\sign(r_{k-1})\sign(p_k)$. 
\end {proof}

\begin {lemma}
\label {lemma:LastRho}
$\rho (n) =1$.
\end {lemma}

\begin {proof}
By the previous lemma, the sign of $r_{n-1} \in \Rect(\x_{n-1}, \y_{n-1})$
is $(-1)^n$. There is also a thin vertical rectangle $r' \in  
\Rect(\x_{n-1}, \y_{n-1})$ supported in the last column, whose sign is 
$(-1)^{n-1}$ by \eqref{eq:LastColumn}. For the little square $r'' \in
\Rect(\y_{n-1}, \x_{n-1})$ supported in the top right corner we have 
$\sign(r'') = - \sign(r')$. Since $r_{n-1} * r'' $ is a horizontal annulus 
and  $\sign(r'') =\sign(r_{n-1})$, we get $\rho(n) =1$, as desired. 
\end {proof}

We can now complete Step~6:

\begin{proof}[Proof of Theorem~\ref{thm:SignAssignments}]
We proved existence of vertical sign assignments in Step 4. According to 
Lemma~\ref{lemma:RowFunction}, the resulting signs give a true sign 
assignment provided that the function $\rho$ defined there is identically 
one. This was checked in Lemmas~\ref{lemma:RhoSquare} and 
\ref{lemma:LastRho}.

To see uniqueness, let $\sign_1$ and $\sign_2$ be two true sign 
assignments. We can restrict them to the square $\Square$ and get 
sign assignments on the Cayley graph (cf.\ 
Lemma~\ref{lemma:CayleyGraphSignChar}). Using the uniqueness part of  
result of Proposition~\ref{prop:sign-cayley}, we obtain a function 
$f\colon \S\longrightarrow \{\pm 1\}$ such that the composite
$$\B: \EmptyRect \longrightarrow \{\pm 1\}, \ \ \B(r) = f(\x) \cdot f(\y) 
\cdot \sign_1(r) \cdot \sign_2(r) \text{ for } r \in \EmptyRect(\x, \y)$$
satisfies the following properties:

\begin {itemize}
\item $\B(r) = 1$ for any rectangle of width one supported in $\Sigma;$
\item (Commutation rule) $\B(r_1)\cdot \B(r_2) = \B(r_1') \cdot \B(r_2')$ 
whenever $r_1, r_2, r_1', r_2' \in \EmptyRect$ are distinct and satisfy 
$r_1 * r_2 = r_1' * r_2';$
\item $\B(r_1) = \B( r_2)$ if $r_1 * r_2$ is a vertical annulus;
\item $\B(r_1) = \B(r_2)$ if $r_1 * r_2$ is a horizontal annulus.
\end {itemize}

We claim that $\B$ is identically~$1$. Indeed, the third property above 
implies that $\B(r) = 1$ whenever $r$ has width one and is not supported 
in the last column. The fourth property implies that $\B(r) = 1$ when $r$
is a square of side length one supported in the last column; the same must 
be true for all width one rectangles in the last column by induction on 
height, applying the commutation rule. Thus $\B$ takes the value one on 
all vertical thin rectangles. The fact that $\B(r) = 1$ for all $r$ 
follows by induction on width, again using the commutation rule. This 
shows that $f$ satisfies the property required in the statement of
Theorem~\ref{thm:SignAssignments}.
\end {proof}

\subsection{Properties of the sign-refined chain complex, and the proof
of Theorem~\ref{thm:WithSigns}}
\label{subsec:DiscussionOverZ}

\begin{prop}
  \label{prop:DSquaredZeroZ}
  Let $\sign$ be a sign assignment.
  The $\Z[U_1,\ldots,U_n]$-module $\Cm(G)$, endowed with the endomorphism
  $\dm_\sign$, is  a chain complex. Moreover, if $\sign_1$ and $\sign_2$
  are two different sign assignments, then there is an isomorphism
  of chain complexes 
  $(\Cm,\dm_{\sign_1})\cong (\Cm,\dm_{\sign_2})$.
\end{prop}

\begin{proof}
  In the expression $\partial_\sign\circ\partial_\sign(\x)$, terms
  can be paired off as in the proof of
  Proposition~\ref{prop:DSquaredZero}. These terms cancel, according
  to the axioms on $\sign$. 
    
  Suppose we are given sign assignments $\sign_1$ and $\sign_2$.
  Consider the map 
  $$\Phi\colon (\Cm(G),\partial_{\sign_1}) \longrightarrow
  (\Cm(G),\partial_{\sign_2})$$
  defined by $\Phi(\x)=f(\x)\cm \x$,
  where $f$ is the function from
  Theorem~\ref{thm:SignAssignments}. It is straightforward to see
  that $\Phi$ is an isomorphism of chain complexes.
\end{proof}

Other algebraic properties from Section~\ref{sec:BasicProperties}
have straightforward generalizations to this context. For example:

\begin{lemma}
  \label{lemma:ChainHomotopiesZ}
  Suppose that $O_i$ and $O_k$ correspond to the same component of
  $\oL$.  Then multiplication by $U_i$ is filtered chain homotopic to
  multiplication by $U_k$.
\end{lemma}

\begin{proof}
  The chain homotopy from Lemma~\ref{lemma:ChainHomotopiesZ} works to
  establish the present lemma.  It is important here that for
  $\x\in\S$, and $r_1$ and $r_2$ are the decompositions of the
  horizontal annulus containing $X_j$, and $r_3$ and $r_4$ are the
  analogous decompositions of the vertical annulus, then
  $\sign(r_1)\sign(r_2)=-\sign(r_3)\sign(r_4)$, but this is
  ensured by the axioms of the sign assignment. This ensures that
  $U_i$ is chain homotopic to $U_k$ (rather than, say, $-U_k$).
\end{proof}

Again, we view the chain complex $(\Cm(G),\dm)$ 
as a module over $\Z[U_1,\ldots,U_\ell]$
where the $\{U_i\}_{i=1}^\ell$ correspond to the $\ell$ components
of our link $\orL$. As before, we have the following:

\begin{prop}
  \label{prop:HLaa}
  Suppose that the oriented link $\oL$ has $\ell$ components.  Choose
  an ordering of $\Os=\{O_i\}_{i=1}^n$ so that for
  $i=1,\ldots,\ell$, $O_i$ corresponds to the $i\th$ component of
  $\oL$.  Then the filtered chain homotopy type of $\Cm(G)$,
  viewed as a chain complex over $\Z[U_1,\ldots,U_\ell]$, is
  independent of the ordering of $\Os$. Then
  $\HLa(G)$ and $\HLaa(G)$ are finitely generated Abelian
  groups. Moreover,
  \[H_*(\CLaa(G)) \cong \HLa(G)\otimes
  \bigotimes_{i=1}^\ell V_i^{\otimes (n_i-1) }, \] where $V_i$ is the
  two-dimensional vector space spanned by two generators, one in zero
  Maslov and Alexander multi-gradings, and the other in Maslov grading
  minus one and Alexander multi-grading corresponding to minus the
  $i\th$~basis vector.
\end{prop}

\begin{proof}
  This is a routine adaptation of
  Proposition~\ref{prop:DoesntDependOnOrdering},
  Lemma~\ref{lemma:HLaFinGen}, and Proposition~\ref{prop:HLaaNoSigns}
  from Section~\ref{sec:BasicProperties}.
\end{proof}

We now turn to the proof of the main theorem,
Theorem~\ref{thm:WithSigns}, with signs.

We adopt the strategy from Section~\ref{sec:Invariance}; however, we
must specify the signs used in defining our various chain maps and
chain homotopies, and verify that they are indeed chain maps and chain
homotopies with appropriate signs. As in Section~\ref{sec:Invariance},
we begin with the case of commutation.

We adopt notation from Subsection~\ref{subsec:Commutation}. Consider
the pentagons $\Pent_{\beta\gamma}(\x,\y)$ used there.
Straightening out the $\beta\cap\gamma$-corner of the
pentagons naturally induces rectangles in $G$. (We could
have in fact defined the map in Section~\ref{subsec:Commutation} as
counts of rectangles, where the $O$'s and $X$'s in the central column
are moving, but then it would have been a little confusing to write
down exactly when they are counted with powers of the $U$'s.)

Formally, we obtain a ``straightening map''
$$\straight\colon \Pent_{\beta\gamma}(\x,\y)\longrightarrow \Rect(\x,\y'),$$
where $\y'$ is the generator corresponding to $\y$, where we slide
horizontally from the $\gamma$ back to the $\beta$ component.
Clearly, the image of $\straight$ consists of rectangles with a vertical
segment along $\beta$. There are two possibilities: either the rectangle
lies to the right of this vertical segment (i.e., the segment is a left
edge of the rectangle), or it lies to the left of this vertical segment.
In the first case, we say the pentagon is a {\em right pentagon}, in the
latter, we say it is a {\em left pentagon}. In Figure~\ref{fig:Pentagon},
the one pictured on the left is a left pentagon, and the one on the right is
a right pentagon.

We define
$$\Phi_{\beta\gamma}(\x)=\sum_{\y\in\S(H)}\,
\sum_{\substack{p\in \Pent_{\beta\gamma}(\x,\y)\\ \x\cap\Interior(p)=\emptyset}}
\epsilon(p)\cm U_1^{O_1(p)}\cdots U_n^{O_n(p)}\cm \y,$$
where 
$$\epsilon(p)=\begin{cases}
\sign(\straight(p)) & {\text{if $p$ is a left pentagon}} \\
-\sign(\straight(p)) & {\text{if $p$ is a right pentagon}}. \\
\end{cases}$$

We obtain the following analogue of Lemma~\ref{lemma:PhiChainMap}:

\begin{lemma}
  The map $\Phi_{\beta\gamma}$ is a filtered anti-chain map, i.e.,
  $$\partial\circ \Phi_{\beta\gamma}+\Phi_{\beta\gamma}\circ\partial=0.$$
\end{lemma}

\begin{proof}
  Again, the proof follows from the proof of
  Lemma~\ref{lemma:PhiChainMap}.  In fact, the fact that the terms
  cancel in pairs typically follows from the same pairing which we see
  in Proposition~\ref{prop:DSquaredZeroZ}.  There are two cases which
  look different, though.  One of these corresponds to the two
  decompositions pictured as in Figure~\ref{fig:ChainMap}. After projecting
  via $\straight$, both decompositions of the composite region in this case
  correspond to the same decomposition of the composite region into two
  rectangles. However, in one case, the rectangle corresponding to the pentagon
  is on the left, in the other, it is on the right. Thus, the signs
  given by $\epsilon$ are opposite. The other case, the rotation by
  $180^\circ$ of Figure~\ref{fig:ChainMap}, works similarly.
\end{proof}

It seemed more natural in the above proposition to consider anti-chain
maps, rather than the more traditional chain maps. Just as chain maps
induce maps on homology, so do anti-chain maps.  One could
alternatively consider the chain map ${\widetilde \Phi}$ defined by
$${\widetilde\Phi}_{\beta\gamma}(\x)=(-1)^{M(\x)} \cm \Phi_{\beta\gamma}.$$

We now turn to the chain homotopies gotten by counting hexagons.

Once again, there is a straightening map
\begin{align*}
\straight'\colon \Hex_{\beta\gamma\beta}(\x,\y)&\longrightarrow \Rect(\x,\y),
\end{align*}
and we can define a homotopy operator
$H_{\beta\gamma\beta}\colon C(G) \longrightarrow C(G)$ by
\begin{align*}
  H_{\beta\gamma\beta}(\x) =& \sum_{\y\in \S(G)}\,
  \sum_{\substack{h\in\Hex_{\beta\gamma\beta}(\x,\y)\\ \x\cap\Interior(h)=\emptyset}}
  \epsilon'(h)\cm U_1^{O_1(h)}\cdots U_n^{O_n(h)} \cm \y,
\end{align*}
where
$$\epsilon'(h)=\sign(\straight'(h)).$$
Similarly define $H_{\gamma\beta\gamma}$.

\begin{prop}
\label{prop:CommuteZ}
With respect to the sign refinements, the map $\Phi_{\beta\gamma}$
induced by commuting two columns induces an isomorphism in homology.
\end{prop}

\begin{proof}
  The proof of Proposition~\ref{prop:Commute} adapts readily to show that
  \begin{align*}
    \Id+\Phi_{\gamma\beta}\circ \Phi_{\beta\gamma} + \partial \circ H_{\beta\gamma\beta} + H_{\beta\gamma\beta}\circ \partial  &= 0 \\
    \Id+\Phi_{\beta\gamma}\circ \Phi_{\gamma\beta} + \partial \circ H_{\gamma\beta\gamma} + H_{\gamma\beta\gamma}\circ \partial &=0.
\end{align*}
Note that in the terms $\Phi_{\gamma\beta}\circ \Phi_{\beta\gamma}$
and $\Phi_{\beta\gamma}\circ \Phi_{\gamma\beta}$, the two pentagons
that appear are either both right pentagons or both left pentagons, so
the extra minus sign for right pentagons has no effect.  The
proposition now follows.
\end{proof}

With commutation invariance in hand, we now turn to stabilization invariance,
following the steps in Subsection~\ref{subsec:Stabilization}.

As a first step, we need a sign in the definition of $C'$,
the mapping cone of the chain map
$$U_1-U_2\colon B[U_1]\longrightarrow B[U_1],$$
to ensure it is a chain complex. One way of doing this
is to define
$$\partial'(a,b)=(\partial a,(U_2-U_1)\cm a - \partial b).$$

We will find the following terminology useful.

\begin{definition}
  Suppose $p\in\pi(\x,\y)$ can be decomposed as $p=r_1*\cdots*r_m$ for
  some $n$ and $r_i\in\EmptyRect$. Suppose moreover that for some $i$ we have
  $r_i*r_{i+1}=r'_i*r_{i+1}'$, for some $r_i',
  r'_{i+1}\in\EmptyRect$. Then we say that the decompositions
  $r_1*\cdots*r_m$ and $r_1*\cdots*r'_i*r'_{i+1}*\cdots*r_m$ differ by an {\em
    elementary move}.
\end{definition}

Recall that in Subsection~\ref{subsec:Stabilization}, we identified
$C(G)$ with a chain complex whose generators are $\Interval\subset
\S(H)$, which contain the distinguished point $x_0$. The differentials
count either empty rectangles which do not contain $x_0$ as corner
points, or those of ``Type 2'', i.e., those which include $x_0$ (and hence
also $O_1$ and $X_1$) in their interior. These rectangles are counted,
but not with a power of $U_1$. 

\begin{figure}
\begin{center}
\mbox{\vbox{\epsfbox{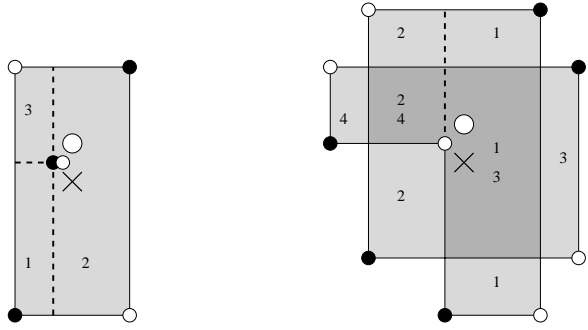}}}
\end{center}
\caption {{\bf Decomposing polygons.}
        On the left, we a have a standard decomposition of a Type 2
        rectangle: more precisely, the decomposition consists of
        $r_1*r_2*r_3$, and the number on a region indicates which
        rectangle it belongs to. On the right, we have indicated
        the standard decomposition of a complexity
        $5$ polygon of type $L$, $r_1*\cdots*r_4$. Some regions
        (which have local multiplicity $2$) are contained in the support
        of more than one rectangle, and hence are labelled with more than
        one integer.}
\label{fig:DecomposePolygons}
\end{figure}

Starting with a sign assignment~$\sign$ on rectangles in~$H$, we
can induce one on~$\Interval$ with the above two types of
differential.  To do
this, we must give an explicit decomposition of each Type 2 rectangle
$r'$ as a product of three empty rectangles. There are four ways of
doing this. We choose the following one: the initial rectangle is
lower left (and involves $x_0$), and the third (final) one uses the
upper left corner, cf.\ Figure~\ref{fig:DecomposePolygons}, and call it
the {\em standard decomposition}~$D_0(r')$ of the rectangle~$r'$. For
consistency, if $r$ is a Type 1 rectangle in $G$, and $r'$ is the
corresponding rectangle viewed as a rectangle in $H$, we let $D_0(r')$
denote the length $1$ decomposition $D_0(r')=r'$.

\begin{lemma}
  Fix $\x, \y\in\S(G)$ and 
  let $r\in \EmptyRect(\x,\y)$
  correspond to the rectangle $r'$ connecting $\x',\y'\in\S(H)$
  under the correspondence above between $\S(G)$ with $\Interval
  \subset S(H)$.  For any sign assignment~$\sign$ for $H$,
  define
  $$\sign_0(r)=
    \begin{cases}
      \sign(r') & {\text{if $r'$ has Type 1}} \\
      \sign(r_1)\sign(r_2)\sign(r_3)  & {\text{if $r'$ has Type 2 and
          $D_0(r') = r_1*r_2*r_3$}}.
    \end{cases}
    $$
    Then $\sign_0$ induces a sign assignment in the sense of
    Definition~\ref{def:SignAssignment} for rectangles in $G$.
\end{lemma}

\begin{proof}
  We must show that if $r_1*r_2=r_1'*r_2'$ in $G$, then the
  decompositions $D_0(r_1)*D_0(r_2)$ and $D_0(r_1')*D_0(r_2')$ can be
  connected by an odd number of elementary moves. This follows from a
  routine case analysis of how~$r_1$ and~$r_2$ can interact.  When one
  rectangle is Type~1 and the other is not, it takes $3$, $5$, or $7$
  elementary moves to connect the two decompositions.  The most
  complicated case is the one where both are Type 2 rectangles; in
  that case, we can connect the two decompositions by nine moves:
  Write the standard decomposition $D_0(r_1)=s_1*s_2*s_3$ and
  $D_0(r_2)=t_1*t_2*t_3$, then we have a decomposition
  $D_0(r_1)*D_0(r_2)=s_1*s_2*s_3*t_1*t_2*t_3$.  In three moves, we
  commute $t_1$ to the beginning, then in three more moves we commute
  $t_2$ to the second place, and in three more moves we commute $t_3$
  to the third spot. The resulting decomposition can be easily seen to
  agree with $D_0(r_1')*D_0(r_2')$.  The other two axioms of a sign
  assignment are also easily verified.
\end{proof}

\begin{remark}
	By ``commuting the $t_1$ to the beginning of
	$s_1*s_2*s_3*t_1*t_2*t_3$'', we mean the following string of
	operations: apply three elementary moves, the first of which
	replaces the consecutive terms $s_3*t_1$ by an alternative
	pair $t_1'*s_3'$, then apply another elementary move to the
	consecutive pair $s_2*t_1'$, to get $t_1''*s_2'$, and finally
	apply an elementary move to the pair $s_1*t_1''$ to get
	$t_1'''*s_1'$. We will use this shorthand in several
	future proofs, as well.
\end{remark}

We need now to introduce signs in the definition of the stabilization
map $F$ of Equation~\eqref{eq:StabilizationMap} to ensure that it is,
in fact, a chain map.

As a first step, we define a function
$$
\mu\colon \pF \longrightarrow \{\pm 1\}.
$$
For this, we give a specific decomposition of $p\in
\pF(\x,\y)$ as an ordered juxtaposition of
rectangles. Specifically, recall that $\partial p$ can be thought of
as an oriented, connected, curve. 
Order now the $\beta$-arcs
$\{v_i\}_{i=1}^m$ so that they inherit the cyclic ordering 
from the orientation of $\partial p$, and so that
$v_m$ contains the stabilization point $x_0$
(which in turn is a component of $\y$).  We can decompose
$$p=r_1*\cdots*r_{m-1},$$
where $r_i$ is a rectangle containing $v_i$,
compare Figure~\ref{fig:DecomposePolygons}, and define
\[
\mu(p) = S(r_1)\cdots S(r_{m-1}).
\]
Note that the left edge of
each odd rectangle is contained in the circle $\beta_1$ containing
$x_0$, while the right edge of each even rectangle is contained
in~$\beta_1$.  We call this
decomposition the {\em standard decomposition}~$D(p)$.
For polygons with complexity $m$, there are $m-1$ rectangles in
this decomposition; $m$ is odd if the polygon is of type $L$
and even if the polygon is of type $R$.

We will analyze the signs according to the cases in the proof of
Lemma~\ref{lemma:FIsChain}.

\begin{lemma}
  \label{lemma:CancelSignsType1} Fix a complexity $m$ domain $p\in\pF$
  and a rectangle $r$ as in cases~I(0) or~I(1); that is, they are
  either disjoint or share one corner, with $r$ disjoint from~$x_0$.
  Then we can either compose $r*p$ or $p*r$, and this composite has an
  alternate decomposition which is either of the form $r'*p'$ or
  $p'*r'$ where $r'$ is an empty rectangle distinct from~$r$ and
  $p'$ is a domain of type~$F$ distinct from~$p$. We have the
  following cases:
  \begin{itemize}
  \item if $p*r=r'*p'$ (or $r*p=p'*r'$), then
    $\mu(p)\sign_0(r)+(-1)^m\sign_0(r')\mu(p')=0$;
  \item if $r*p=r'*p'$ or $p*r = p'*r'$, then
    $\sign_0(r)\mu(p)+\sign_0(r')\mu(p')=0$.
  \end{itemize}
\end{lemma}

\begin{proof}
  If $r*p=p'*r'$ and $p$ has complexity $m$,
  then the decomposition $r*D(p)$ can be obtained from $D(p')*r'$
  by $m-1$ elementary moves: we successively commute the rectangles in
  $D(p)$ past the rectangle $r$. The case where $r'*p'=p*r$ is symmetric.

  Otherwise, $r$ shares an edge with some rectangle $r_i$ contained in
  $D(p)$.  With some number $k$ of elementary moves, we can change to
  a composition series where some rectangle~$s$, with the same support
  as~$r$, appears next to~$r_i$.  Then we can perform one elementary
  move to change these two rectangles ($s$ and $r_i$) to $s'$ and
  $r_i'$, respectively, where $s$ has the same support as $r'$; then
  $k$ more elementary moves returns us to the composition series for $r'$
  and $D(p')$, for a total of $2k+1$ elementary moves, which is odd, as
  desired.
\end{proof}

\begin{lemma}
  \label{lemma:CancelSignsType2}
  Suppose that $p\in\pF(\x,\y')$ is a domain of complexity~$m$,
  $r\in\Rect(\y,\z)$ is a Type 2 rectangle, and the corners of $p$ and
  $r$ are distinct; that is, they are in
  case~II(0).  This case matches with case~I(2), so either $p*r$ has
  an alternate decomposition $r'*p'$, in which case
  \begin{equation}
    \label{eq:TypeTwoCommutation}
    \mu(p)\sign_0(r)+(-1)^m\mu(p')\sign_0(r')=0,
  \end{equation}
  or $p*r$ has
  an alternate decomposition as $p'*r'$, where $r'$ is a Type 1
  rectangle, in which case
  \begin{equation}
  \label{eq:TypeTwoAntiCommutation}
  \mu(p')\sign_0(r')+\mu(p)\sign_0(r)=0.
  \end{equation}
\end{lemma}

\begin{proof}
        In the first case, we find it convenient to
        start with the decomposition $r'*p'$,
        and write the standard decomposition
        $$D(p')=r_1'*\cdots*r_{m+1}'.$$
        Recall that $r'$ shares two corner points with~$p'$. There are two
        cases, according to whether these two corner points are upper
        right or lower left corners of $p'$. We consider the upper right
        case.  In this case, the boundary of $r'$ meets the boundary
        of two consecutive odd rectangles, as in Figure~\ref{fig:TypeTwoCase}.
\begin{figure}
\begin{center}
\mbox{\vbox{\epsfbox{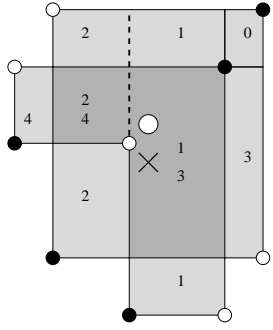}}}
\end{center}
\caption {{\bf $p*r=r'*p'$, where $r$ has Type 2.}
        The rectangle labeled by $0$ is the rectangle $r'$,
        and $p'$ is a polygon with complexity 5, whose standard decomposition
        is indicated by the numbers. This decomposition can be transformed
        into a decomposition $D(p)*D(r)$, where $r$ has Type 2, in 
        an odd number of steps.}
\label{fig:TypeTwoCase}
\end{figure}
        Write the first as
        $r_{2i-1}'$. Starting from $r'*D(p')$, we perform $2i-1$
        elementary moves, to commute $r'$ past $r'_{2i-1}$, to obtain
        a new decomposition
        $r_1*\cdots*r_{2i-1}*r''*r'_{2i}*r'_{2i+1}*\cdots*r'_{m-1}$.  The
        support of the union of the three consecutive rectangles
        $r''*r'_{2i}*r'_{2i+1}$ is a rectangle, which is decomposed so
        that the upper right rectangle $r''$ comes first. 

        We now need some terminology, and then a simple observation.
        Suppose that $s_1*s_2*s_3$ are three rectangles, whose union
        is a rectangle (so that some point $x$ in its interior is a corner of two of the three original rectangles). Suppose also that $s_1$, $s_2$, and $s_3$ are ordered so that the lower left corner of
        the total rectangle is the lower left corner of $s_1$, while
        the upper left corner in the total rectangle is the upper left
        corner of $s_3$. We then say that $s_1*s_2*s_3$ is the
        standard decomposition of a rectangle. (This notion coincides
        with the earlier standard decomposition of a Type 2 rectangle,
        when the central point $x=x_0$.)  Suppose now that $r$ is some
        rectangle which can be post-composed with $s_1*s_2*s_3$, and
        which has two corners inside the support of $s_1*s_2*s_3$. Then after an odd number of elementary moves
        (actually, $3$ or $5$), we can transform $s_1*s_2*s_3*r$ to
        $r'*s_1'*s_2'*s_3'$, so that the supports of $r$ and $r'$
        coincide, and $s_1'*s_2'*s_3'$ is a standard decomposition of
        the rectangle.

Starting from the composition $r''*r'_{2i}*r'_{2i+1}$, we can apply two elementary moves to transform it into the standard decomposition $s_1*s_2*s_3$ of a rectangle. Then, applying the principle in the previous paragraph $m-2i$ times, we can commute $s_1*s_2*s_3$ to the end of the decomposition
$$D'=r_1*\cdots*r_{2i-1}*s_1*s_2*s_3*r'_{2i+1}*\cdots*r'_{m-1},$$  
turning it into the desired decomposition $D(p)*D_0(r)$, where $r$ is a Type 2 rectangle. In all, the number of elementary moves has the same parity as
        $m+1$, verifying Equation~\eqref{eq:TypeTwoCommutation} in the
        case where $r'$ is on the upper right side of the support
        of~$p'$.  The case when $r'$ meets~$p'$ in lower left corners
        of~$p'$ is similar.

        For Equation~\eqref{eq:TypeTwoAntiCommutation}, again there
        are two cases, according to whether $r'$ shares two lower
        right or two upper left corner points of $p'$. Assume they are
        lower right, and write the standard decomposition of~$p'$,
        $D(p')=r'_1*\cdots*r'_{m+1}$. Now, one edge of $r'$ is contained
        in $r'_{2i-1}$, while the other is contained in $r'_{2i+1}$.
        Consider now the decomposition $D(p')*r'$. In $m-2i$
        elementary moves, we commute $r'$ before $r'_{2i-1}$, to
        obtain a new decomposition $$r'_1*\cdots*r'_{2i-2} * r_{2i-1} *
        s_1 * s_2 * s_3 * t_{2i+2} * \cdots * t_{m+2},$$ where
        $s_1*s_2*s_3$ is a decomposition of a rectangle. After two
        elementary moves, we can change it to a standard decomposition
        of the rectangle. Applying the observation about commuting
        rectangles discussed above, we can now commute this
        decomposition of the rectangle to the end of our discussion;
        the number of steps is congruent to $m+1$ modulo two. The new
        decomposition is the decomposition $D(p)*r$, and it was
        obtained from $D(p')*r$ by an odd number of elementary moves.
        Once again, the case when $r'$ meets~$p'$ in upper left
        corners of~$p'$ works similarly.
\end{proof}

\begin{lemma}
\label{lemma:CancelSignsThin}
Let $r\in\EmptyRect(\x,\y)$ and $p\in\pF(\y,\z)$
a domain with complexity $m$. Suppose
that $r$ and $p$ share one corner, and suppose that $x_0$
appears in the interior of the boundary of $r$ (Case I(1\/$'$)).  Then
there is a horizontal or vertical annulus so that the domain~$p_0'$
obtained by adding the annulus to $r * p$ has an alternate
decomposition.
We have the following cases:
\begin{enumerate}
\item
\label{item:CaseGenericI3}
if there is a $\y'\neq\y$ and 
$r'\in\EmptyRect(\x,\y')$ and $p'\in\pF(\y',\z)$ so that $r'*p' = p_0'$
(Case I(3)),
then $\sign_0(r)\mu(p)+\sign_0(r')\mu'(p)=0$.
\item
\label{item:OtherCase}
if there is a $\y'\neq\y$ and 
$p'\in\pF(\x,\y')$ and $r'\in\EmptyRect(\y',\z)$ so that $p'*r' = p_0'$
(Case I(3) if $r'$ has Type 1, Case II(1) if $r'$ has Type 2)
then 
$\sign_0(r)\mu(p)+(-1)^{m}\mu(p')\sign_0(r')=0$.
\end{enumerate}
\end{lemma}

\begin{proof}
Note that in the alternate decomposition, in the case where $r'$
is of Type 1, the rectangles $r'$ and $p'$ meet in three points,
and the composite domain of $r'$ and $p'$ contains an annulus.

Consider Case~\eqref{item:CaseGenericI3}. This can be subdivided
into two subcases: either the annulus is vertical or horizontal.
Suppose the annulus is horizontal.  In this case, 
write $D(p')=r'_1*\cdots*r'_{m+1}$. Consider the decomposition
$r'*D(p')$. After $m$ elementary moves (commuting $r'$ so that it is
next-to-last), we obtain an alternate decomposition, where the last
two rectangles compose to the row through~$O_1$. Cancelling these last
two, and performing $m-1$ more elementary moves (bringing the last
term to the first), we obtain the decomposition $r*D(p)$. The total
change in sign is $(-1)^{2m-1}=-1$, so
$\sign_0(r)\mu(p)+\sign_0(r')\mu(p')=0$ as claimed.

Consider next the case that the annulus is vertical. Write
$D(p')=r'_1*\dots*r'_{m+1}$, and consider the decomposition
$r'*r'_1*\dots*r'_{m+1}$. The rectangles~$r'$ and~$r'_1$ together form
the column through~$O_1$.  Thus, cancelling the first two terms (and
introducing a minus sign by Property~(V)), we obtain a decomposition
$r'_2*\dots*r'_{m+1}$. Now $r'_2=r$, and $r'_3*\dots*r'_{m+1}$ is a
decomposition of $p$. We have changed signs only once, so again
$\sign_0(r)\mu(p)+\sign_0(r')\mu(p')=0$.

Examples of both possibilities for Case~\eqref{item:CaseGenericI3} are
given in Figure~\ref{fig:Case1Sign}.

\begin{figure}
\begin{center}
\input{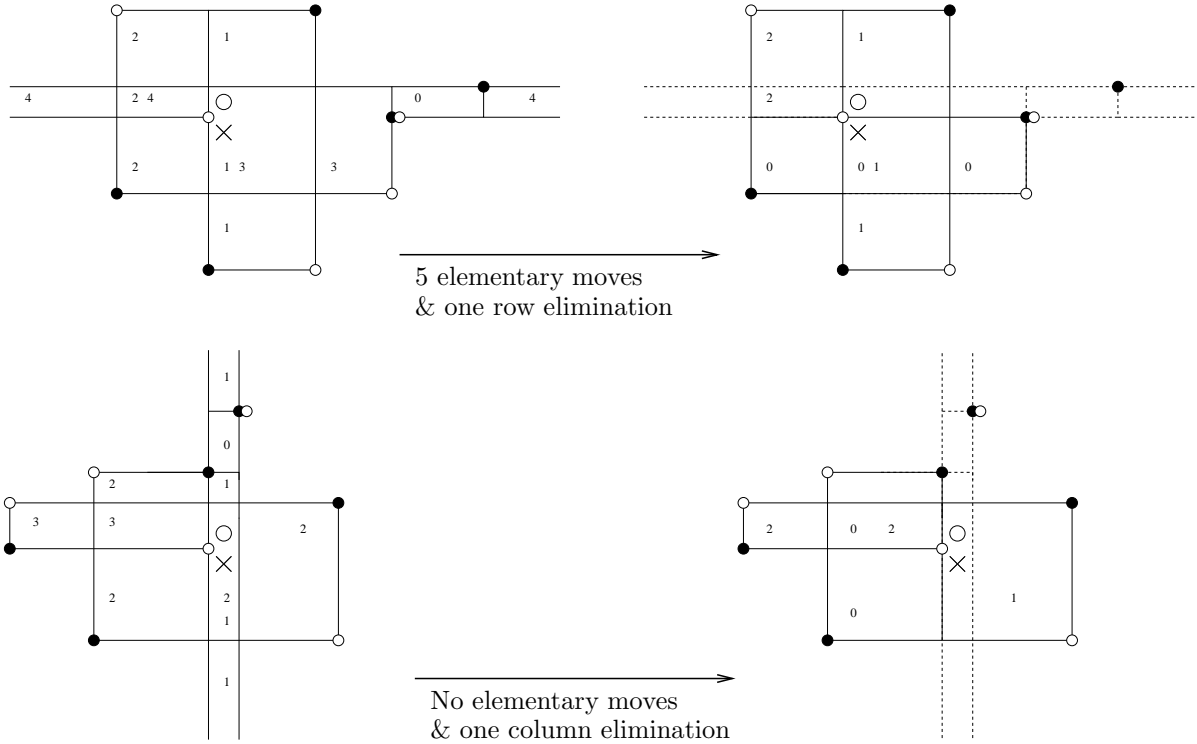}
\end{center}
\caption {{\bf Case~(1) of Lemma~\ref{lemma:CancelSignsThin}.}
\label{fig:Case1Sign}
We have illustrated examples of Case~(1) 
in Lemma~\ref{lemma:CancelSignsThin}. The domains on the right
column correspond to decompositions of the form $r*p$, where
$x_0$ is contained in a boundary of $r$. On the left, we have
corresponding alternate composite domains $r'*p'$.
The domains here are decomposed into ordered rectangles; the integers
in a region give the number of the rectangle the given region belongs to.}
\end{figure}

Consider Case~\eqref{item:OtherCase}, which we divide into subcases:
either $r'$ has Type 1 or Type 2. Consider first the case where $r'$
has Type~1. This again can be divided into two subcases, according to
whether the annulus in the decomposition $r*p$ is horizontal or
vertical. Suppose first that it is horizontal.  Write
$p'=r'_1*\dots*r'_{m+1}$, and consider the decomposition
$D(p')*r'$. Performing an elementary move on the last two rectangles,
and then on the next-to-last two, we obtain a new decomposition
$r_1'*\dots*r_{m-1}'*s_m*s_{m+1}*s_{m+2}$, with the property
that~$s_{m}$ and~$s_{m+1}$
form the row through~$O_1$. Thus, they can be cancelled; performing
$m-1$ elementary moves (commuting $s_{m+2}$ to the beginning of the
decomposition), we obtain the decomposition $r*D(p)$ with total
sign change $(-1)^{m+1}$, verifying the
claim. In the case where the annulus is vertical, write
$D(p')=r_1'*\dots*r_{m+1}'$ and consider the
decomposition $D(p')*r'$. We commute the last rectangle~$r'$ to the second
place in $m$ moves, then cancel the first two rectangles (which, since
they form a vertical annulus, introduces the sign $-1$) to obtain the
alternate decomposition $r*D(p)$ (with total sign change
$(-1)^{m+1}$).  This verifies the stated relation when $r'$ has
Type~1.

\begin{figure}
\begin{center}
\input{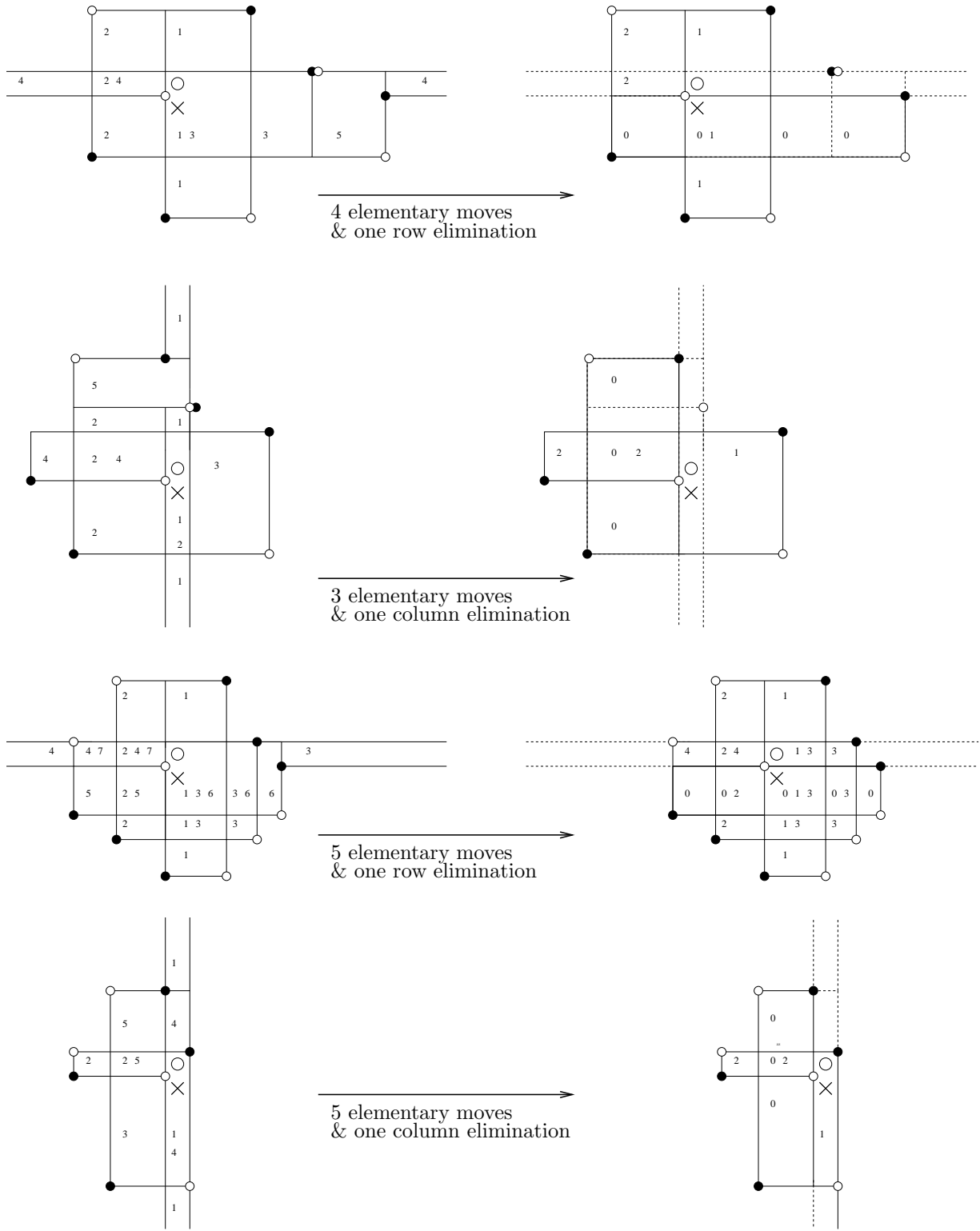}
\end{center}
\caption {{\bf Case~(2) of Lemma~\ref{lemma:CancelSignsThin}.}
\label{fig:Case2Sign}
We have illustrated examples of Case~(2) 
in Lemma~\ref{lemma:CancelSignsThin}. The conventions are
the same is an Figure~\ref{fig:Case1Sign}.}
\end{figure}

We turn to the case where $r'$ has Type 2.  Again, we have two cases,
according to whether the annulus is vertical or horizontal. Assume the
annulus is horizontal. Write the standard decompositions
$D(p')=r_1'*\dots*r'_{m-1}$ and $D_0(r')=r'_{m}*r'_{m+1}*r'_{m+2}$, and
consider the decomposition $D(p')*D_0(r')=r_1'*\dots*r'_{m+2}$.
Performing one elementary move, replacing $r'_m*r'_{m+1}$ by
$s_m*s_{m+1}$, we
obtain a new decomposition in which the rectangles~$r'_{m-1}$
and~$s_m$ form a row, and hence can be cancelled. Finally, in $m-2$
elementary moves (commuting $s_{m+1}$ to the beginning of the
decomposition), we obtain the decomposition $r*D(p)$, with total sign
change $(-1)^{m-1}$.

Consider the final case, where $r'$ has Type~2 and the annulus is vertical.
Again, write the decomposition $D(p')=r_1'*\dots*r'_{m-1}$, and
$D(r')=r'_m*r'_{m+1}*r'_{m+2}$.  Moving $r'_{m+1}$ to the second place
in $m-1$ steps, we obtain a new decomposition whose first two terms
make up a column. Cancel this, at the cost of introducing one more
$-1$. Next, commute the last two rectangles, and then 
move the pair to the first two
spots in an even number of steps. In this way,
we end up with the decomposition $D(r)*D(p)$, with a change
in sign of $(-1)^{m+1}$, as needed.

Examples for all possibilities of Case~\eqref{item:OtherCase} are shown
in Figure~\ref{fig:Case2Sign}.

\end{proof}

We can now define the stabilization map.  In the same way as in
Subsection~\ref{subsec:Stabilization}, define
\begin{align*}
\FL(\x)&=\sum_{\y\in\S}\,\sum_{p\in\pL(\x,\y)}
\mu(p)\cm U_2^{O_2(p)}\cdots U_n^{O_n(p)}\cm\y \\
\FR(\x)&=\sum_{\y\in\S}\,\sum_{p\in\pR(\x,\y)}
\mu(p)\cm U_2^{O_2(p)}\cdots U_n^{O_n(p)}\cm\y.\\
\end{align*}
and put them together to give:
\[
F=\begin{pmatrix}
\FL \\
\FR 
\end{pmatrix} \colon C \longrightarrow C'.
\]

We have the following sign-refinement of Lemma~\ref{lemma:FIsChain}:

\begin{lemma}
  \label{lemma:FIsChainZ}
The map $F\colon C \longrightarrow C'$ preserves
Maslov grading, respects Alexander filtrations, and is a chain map
with coefficients in~$\Z$.
\end{lemma}

\begin{proof}
  Our goal is to show that
  \[
  F \circ \partial_C - \partial_{C'} \circ F = 0.
  \]
  Recall that $\partial_C'$ has three terms: rectangles in~$\Left$,
  rectangles in~$\Right$ (counted with the opposite sign) and the
  differential from~$\Left$ to~$\Right$, multiplication by
  $(U_2-U_1)$.  Again, we can collect the terms in $F \circ \partial_C
  - \partial_{C'} \circ F$ into terms of Types~I(0), I(1), I(1$'$),
  I(2), I(3), II(0), II(1), and~(S). In the proof of
  Lemma~\ref{lemma:FIsChain}, we have seen that these terms can be
  grouped into pairs. We must show that in each pair, the associated
  signs cancel.
  
  Lemma~\ref{lemma:CancelSignsType1} ensures that the terms in
  Case~I(0) and Case~I(1) drop out in cancelling pairs.  The
  interesting case is when we have alternate decompositions $p*r =
  r'*p'$; in this case, if $p$ is of type~$L$, $m$ is odd and the
  differential in~$\Left$ corresponding to~$r$ is taken with the usual
  sign~$\sign_0(r)$, while if $p$ is of type~$R$, $m$ is even and the
  differential in~$\Right$ is taken with sign $-\sign_0(r)$.  In both
  cases the terms cancel.  Similarly,
  Lemma~\ref{lemma:CancelSignsType2} ensures that all terms
  with complexity $m\geq 3$ in Case~I(2) drop out with their
  corresponding terms in~II(0).

  Lemma~\ref{lemma:CancelSignsThin} ensures that all terms in I($1'$)
  cancel with their corresponding terms of types I(3) or II(1),
  leaving possible terms of Type (S). Specifically, a term of
  type I(1$'$) corresponds to a decomposition $p_0=r_1*p_1$, where
  $r_1$ is an empty rectangle and $p_1$ is a term of type $F$ and
  complexity $m$. Adding an annulus to $p_0$ as in the proof of
  Lemma~\ref{lemma:FIsChain}, we obtain a new domain $p_0'$, which
  in turn decomposes as $r_2*p_2$ or $p_2*r_2$ as in case I(3) or II(1).
  In cases where $O_1\not\in r_1$, these
  terms appear in cancelling pairs according to
  Lemma~\ref{lemma:CancelSignsThin}.

  In the cases where $O_1\in r_1$, the decomposition $r_1*p_1$
  contributes once counted with a multiple of~$U_1$
  (as $r_1$, which contains~$O_1$, is thought of as a differential in $C$),
  the decomposition $p_2*r_2$ or $r_2*p_2$ contributes with a multiple
  of $U_2$ (as it contains the row through $X_1$), but there is also
  a contribution coming from the composite domain
  $p_0$, thought of as a domain of type $L$, times $(U_2-U_1)$, the
  differential from~$\Left$ to~$\Right$ within $C'$. Cancellation of the terms
  involving $U_1$ follows from the observation that $r_1*D(p_1)$ differs
  from the standard decomposition of $p_0$ by a sequence of $m-1$
  elementary moves (commuting $r_1$ to the very end), where here $m$
  denotes the complexity of~$p_1$
  (which is necessarily even).  Cancellation of the terms
  involving~$U_2$ follows since they have the opposite sign from the
  terms involving~$U_1$.
  
  It is straightforward to see that the remaining possible $m=2$ 
  terms in~II(2) cancel the remaining possible two $m=1$ terms of type~(S).
\end{proof}

Putting everything together, we have the following:

\medskip
\begin{proof}[Proof of Theorem~\ref{thm:WithSigns}]
  This result now is an immediate consequence of Cromwell's theorem
  and the sign refinements discussed above. Specifically, independence
  of the choice of sign assignment is established in the uniqueness
  statement of Theorem~\ref{thm:SignAssignments}. Commutation
  invariance follows from Proposition~\ref{prop:CommuteZ}.
  Stabilization invariance follows from Lemma~\ref{lemma:FIsChainZ},
  together with a straightforward adaptation of the 
  proof of Proposition~\ref{prop:Stabilization}.
\end{proof}

It is sometimes convenient to consider the chain complex $\CLm(G)$
which is the graded object associated to the Alexander filtration of
$\Cm(G)$. Explicitly, it is the group with the same underlying chain
complex, endowed with a differential as in Equation~\eqref{eq:DefDm}.
It is a formal consequence of Theorem~\ref{thm:WithSigns} that the
homology $\HLm(G)$ of $\CLm(G)$, thought of as a module over
$\Z[U_1,\dots,U_\ell]$, is a link invariant.


\section{More properties}
\label{sec:MoreProperties}

We next give a few of the basic properties of knot and link Floer
homology.  Again, most of these properties are
well-known~\cite{Knots,RasmussenThesis,Links}; but again, we can give
a self-contained derivation here.  Let $\HLa_d(\orL,\s)$ be the part
of $\HLa(\orL)$ with Alexander grading~$\s$ and Maslov grading~$d$.

\begin{prop}
\label{prop:TotalHomology}
The total homology groups of the chain complex $\Cm(G)$ are isomorphic
to the module $\Z[U]$, where all the $U_i$ act as multiplication by
$U$.  The homology groups of $\Ca(G)$ are isomorphic to $\Z$.
\end{prop}

\begin{proof}
  The chain complex $\Cm(G)$ refers to the $\{X_i\}_{i=1}^n$ only
  through its Alexander filtration; in particular, the homology of
  $\Cm(G)$ makes no reference to this placement, and it is unchanged
  by a rearrangement of these decorations (though it does appear to
  depend on the placement of the $\{O_i\}_{i=1}^n$). Now, given any
  grid diagram $G$, we can consider instead the alternate grid diagram
  $H$ gotten by placing $X_i$ in the square immediately under each
  $O_i$. This new diagram clearly represents a suitably stabilized
  diagram for the unknot.  Indeed, after destabilizing sufficiently
  many times, we can reduce to the $2\times 2$-grid diagram $J$ for
  the unknot. A direct calculation in this case gives that
  $H_*(\Cm(J))\cong \Z[U]$ (or $\Field[U]$ with coefficients
  modulo~$2$).

The analogous statement for $\Ca(G)$ follows similarly.
\end{proof}

Proposition~\ref{prop:TotalHomology} allows us to define the invariant $\tau(K)$
for a knot $K$ (see~\cite{4BallGenus,RasmussenThesis}, compare
also~\cite{RasmussenSlice}): If we consider the natural inclusion map
$\iota_m\colon
\Filt_m(\Ca(G))\longrightarrow \Ca(G)$, then $\tau(K)$ is the smallest
integer $m$ for which the map induced on homology by $\iota_m$ is
non-trivial, as a map to $H_*(\Ca(G))\cong\Z$.

The Alexander polynomial of a link remains unchanged under overall
orientation reversal, it is a symmetric polynomial, and it is
invariant under mirror.  These three facts are reflected in
Propositions~\ref{prop:ReverseLink},
\ref{prop:JSymmetry}, and~\ref{prop:Mirror}, respectively.

\begin{prop}
  \label{prop:ReverseLink}
  The filtered quasi-isomorphism type of the complex~$\Cm(G)$ does not
  change if we reverse the orientation of all components of the
  link~$\orL$.
\end{prop}

\begin{proof}
  Consider the diagram~$G'$ obtained by switching the $x$ and $y$
  coordinates, thus flipping $G$ along the diagonal from the bottom
  left to upper right corner.  Switching the $x$ and $y$ coordinates
  also gives a map from the original set of generators~$\S$ to the new
  set of generators~$\S'$ which preserves both degrees and is a chain
  map.  The new diagram~$G'$ is a diagram for $\orL$ with the
  orientation of each component reversed.

  A few more remarks are needed when working over $\Z$, since the
  pre-composition of a sign assignment with reflection through the
  diagonal is not quite a sign assignment, in the sense of
  Definition~\ref{def:SignAssignment}: the roles of rows and columns
  are reversed. However, this is remedied by substituting $-U_i$ in
  place of $U_i$.
\end{proof}

\begin{prop}
  \label{prop:JSymmetry}
  Given ${\mathbf s}=(s_1,\ldots,s_\ell)\in\halfzl$, 
  we have that 
  $$\HLa_{d}(\orL,{\mathbf s})\cong \HLa_{d-2S}(\orL,-{\mathbf s}),$$
  where $S=\sum_{i=1}^\ell s_i$.
\end{prop}

\begin{proof}
  Fix a grid diagram for $\orL$, and let $A^1$ and $M^1$ denote its
  total Alexander filtration and Maslov grading.  (By total Alexander
  filtration, we mean the sum of the components of the Alexander
  multi\hyp filtration.) Switching the roles of
  $\Os$ and $\Xs$, we obtain a grid diagram for $-\orL$. Differentials
  within $\Caa$ are the same for the two diagrams, but the Alexander and
  Maslov gradings are different.  We let $A^2$ and $M^2$ denote the
  Alexander and Maslov gradings of the new diagram. We find it
  convenient to symmetrize $A^i$, defining
  \[
    {\widetilde A}^i(\x)= A^i(\x)+\Big(\frac{n_1-1}{2},\ldots,\frac{n_\ell-1}{2}\Big)
  \]
  for $i=1,2$.
  It is a straightforward calculation from
  Equations~(\ref{eq:MaslovFormula})
  and~\eqref{eq:AlexanderFormulaTwo} that
  \begin{align*}
    M^1-2\sum_{i=1}^\ell \widetilde{A}^1_i &= M^2 \\
    -{\widetilde A}^1&={\widetilde A}^2.
  \end{align*}
  The result now follows from Proposition~\ref{prop:ReverseLink}
  together with Proposition~\ref{prop:HLaa}.
\end{proof}

Indeed, we have the following more general version:

\begin{prop}
  \label{prop:ReverseComps}
        Let $\orL$ be an oriented, $\ell$-component link, and 
        let $\orL'$ be the oriented link obtained from $\orL$ by reversing
        the orientation of its $i\th$ component. Then,
        writing ${\mathbf s}=(s_1,\ldots,s_\ell)$, 
        $$\HLa_{d}(\orL,(s_1,\ldots,s_{\ell}))
        \cong 
        \HLa_{d-2s_i+\ell_i}(\orL',(s_1,\ldots,s_{i-1},-s_i,s_i,\ldots,s_\ell)),$$
        where here $\ell_i$ denotes the total linking number of the
        $i\th$ component
        of $\orL$ with the remaining components.
\end{prop}

\begin{proof}
  From a grid diagram for~$\orL$ we can obtain a grid
  diagram for $\orL'$ by
  switching the roles of $\Os_i$ and $\Xs_i$, i.e., those markings
  which correspond to the $i\th$ component of the link. As in the
  proof of Proposition~\ref{prop:JSymmetry}, the complexes~$\Caa$
  agree, but the Maslov and Alexander functions
  change, as follows. Let ${\widetilde A}^1_j$ and $M^1$ be the
  $j\th$ symmetrized Alexander filtration and Maslov grading for
  $\orL$, respectively, and let ${\widetilde A}^2_j$ and $M^2$ be
  the corresponding
  functions for~$\orL'$.  Let $\widetilde \Os_i = \Os \setminus
  \Os_i$ and $\widetilde \Xs_i = \Xs \setminus \Xs_i$.

  By a direct application of
  Equation~\eqref{eq:AlexanderFormulaTwo},
  \begin{align*}
    {\widetilde A}_j^2(\x)=\begin{cases}
      {\widetilde A}_j^1(\x)  & i\neq j \\
      -{\widetilde A}_i^1(\x)  & i=j
    \end{cases}
  \end{align*}
  and, using Equations~\eqref{eq:MaslovFormula} and~\eqref{eq:AlexanderFormulaTwo},
        \begin{align*}
          M^2(\x)-M^1(\x)
          &=-2\NESW(\x-\Os,\Xs_i-\Os_i) + \NESW(\Xs_i - \Os_i,\Xs_i - \Os_i)\\
          &= -2\widetilde A_i(\x) - \NESW(\Xs - \Os,\Xs_i - \Os_i) +
          \NESW(\Xs_i - \Os_i, \Xs_i - \Os_i)\\
        &= -2{\widetilde A}_i^1(\x)+\NESW({\widetilde \Xs}_i-{\widetilde \Os}_i,\Xs_i-\Os_i).
        \end{align*}
        Moreover, it is straightforward to see that $\NESW({\widetilde \Xs}_i-{\widetilde \Os}_i,\Xs_i-\Os_i)=\ell_i$.
\end{proof}

\begin{prop}
  \label{prop:Mirror}
  Let $\oL$ be a link, and let $r(\oL)$ denote its mirror.
  In this case, we have an identification
  $$\HLa_d(\orL,{\mathbf s})\cong \HLa^{2S-d}(\Ca(r(\orL),{\mathbf s}));$$
  note the right-hand-side denotes cohomology.
\end{prop}

\begin{proof}
  Rotating the grid diagram $G$ ninety degrees to get a new diagram
  $G'$ corresponds to passing from the knot to its mirror. There is an
  induced map $\phi$ from $\S(G)$ to $\S(G')$.  Letting ${\widetilde
    A}$, ${\widetilde M}$ and ${\widetilde A}'$, ${\widetilde M}'$ be
  the Alexander and Maslov gradings of $G$ and $G'$ respectively, it
  is clear that ${\widetilde A}(\x)=-{\widetilde A'}(\phi(\x))$,
  ${\widetilde M}(\x)=-{\widetilde M'}(\phi(\x))$.  Indeed, if we
  think of $\phi$ as taking $\x$ to $\x^*$, the dual basis element of
  $\Caa(H)$ which is one on $\x\in\S(H)$ and zero on all other
  $\y\in\S(H)$, then $\phi$ induces an isomorphism of chain complexes.
  The shift in absolute grading now follows from
  Proposition~\ref{prop:JSymmetry}.  
\end{proof}


\section{Relation to the Alexander polynomial}
\label{sec:Alexander}

In this section we will show that the Euler characteristic of the
multi-graded complex~$\CLa$ with respect to the Maslov grading is the
Alexander polynomial.  Precisely, fix a grid diagram~$G$ of a
link~$\orL$ with $\ell$~components. Given $\mathbf{s} \in \halfzl$,
let $t = (t_1,\ldots,t_\ell)$ be a
collection of variables, and for $\s = (s_1,\ldots,s_\ell)$ an element
of $\halfzl$, define
\[
t^\s = t_1^{s_1}\cdots t_\ell^{s_\ell}.
\]
For any multi-graded groups~$C_i(\s)$ with Maslov grading~$i$ and Alexander
grading~$s$, define
$$\chi(C;t) = \sum_{i,\s} (-1)^it^\s\rank(C_i(\s)).$$
\begin{theorem}
  \label{thm:alexander}
  For any link~$\orL$, the Euler characteristic of $\HLa$ determines the
  multi-variable Alexander polynomial up to sign.  Precisely,
  \[
  \chi(\HLa(\orL)) =
  \begin{cases}
     \pm\prod_{i=1}^\ell(t_i^{1/2}-t_i^{-1/2}) \Delta_A(\orL;t)&\ell > 1\\
     \pm\Delta_A(\orL;t)& \ell = 1
  \end{cases}
  \]
  where $\Delta_A(\orL;t)$ is the multi-variable Alexander polynomial,
  normalized so that it is symmetric up to sign under the involution
  of sending all $t_i$ to their inverses.
\end{theorem}

We will prove this by taking the Euler characteristic of the alternate
complex~$\Caa(\orL)$.  The Maslov grading of a generator~$\x\in\S$ is, up
to an overall sign, the sign of~$\x$ considered as a permutation.  The
Alexander grading is, up to an overall shift, minus the sum of the
winding numbers around the generators.  Summing over all
generators~$\x$, we get a ``minesweeper determinant'' as illustrated
below:
\[
\left\vert\mfigb{draws/gridlink.1}\right\vert
\]
This turns out to give the Alexander polynomial, up to a sign, powers
of the variables~$t_i$, and factors of $(1-t_i)$.

More formally, given a grid diagram~$G$ of size $n$, define an
$n\times n$ matrix~$M(G)$ by
\[
M(G)_{ij} = t^{a(i,j)}
\]
where $a(i,j)\in\halfzl$ is the vector whose $k\th$ component is the
minus the winding number of the $k\th$ component of the link around the point
$(i,j)$.  (Here we use the convention that the links runs between the
$O$'s and the $X$'s, which have half-integer coordinates, so this
winding number is well\hyp defined.)
Then we have:
\begin{prop}
  \label{prop:alexander-det}
  For any grid diagram~$G$ of a link~$\orL$ with
  $\ell$~components, let $n_i$ be the number of vertical segments
  corresponding to component~$i$.
  Then
  \[
  \det M(g) = 
  \begin{cases}
    \pm t^{k} (1-t)^{n-1}\Delta(\orL;t)&\ell = 1\\
    \pm \biggl(\prod_{i=1}^\ell t_i^{k_i}(1-t_i)^{n_i}
      \biggr)\Delta(\orL;t) & \ell > 1
  \end{cases}
  \]
  for some integers~$k_i$.  In the case $\ell=1$, we write $t,n,k$ for
  $t_1,n_1,k_1$ for convenience.
\end{prop}

This proposition implies Theorem~\ref{thm:alexander}:

\begin{proof}[Proof of Theorem~\ref{thm:alexander}]
  It follows from Proposition~\ref{prop:HLaa} that
  \[
  \Bigl(\prod_i(1-t_i)^{n_i-1}\Bigr)\chi(\HLa(\orL)) = \chi(\HLaa(G))
    = \chi(\Caa(G)).
  \]
  The theorem follows by Proposition~\ref{prop:alexander-det} up to an
  overall sign and powers of the $t_i$.  The powers of $t_i$ are
  determined by Proposition~\ref{prop:JSymmetry} and the chosen
  normalization of the Alexander polynomial.
\end{proof}

\begin{proof}[Proof of Proposition~\ref{prop:alexander-det}]
We will use Fox's
free differential calculus~\cite{Fox53:FreeDiffCalc} with respect to a
presentation associated to a grid diagram for a link.  This presentation was
described by Neuwirth~\cite{Neuwirth84:StarProjections}, who also
proved that it is actually a presentation of the knot group.

The presentation, as shown in Figure~\ref{fig:gener-relat-gridlink},
has one generator for each vertical segment of the link, starting
from the basepoint (outside of the page at the position of the
reader), coming down to the left of a vertical segment, going under
the segment, and coming back out of the page.  There are $n$
generators, one for each vertical segment.  There are $n-1$ relations,
one for each position between horizontal segments.  The path of the
relation comes down to the left of the diagram, runs across the
diagram at a constant horizontal level, and comes back up out of the
page.  On the one hand this loop is contractible (pull it
beneath the diagram); on the other hand it is equal to the product of
the generators corresponding to the vertical segments that we cross.

\begin{figure}
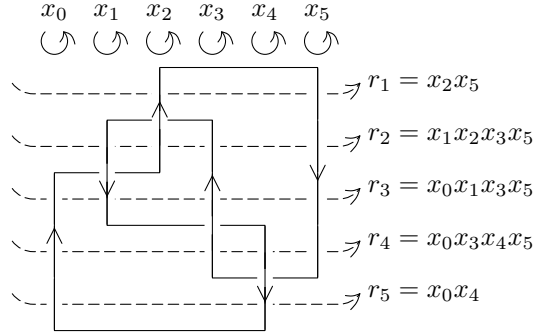

  \[
  \mfig{draws/gridlink.2}
  \]
  \caption{\textbf{Generators and relations from a gridlink presentation.}}
  \label{fig:gener-relat-gridlink}
\end{figure}

In the example in Figure~\ref{fig:gener-relat-gridlink}, there are six
generators, $x_0$ through $x_5$, and five relations:
\begin{align*}
r_1 &= x_2 x_5\\
r_2 &= x_1 x_2 x_3 x_5\\
r_3 &= x_0 x_1 x_3 x_5\\
r_4 &= x_0 x_3 x_4 x_5\\
r_5 &= x_0 x_4.
\end{align*}
The Fox derivative matrix is particularly easy to compute in this
case, since the
generators appear only positively in the relations and at most once in
each relation.  The entries are the initial portions of the relations,
and they appear in the positions given by the vertical strands.  In
our running example, we get
\[
\Bigl(\frac{\partial r_i}{\partial x_j}\Bigr)_{ij} =
\begin{pmatrix}
      0 &    0 &    1 &     0 &     0 &   x_2\\
      0 &    1 &  x_1 & x_1x_2&     0 & x_1x_2x_3\\
      1 &  x_0 &    0 & x_0x_1&     0 & x_0x_1x_3\\
      1 &    0 &    0 &   x_0 & x_0x_3& x_0x_3x_4\\
      1 &    0 &    0 &     0 &   x_0 &     0\\
\end{pmatrix}
\]

To find the Alexander polynomial, map this matrix to the
abelianization of the knot group, mapping each $x_j$ to $t_i^{\pm 1}$,
depending on which component the vertical segment belongs to and
whether the corresponding vertical strand is oriented upwards or
downwards.  For knots, the Alexander polynomial is the determinant of
a maximal minor of the resulting $n$ by $n-1$ matrix, up to a factor
of $\pm t^k$.  For links, the multi-variable Alexander polynomial is,
up to the same factor, the determinant of a maximal minor divided by
$(1-t_i)$ for each component~$i$ that is not the component of the
deleted column.

In the example, we get
\[
\Delta(t) = \pm t^k
\begin{vmatrix}
      0 &    0 &    1 &     0 &     0\\
      0 &    1 &t^{-1}&     1 &     0\\
      1 &    t &    0 &     1 &     0\\
      1 &    0 &    0 &     t &   t^2\\
      1 &    0 &    0 &     0 &     t\\
\end{vmatrix}
\]

Now let us turn to computing $\det (M)$ where $M$ is the minesweeper
matrix defined earlier.  Subtract each column from the next
one.  The winding numbers change by at most one when we move from one
square to a neighbor.  Therefore in every column but the first we have
zero entries where the vertical segment does not intervene, and where
a vertical segment of component~$i$ does intervene every entry is
divisible by $t_i^{\pm 1}-1$.
Thus, for each column but the first, we can factor out $t_i-1$ if the
column is oriented upwards or $t_i^{-1}-1$ if it is oriented downwards.
Furthermore, after this operation the last row contains only a single
non-zero entry, $1$~in the first column, so we can delete the first
column and last row without changing the determinant (up to sign).

In the example, we get
\begin{align*}
\begin{vmatrix}
  1 & 1 & 1 & t & t & t \\
  1 & 1&t^{-1}&1& t & t \\
  1 & t & 1 & 1 & t & t\\
  1 & t & t & t &t^2& t\\
  1 & t & t & t & t & 1\\
  1 & 1 & 1 & 1 & 1 & 1
\end{vmatrix}
&=
\begin{vmatrix}
   1 &   0 &  0     & t-1      &  0   & 0\\
   1 &   0 &t^{-1}-1& 1-t^{-1} & t-1  & 0\\
   1 & t-1 & 1-t    & 0        & t-1  & 0\\
   1 & t-1 & 0      & 0        & t^2-t& t-t^2\\
   1 & t-1 & 0      & 0        & 0    & 1-t\\
   1 &   0 & 0      & 0        & 0    & 0
\end{vmatrix}\\
&=
\pm (t-1)^3 (t^{-1}-1)^2
\begin{vmatrix}
  0 & 0 & 1 & 0 & 0 \\
  0 & 1 & t^{-1} & 1 & 0\\
  1 & t & 0 & 1 & 0\\
  1 & 0 & 0 & t & t^2\\
  1 & 0 & 0 & 0 & t
\end{vmatrix}.
\end{align*}

Up to the expected factors of $1-t_i$ and $\pm t^k$, this is the
same determinant we got from the Fox derivative calculations.
\end{proof}


\bibliographystyle{custom}
\bibliography{biblio}

\end{document}